\documentclass[times,final]{elsarticle}
\pdfoutput=1

\usepackage{jcomp}
\usepackage{framed,multirow}

\usepackage{amssymb}
\usepackage{latexsym}
\usepackage{bm}
 \usepackage{indentfirst}
\usepackage{amsmath,amsthm,amssymb}

\usepackage{graphicx}
\usepackage{subcaption}
\usepackage{float}
\usepackage[framemethod=tikz]{mdframed}
\usepackage{epstopdf}

\usepackage{import}

\newcommand{\velx}{u}

\renewcommand*{\Pr}{Pr}

\newif\ifturb

\newif\ifsutherland

\newcommand{\mask}{\chi}
\newcommand{\IBMparam}{\eta}

\newcommand{\advcoef}{c}

\usepackage{url}
\usepackage{xcolor}
\definecolor{newcolor}{rgb}{.8,.349,.1}

\journal{Journal of Computational Physics}

\begin{document}

\verso{Jiaqing Kou \textit{etal}}

\begin{frontmatter}

\title{A combined volume penalization / selective frequency damping approach for immersed boundary methods applied to high-order schemes}

\author[1,2]{Jiaqing \snm{Kou}\corref{cor1}}
\cortext[cor1]{Corresponding author}
\ead{jiaqingkou@gmail.com}
\author[1,3]{Esteban \snm{Ferrer}}

\address[1]{ETSIAE-UPM-School of Aeronautics, Universidad Politécnica de Madrid, Plaza Cardenal Cisneros 3, E-28040 Madrid, Spain}
\address[2]{NUMECA International S.A., Chaussee de la Hulpe 187, Brussels, B-1170, Belgium}
\address[3]{Center for Computational Simulation, Universidad Politécnica de Madrid, Campus de Montegancedo, Boadilla del Monte, 28660 Madrid, Spain}


\begin{abstract}
There has been an increasing interest in developing efficient immersed boundary method (IBM) based on Cartesian grids, recently in the context of high-order methods. IBM based on volume penalization is a robust and easy to implement method to avoid body-fitted meshes and has been recently adapted to high order discretisations \citep{kou2021IBMFR2}. This work proposes an improvement over the classic penalty formulation for flux reconstruction high order solvers. We include a selective frequency damping (SFD) approach \citep{aakervik2006steady} acting only inside solid body defined through the immersed boundary masking, to damp spurious oscillations. An encapsulated formulation for the SFD method is implemented, which can be used as a wrapper around an existing time-stepping code. The numerical properties have been studied through eigensolution analysis based on the advection equation. These studies not only show the advantages of using the SFD method as an alternative of the traditional volume penalization, but also show the favorable properties of combining both approaches. This new approach is then applied to the Navier-Stokes equation to simulate steady flow past an airfoil and unsteady flow past a circular cylinder. The advantages of the SFD method in providing improved accuracy are reported.
\end{abstract}


\end{frontmatter}


\section{Introduction}
Mesh generation is an essential procedure in computational fluid dynamics (CFD) to simulate flow over complex geometries. To accurately impose boundary conditions, the mesh usually conforms to the geometry \citep{steger1980generation,pruett2003temporally} (known as body-fitted mesh). Generating body-fitted meshes for complex flows is known to be a difficult and time consuming task \citep{bern2000mesh}. As an alternative, immersed boundary method (IBM) \citep{peskin1972IBM} allows the use of simple meshes such as the Cartesian grids to handle complex and moving geometries more easily. Since Cartesian grids are considered, IBM are well suited to fast solvers and hence are often efficient \citep{ye2020openfsi}.

IBM originates from the early work of Peskin \citep{peskin1972IBM}, where singular forces centered at the boundary are modeled as Delta functions that distribute the external force to surrounding grid to mimic the effect of boundary. Since then, different IBM approaches have been proposed, as reviewed by Mittal and Iaccarino \citep{mittal2005immersed}, Sotiropoulos and Yang \citep{sotiropoulos2014immersed}, and Griffith and Patankar \citep{griffith2020immersed}. Currently, IBM approaches mainly include cut-cell \citep{ye1999accurate,udaykumar2001sharp},  direct forcing \citep{fadlun2000combined,luo2012numerical,tian2014fluid}, ghost-cell \citep{majumdar2001rans,tseng2003ghost}, volume penalization \citep{angot1999penalization,carbou2003boundary,brown2014CBVP,abgrall2014IBM}, etc. Among all these approaches, volume penalization has shown advantages in terms of robustness, ease of implementation and rigorous theoretical background of convergence properties~\cite{arquis1984conditions,angot1999penalization}. It models the solid obstacle as a porous medium with very small permeability \citep{kadoch2012volume,schneider2015immersed}. This is achieved by introducing an external penalization source term to the governing equations to impose the boundary conditions. So far, volume penalization has been applied to different problems, such as flapping wings~\cite{kolomenskiy2009fourier}, two-phase flow~\cite{horgue2014penalization}, aeroacoustics \citep{komatsu2016direct}, fluid-structure interaction~\cite{engels2015numerical} and thermal flows~\cite{cui2018coupled}.

In the last decade, high-order methods for computational fluid dynamics (CFD) \citep{wang2013high} have gained attentions due to their higher-accuracy and lower-dissipation properties when compared to low-order numerical techniques at comparable computational cost. These methods include discontinuous Galerkin (DG) \citep{hesthaven2007nodal,karniadakis2013spectral}, spectral difference (SD) \citep{kopriva1996conservative,liu2006SD}, flux reconstruction (FR) \citep{huynh2007FR}, or correction procedure via reconstruction (CPR) \citep{wang2009unifying}. Applying high-order methods on unstructured grids remains a challenge due to the difficulty of generating meshes for complex geometries. This makes the development of IBM under high-order frameworks an appealing alternative. One of the earliest works considering IBM and high order was reported in ~\cite{lew2008,lew2011}, where different possibilities of combining DG with IBM were discussed, and the advantages of using the DG method over standard finite-differences for a $2$D Poisson problem were shown. Fidkowski and Darmofal~\cite{fidkowski2007triangular} were the first to report the use of cut-cell method to solve steady compressible flows over two-dimensional geometries based on DG. The high-order cut-cell approach was also studied in ~\cite{muller2017high,geisenhofer2019discontinuous,cutCellDGAcoustics}. While the cut-cell approaches show superior convergence property on static grids, the extension to moving grids is not straightforward. Therefore, the IBM approaches based on external forcing are worth investigating in the high-order framework, where volume penalization has potential due to the aforementioned advantages. Recently, volume penalization has been applied to the high-order FR scheme to simulate flows past both static and moving bodies \citep{kou2021IBMFR1,kou2021IBMFR2}. However, the performance of volume penalization in the framework of high-order methods remains to be explored and improved. 

Volume penalization (and other IBM approaches based on external forcings) can be improved from the physical intuition that the solid region is modeled as a porous medium which introduces localised damping to recreate the solid \citep{kolomenskiy2009fourier,kadoch2012volume}. When the permeability (or penelization parameter) tends to zero, i.e., $\IBMparam \rightarrow 0$, the obstacle becomes totally solid and the flow solution inside the solid is completely damped. This perspective has been discussed and followed in different works. An example is shown in \citep{kadoch2012volume}, where the physical diffusivity is replaced by a much smaller penalization parameter inside the solid, to control the damping and model the homogeneous Neumann boundary conditions. In addition, in the original IBM approach by Peskin \citep{peskin1972IBM}, the immersed boundary is connected by springs with the damping coefficient based on the material property \citep{kim2007penalty}. A feedback scheme of original IBM approach is proposed by Goldstein et al. \citep{goldstein1993modeling}, where the forcing is modeled as a second-order damped oscillator in order to model the no-slip condition at the surface. This feedback approach was extended to simulate the flexible filaments in a uniform flow by Huang et al. \citep{huang2007simulation}. It is then improved in \citep{margnat2009behaviour} by introducing a new pair of parameters which provides better insight in such a feedback control system driving the velocity to zero. Our recent work applied eigensolution analysis to volume penalization based on FR, indicating that volume penalization imposes the boundary condition by adding additional dissipation inside the solid \citep{kou2021VonNeumann}. Since high-order methods provide low numerical dissipation errors, it becomes critical that IBM impose the correct flow field inside the solid body, through the volume penalization term. As will be shown in Section \ref{sec:naca}, if insufficient penalization is imposed, the solution exhibits oscillations inside the solid, thus leading to nonphysical oscillations in the solid region (and outside the solid). To correct this unphysical behaviour, we propose a new approach to increase damping inside the solid. Selective frequency damping (SFD) \citep{aakervik2006steady} is introduced inside the solid to damp spurious oscillations, which in turn correct the unphysical oscillations outside the solid. The SFD method is originated from capturing the steady solutions of the Navier-Stokes equations in globally unstable configurations \citep{theofilis2011global}. It introduces a temporally low-pass filter in the dynamic system and obtains the filtered steady solution by damping the unstable (temporal) frequencies. We use the low-pass filter given by the SFD method, to damp the oscillations of penalized solutions inside the solid region and avoid the undesired solution in the flow.

In this work, the SFD method is used as an alternative and a complementary approach to volume penalization. The FR scheme, which is a unified framework that can recover a particular SD and nodal DG schemes \citep{huynh2007FR,vincent2011ESFR}, is selected as a typical high-order method. To implement the method efficiently, an operator splitting approach is considered for both the penalization source term \citep{strang1968construction} and the SFD subsystem \citep{jordi2014encapsulated}. Through numerical experiments and eigensolution analyses based on a linear advection equation, we can find guidelines for the numerical parameters of the methods (e.g., volume penalization damping, filter damping and width), and show the advantages of combining the SFD method and volume penalization. To show the efficacy of the proposed approach in real flow simulations, steady and unsteady flow based on Navier-Stokes equations are considered, including flows past a NACA0012 airfoil (steady case) and a circular cylinder (unsteady case).

The remainder of this paper is organized as follows. In Section 2, all methodologies, including the FR method for space discretisation, volume penalization, selective frequency damping and the implementation of the proposed method, are introduced. Next, Section 3 introduces the numerical experiments, including the linear advection equation and flow past cylinder and airfoil based on the Navier-Stokes equations. The numerical behaviors of the proposed approach are investigated based on the eigensolution analyses. Finally, conclusions are drawn in Section 4.

\section{Methodology}
\subsection{Governing equations}
We will use a one-dimensional advection equation to study the properties of the methods (and derive user guidelines) through eigensolution analyses to then apply the method to the compressible Navier-Stokes equations.
\subsubsection{One-dimensional advection}
The method will be first applied to the one-dimensional advection equation defined in space $x$ and time $t$: 
\begin{equation}
    \frac{\partial \velx}{\partial t} + \frac{\partial f}{\partial x} = 0,
\end{equation}
where $\velx(x,t)$ is the transported unknown and $f = \advcoef \velx$ is the flux with advection speed $\advcoef$. Considering the periodic boundary condition with initial condition $u(x,0) = \text{exp}(ikx)$, the analytical solution is $ u(x,t) = \text{exp} [i(kx -\omega t)]$, where $k$ is the wavenumber and $\omega = \advcoef k$ is the (angular) frequency. $i = \sqrt{-1}$ is the imaginary unit. Details on the discretisation is given in \ref{append2}. 

\subsubsection{Navier-Stokes equations}
To test the proposed method in more realistic flow simulations, the governing equations for a compressible viscous fluid are solved:
\begin{equation}
    \frac{\partial \bm{U}}{\partial t} + \mathbf{\nabla} \cdot \bm{F} =  \frac{\partial \bm{U}}{\partial t} + \frac{\partial \bm{F}_x}{\partial x} + \frac{\partial \bm{F}_y}{\partial y} + \frac{\partial \bm{F}_z}{\partial z} = 0\,,
\end{equation}
where $\bm{U}$ denotes the vector of conserved variables $\bm{U} = (\rho , \rho u, \rho v, \rho w, E)^T$. $\rho$ is the density, $u$, $v$,  and $w$ are the velocity components and $E$ is the total energy. The flux vectors are $\bm{F}_{x}$, $\bm{F}_{y}$, $\bm{F}_{z}$ containing the inviscid and viscous fluxes. Details on the governing equations and the FR discretisation for general conservation laws are given in \ref{append1} and in \ref{append3}. 

\subsection{Volume penalization}
\label{sec:vp}
Volume penalization introduces penalization source terms inside the solid region. Firstly, a mask function $\mask(\boldsymbol{x},t)$ is defined, which distinguishes between the fluid region $\Omega _{f}$ and solid region $\Omega _{s}$:
\begin{equation}
\mask  (\boldsymbol{x},t) = \left\{\begin{matrix}1, \, \, \text{if}\, \, \boldsymbol{x} \in \Omega _{s}
\\0,\, \, \text{otherwise}
\end{matrix}\right. .
\end{equation}

The main parameter for volume penalization is the penalization parameter $\IBMparam$, which is used to determine the magnitude of penalization. It has been proved that as $\eta \rightarrow 0$ the solution will converge to the penalized solution \citep{angot1999penalization,carbou2003boundary}. In the present work, we consider Dirichlet boundary conditions (in particular no-slip walls), which lead to the following penalized formulation of the advection equation:
\begin{equation}
    \frac{\partial \velx}{\partial t} + \advcoef \frac{\partial \velx}{\partial x} + \frac{\mask}{\IBMparam}(\velx - \velx_s) = 0,
    \label{eq:adv}
\end{equation}
where $\velx_s$ refers to the penalized boundary condition in the solid, which is fixed to $0$ (non-moving solids) in this study to approximate the no-slip wall boundary condition. The penalization parameter $\IBMparam$ drives the solution to $\velx_s$ as $\IBMparam \rightarrow 0$.

The Navier-Stokes equations with volume penalization for Dirichlet conditions can be written as:
\begin{equation}\label{eq:RHS}
\frac{\partial \bm{U}}{\partial t} = \boldsymbol{RHS} (\bm{U}) + \bm{S} (\bm{U})
\end{equation}
where $\bm{S}$ and $\boldsymbol{RHS}$ refer to the IBM forcing term and the right hand side term of the Navier-Stokes equations $\boldsymbol{RHS} = -\mathbf{\nabla} \cdot \bm{F} (\bm{Q},\bm{U}) $. The IBM source term is written as:
\begin{equation}
\bm{S} (\bm{U}) = \frac{\chi}{\eta} \times \begin{pmatrix}
0\\ 
\rho u_{s}-\rho u\\ 
\rho v_{s}-\rho v\\ 
\rho w_{s}-\rho w\\ 
\frac{\rho}{2}(u_{s}^2+v_{s}^2+w_{s}^2)-\frac{\rho}{2}(u^2+v^2+w^2)
\end{pmatrix}
\label{IBMsource}
\end{equation}
The penalization terms were proposed in \cite{abgrall2014IBM} for compressible Navier-Stokes equations. It should be noted that this formulation also allows nonuniform Dirichlet conditions with additional source term for temperature \cite{abgrall2014IBM}. The penalized velocity inside the solid body is represented by $\mathbf{u}_s = (u_{s},v_{s},w_{s})^T$. Generally, the penalization parameter $\eta$ should be sufficiently small to ensure good accuracy. However, the stiffness of the equation increases as $\eta$ decreases, therefore it should be limited as well. When no-slip boundary condition is considered, we have $\mathbf{u}_s = (0, 0, 0)^T$.

The main advantage of volume penalization over other IBM approaches is that this method is based on rigorous proofs of the convergence \cite{angot1999penalization,carbou2003boundary}. Therefore, the numerical error introduced from the penalization term can be controlled a-priori \cite{brown2014CBVP}. The theory is briefly detailed here, while more details can be found in \citep{engels2015numerical,schneider2015immersed}. The error of the numerical solution of the penalized problem corresponding to the original problem includes two parts \cite{engels2015numerical}, the penalization (modeling) error and the discretisation error:

\begin{equation}
\left \|  u^{exact}  - u_{\eta} ^ {N} \right \| \leq \left \|  u^{exact}  - u_{\eta} \right \| + \left \|  u_{\eta} - u_{\eta} ^ {N} \right \|,
\label{eq:error}
\end{equation}
where $u^{exact} $ is the exact analytical solution of the governing equations, $u_{\eta} $ and $u_{\eta} ^ {N}$ are the exact and numerical solutions of the penalized problem. $\left \|  \cdot \right \|$ is the norm used for quantifying the error. The first part (penalization error) depends on the penalization parameter \cite{schneider2015immersed}, which has a decay rate of $\mathcal{O}(\sqrt{\eta})$ for Dirichlet boundary conditions \cite{angot1999penalization,carbou2003boundary} and $\mathcal{O}(\eta)$ for Neumann boundary conditions \cite{kolomenskiy2015analysis}. In addition, the second part of error (discretisation error) refers to the error between the exact solution and the numerical solution of the penalized equations. With consistent discretisation and a stable numerical scheme, the discretisation error converges with mesh refinement. However, the rate of convergence is not only determined by the numerical scheme, but also limited by the regularity (continuity) of the solution \citep{kadoch2012volume,schneider2015immersed,kou2021IBMFR2}. 

From the error estimate in Equation \ref{eq:error}, it is suggested to use a very small penalization parameter $\eta$ to minimize the penalization error. However, small $\eta$ will lead to very stiff source term, thus requiring a smaller time step. When the penalization source term is treated explicitly, from the linear stability analysis \cite{kolomenskiy2009fourier}, the time step for explicit time integration must be smaller than the penalization parameter $ \Delta t < \eta $. Also in \cite{jause2012simulation}, the intimate coupling between $\Delta t$ and $\eta$ is reported, indicating the error will saturate when $ \Delta t \approx \eta $. These explains why it is usually suggested to use $ \Delta t = \eta $ in the volume penalization method. Another interpretation is that the penalty term acts as a strong damping term with order $\eta$ on the velocity, which has to be resolved by the time discretisation scheme \cite{engels2015numerical}. From a practical point of view, one can select the maximum $\Delta t$ which allows $ \Delta t = \eta $ to maintain both computational efficiency and accuracy. Note that in high-order methods, the time step $\Delta t$ is determined by the Courant-Friedrichs-Levy (CFL) condition, which scales with the polynomial order  \cite{cockburn1989tvb,hesthaven1999stable}. Therefore a practical guideline is to fix $ \Delta t = \eta $ first and determine $\Delta t$ according to stability criterion.

\subsection{Selective frequency damping}
The SFD method was proposed by Åkervik et al. \citep{aakervik2006steady}, as an alternative to classic Newton's method to compute steady solutions (unstable fixed points) of the Navier-Stokes equations in globally unstable states. This method has been widely used to obtain the unstable base flow for global instability analysis of different flow configurations \citep{theofilis2011global}. Considering any nonlinear dynamic system with appropriate boundary and initial conditions, like the Navier-Stokes equations, the following system can written
\begin{equation}
    \dot{\boldsymbol{q}} = \boldsymbol{f}(\boldsymbol{q}),
\end{equation}
where $\boldsymbol{f}$ is the nonlinear operator applied to a state variable vector $\boldsymbol{q}$. The steady-state $\boldsymbol{q}_s$ is obtained when $\dot{\boldsymbol{q}}_s = \boldsymbol{f}(\boldsymbol{q}_s) = \boldsymbol{0}$, where $\boldsymbol{q}_s$ can be the base flow around which the system is linearized. SFD is a method to obtain this base flow when the linearized system has unstable eigenvalues. From control theory, the basic idea of SFD is to introduce a proportional feedback control \citep{kim2007linear} as a forcing term to drive the solution towards the target solution $\boldsymbol{q}_s$ with the control coefficient $\chi_f$
\begin{equation}
    \dot{\boldsymbol{q}} = \boldsymbol{f}(\boldsymbol{q}) - \chi_{f} (\boldsymbol{q} - \boldsymbol{q}_s).
\end{equation}

This forcing is a linear reaction term proportional to the high-frequency content of the flow \citep{casacuberta2018effectivity}. However, since in practice the target base flow $\boldsymbol{q}_s$ is not known a-priori, a low-pass filtered version of $\boldsymbol{q}$ (defined as $\bar{\boldsymbol{q}}$) is used as the target solution \citep{aakervik2006steady}. The first order low-pass time-domain filter used for the SFD method leads to the following transfer function
\begin{equation}
    \frac{\bar{\boldsymbol{q}}}{\boldsymbol{q}} = 1 \times \frac{1}{1 + i \omega_f \Delta},
\end{equation}
where $\omega_f$ is the circular frequency and $i = \sqrt{-1}$. $\bar{\boldsymbol{q}}$ is the temporally filtered quantity and $ \Delta $ is the filter width. The cutoff
frequency of the filter is given by $\omega_c = 1 / \Delta$. This low-pass filter is able to damp the most unstable frequencies, thus avoiding the corresponding instability to reach the expected unstable steady state. This method is adapted from the temporal filtered model developed for large eddy simulation \citep{pruett2003temporally,pruett2006temporal}. From an implementation point of view, to avoid high memory requirement of storing temporal evolution of flow quantities, the differential form is used instead
\begin{equation}
    \dot{\bar{\boldsymbol{q}}} = \frac{\boldsymbol{q}-\bar{\boldsymbol{q}}}{\Delta}.
\end{equation}

In summary, SFD solves the following system of equation
\begin{equation}
    \left\{\begin{matrix}
    \dot{\boldsymbol{q}} = \boldsymbol{f}(\boldsymbol{q}) - \chi_{f} (\boldsymbol{q} - \bar{\boldsymbol{q}})
\\  \dot{\bar{\boldsymbol{q}}} = (\boldsymbol{q} - \bar{\boldsymbol{q}}) / \Delta
\end{matrix}\right. .
\end{equation}

The filtered solution $\bar{\boldsymbol{q}}$ is time varying and this approach converges when the steady state $\boldsymbol{q} = \bar{\boldsymbol{q}}$ is reached. This formulation indicates that the performance of SFD depends on two parameters, including the control coefficient $\chi_{f}$ and the filter width $\Delta$. These parameters are key to the stability and the convergence rate of this method \citep{aakervik2006steady,jordi2014encapsulated,casacuberta2018effectivity}. As noted in \citep{aakervik2006steady}, the filter cutoff frequency $\omega_c$ is related to the frequency of the relevant unstable modes and should be smaller than those frequencies at which perturbation growth is expected. The coefficient $\chi_{f}$ should be large enough to stabilize the unstable modes. However, $\chi_{f}$ should not be too large since this will make the convergence very slow. The selection of these parameters has been detailed in \citep{aakervik2006steady,jordi2014encapsulated,jordi2015adaptive,casacuberta2018effectivity}. It should be noted that for the present problem, the SFD method is only applied to the solid region, and therefore the conclusions can be different.

\subsection{The proposed approach and implementation}
In this work, a combining approach based on volume penalization and SFD for IBM treatment is considered and applied to high-order FR schemes. This approach aims to take the advantages of both volume penalization for its robustness and rigorous theoretical foundation, and SFD for the capability of providing a temporal filter to damp the frequencies inside the solid. Therefore, this combination can ensure good accuracy and reduce the oscillations to better satisfy the boundary condition. It should be noted again that both the volume penalization term and the additional SFD term are added only to the solution points immersed in the solid, therefore minimum additional computations are required.

\begin{figure*}[htbp]
	\centering
	\includegraphics[width=350pt]{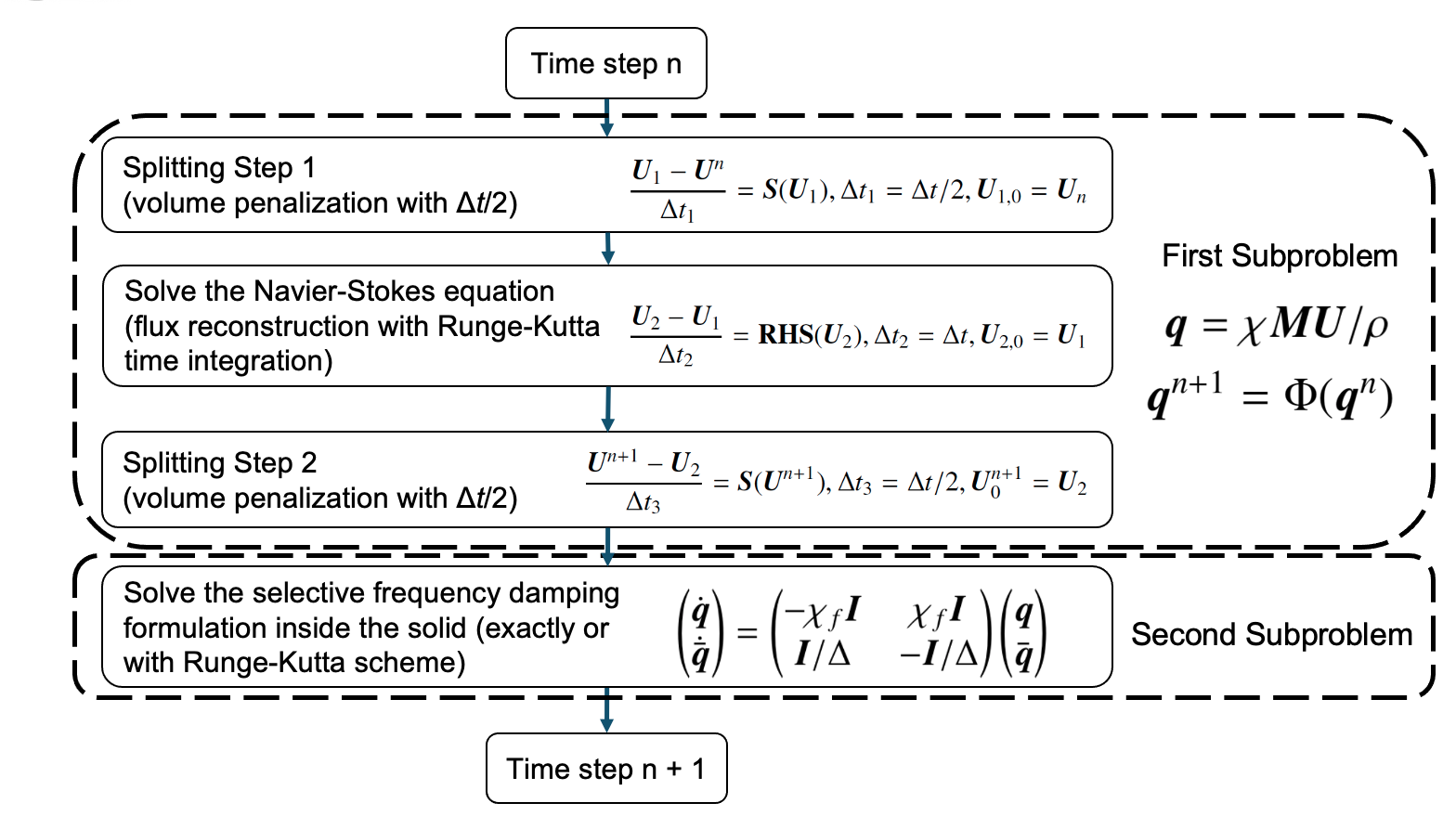}
	\caption{Workflow of the proposed combined volume penalization / selective frequency damping approach for IBM treatment in one time step.}
	\label{fig:workflow}
\end{figure*}

To implement this approach to a Navier-Stokes flow solver with minimal effort, two splitting schemes for the IBM source term and the SFD method are considered. Due to the stiffness of the IBM source term, the first splitting scheme \cite{strang1968construction} is utilized to impose the IBM source term for the Navier-Srokes equations. In addition, the sequential operator-splitting method \citep{farago2004splitting} is considered for the SFD treatment, which is implemented as a wrapper around an existing black-box time-stepping code \citep{jordi2014encapsulated} to avoid modifying the source code for FR discretisation and IBM treatment. The underlying governing equations (written in semi-discrete formulation as $\frac{d \bm{U}}{d t} = \boldsymbol{RHS}(\bm{U})$) can be marched in time by efficient time integration methods. The source term for volume penalization is added by a Strang splitting method \cite{strang1968construction}. This approach is used to handle the stiff source terms of volume penalization \citep{kou2021IBMFR2}, which in practice allows a lower penalization parameter about $\eta \approx \Delta t / 2$. As discussed by Piquet et al. \cite{piquet2016comp}, with Strang splitting, penalization terms are computed exactly for the momentum and energy equations. At time step $n$, the following sequence of operations is performed:

\begin{alignat}{1}
\label{eq:Strang1}
& \text{Splitting step 1}: \frac{\bm{U}_{1} - \bm{U}^{n}}{\Delta t_1} = \bm{S} (\bm{U}_{1}), \Delta t_1 = \Delta t /2, \bm{U}_{1,0} = \bm{U}_n\\
\label{eq:RKevaluation}
& \text{Runge-Kutta time marching}: \frac{\bm{U}_2 - \bm{U}_1}{\Delta t_2} = \mathbf{RHS}(\bm{U}_2), \Delta t_2 = \Delta t, \bm{U}_{2,0} = \bm{U}_1\\
\label{eq:Strang3}
& \text{Splitting step 2}: \frac{\bm{U}^{n+1}-\bm{U}_2}{\Delta t_3} = \bm{S}(\bm{U}^{n+1}), \Delta t_3 = \Delta t/2, \bm{U}_0^{n+1} = \bm{U}_2
\end{alignat}
Currently, in the second step of Equation \ref{eq:RKevaluation}, we use the third-order TVD Runge-Kutta method or the low-storage five-stage fourth-order explicit Runge-Kutta method (LSERK) \citep{carpenter1994fourth,hesthaven2007nodal} method to perform explicit time marching. Numerical tests show that both methods produce identical predictions but the latter allows a relatively larger time step. In each of the splitting step, half of the IBM forcing term is added, which can be achieved by either implicit or explicit forcing methods \citep{kou2021IBMFR1,kou2021IBMFR2}. Here we use a simple explicit treatment, i.e., $\bm{U}_1 =  \bm{U}^n + \Delta t_1 \cdot \bm{S}(\bm{U}^n)$. These three steps in Equation \ref{eq:Strang1}-Equation \ref{eq:Strang3} can be combined into the first subsystem. It should be noted that in this work, the Strang splitting approach is only considered for the Navier-Stokes equations. Analyses of the advection equation are performed in a more general situation where direct imposition of the IBM source term is considered. After that, an encapsulated version of the SFD
\citep{jordi2014encapsulated} is implemented as the second subsystem, where the vector $\boldsymbol{q}$ contains only flow velocities inside the solid. Therefore, for the Navier-Stokes equations, we define $\boldsymbol{q}$ as

\begin{equation}
    \boldsymbol{q} = \mask  [u,v,w]^T = \mask  \boldsymbol{M} \boldsymbol{U} / \rho \ , \  \boldsymbol{U} = [\rho, \rho u, \rho v, \rho w, E]^T,
\end{equation}
where $\boldsymbol{U}$ only contains the solution points with $\mask = 1$, and $\boldsymbol{M}$ is the matrix to select the velocities from the state vector

\begin{equation}
    \boldsymbol{M} = \begin{pmatrix}
    0 & 1 & 0 & 0 & 0 \\
    0 & 0 & 1 & 0 & 0 \\
    0 & 0 & 0 & 1 & 0 
    \end{pmatrix}.
\end{equation}

A function $\Phi$ is introduced for the first subsystem such that the numerical solution at time step $n$ (without applying the SFD method) is given by $\boldsymbol{q}^{n+1} = \Phi (\boldsymbol{q}^{n})$. We now solve the second subsystem with the SFD method that represents the actions of the feedback control and the low-pass time filter. This leads to the following formulation

\begin{equation}
    \left\{\begin{matrix}
    \dot{\boldsymbol{q}} = - \chi_{f} (\boldsymbol{q} - \bar{\boldsymbol{q}})
\\  \dot{\bar{\boldsymbol{q}}} = (\boldsymbol{q} - \bar{\boldsymbol{q}}) / \Delta
\end{matrix}\right. ,
\end{equation}
which has the matrix form:
\begin{equation}
    \begin{pmatrix}\dot{\boldsymbol{q}}\\ \dot{\bar{\boldsymbol{q}}} \end{pmatrix} = \begin{pmatrix} - \chi_{f} \boldsymbol{I} & \chi_{f} \boldsymbol{I} \\ \boldsymbol{I}  / \Delta & - \boldsymbol{I}  / \Delta \end{pmatrix} \begin{pmatrix}\boldsymbol{q}\\ \bar{\boldsymbol{q}} \end{pmatrix} ,
\end{equation}
where $\boldsymbol{I}$ is the identity matrix. At the same time step $n$, this equation can be solved exactly as 

\begin{equation}
\label{eq:SFD}
    \begin{pmatrix}\boldsymbol{q}^{n+1}\\ \bar{\boldsymbol{q}}^{n+1} \end{pmatrix} = e ^ {\mathcal{L} \Delta t} \begin{pmatrix}\Phi (\boldsymbol{q}^{n})\\ \bar{\boldsymbol{q}}^n \end{pmatrix} ,
\end{equation}
or solved by any explicit time marching method. The linear operator $\mathcal{L}$ is defined such that
\begin{equation}
     e ^ {\mathcal{L} \Delta t} = \frac{1}{1+\chi_f \Delta} 
     \begin{pmatrix} \boldsymbol{I} + \chi_f \Delta \boldsymbol{I} e ^ {-(\chi_f + \frac{1}{\Delta}) \Delta t} & \chi_f \Delta \boldsymbol{I} [1 - e^{-(\chi_f + \frac{1}{\Delta}) \Delta t}] \\ 
     \boldsymbol{I} - \boldsymbol{I} e ^{-(\chi_f + \frac{1}{\Delta}) \Delta t} & \chi_f \Delta \boldsymbol{I} + \boldsymbol{I} e ^{-(\chi_f + \frac{1}{\Delta}) \Delta t}
     \end{pmatrix} .
\end{equation}

The whole scheme is called “encapsulated” since $\Phi$ is not modified but simply used as an input of the linear solver in Equation \ref{eq:SFD} \citep{jordi2014encapsulated}. The final $\boldsymbol{q}^{n+1}$ from Equation \ref{eq:SFD} will be the updated velocity for the solution points inside the solid. Moreover, the initial value for $\bar{\boldsymbol{q}}$ is set to the known boundary velocity $\boldsymbol{u}_s$, which may change slightly as the solution evolves. This indicates the importance of using the volume penalization since it should be used to ensure the accuracy of approximating the exact boundary conditions. The entire workflow is summarized in Figure \ref{fig:workflow}.

To apply the proposed method, three parameters should be determined, including the penalization parameter $\IBMparam$ for volume penalization, as well as the control parameter $\chi_f$ and filter width $\Delta$ for SFD. As discussed in previous sections, the selection of $\IBMparam$ depends on the stiffness of the IBM source term and the explicit time step $\Delta t$. Moreover, from the SFD formulation, $\chi_f$ also works similar to a penalization parameter. Practical guidance for the selection of $\chi_f$ and $\Delta$ will be determined by numerical tests in the following section. 

\section{Validation}
In this section, the proposed IBM approach based on the SFD method and volume penalization is tested. Firstly, the one-dimensional advection equation with a no-slip wall in the middle of the computational domain is considered, where numerical experiments and von Neumann analyses are utilized to provide guidelines for the selection of important parameters. Secondly, the proposed method is tested by simulating steady flows past a NACA0012 airfoil based on the Navier-Stokes equations. Thirdly, unsteady flow past a cylinder is simulated in order to validate the efficacy of the method in unsteady problems. In the latter two cases, the main focus is to compare additional damping effects through observing the temporal evolution of velocities inside the solid.

\subsection{Linear advection equation}
We perform numerical experiment and eigensolution analyses for the linear advection equation in Equation \ref{eq:adv}. We consider an advection speed $\advcoef = 1$ and a computational domain defined in $x \in [-1,1]$ discretized by $N$ equispaced elements with mesh size $h$. Periodic boundary conditions are imposed at both sides of the computational domain. The fully upwind flux with $\lambda = 1$ is chosen. The solid region is defined as a no-slip wall, i.e., $\velx_s = 0$ in Equation \ref{eq:adv}. It lies in the middle of the computational domain and starts from $x=0$, whose width is defined as $\delta$:
\begin{equation}
\mask  (x,t) = \left\{\begin{matrix}1, \, \, \text{if}\, \, 0 < x < \delta 
\\0,\, \, \text{otherwise}
\end{matrix}\right. .
\end{equation}

\begin{figure*}[htbp]
		\centering
		\includegraphics[width=300pt]{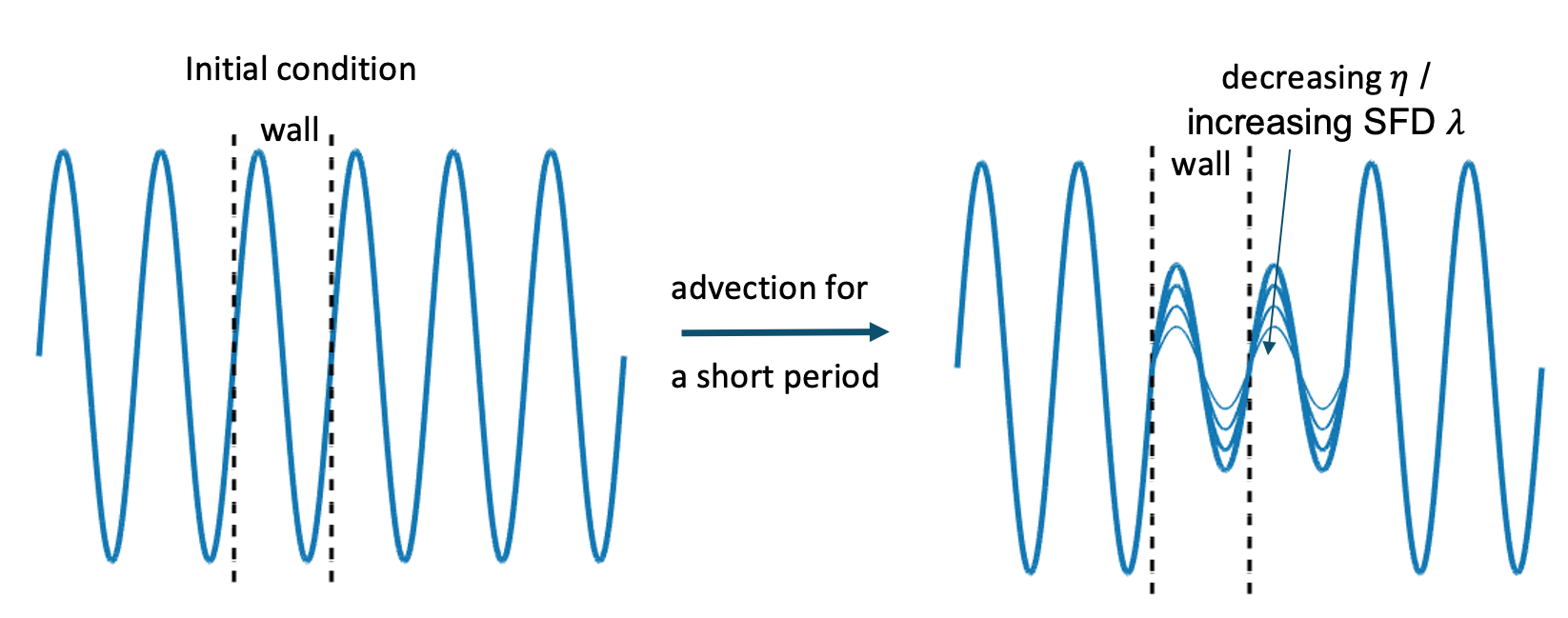}
		\caption{Schematic illustration of the advection problem with IBM.}
		\label{fig:cartoon}
\end{figure*}

We consider $\delta$ as integer multiples of element size $h$ ($\delta = Zh$, where $Z$ is an integer), thus the solid boundaries will lie exactly at the interface between elements. In particular, we have $Z = 1$. This allows us to define the solid ratio $r = Z/N$ as the ratio between the solid region and the computational domain. The matrix form of space discretisation based on FR is detailed in \ref{append2}, which is the basis of eigensolution analysis. A schematic illustration of the present problem is shown in Figure \ref{fig:cartoon}. The initial wavelike solution starts to pass through the no-slip wall in the middle, and the damped solution will move towards right as time evolves and the global solution will eventually become $0$ due to the IBM conditions. In a short time period between $0 < x< t$, the solution is expected to approach $0$ since it comes out from the wall. However, in practice, the value of this solution depends on the damping provided by the volume penalization and the SFD methods. Therefore the accuracy of IBM treatment can be evaluated by comparing the solution with zero (within a short advection time $0 < x< t$ ), which is the expected exact solution.

We consider the numerical experiment of a linear advection equation with a wavelike initial condition. The initial condition is defined by a sinusoidal wave with wavenumber $k$, which is non-dimensionalized by the mesh size $h$ and polynomial order $P$ and defined as $kh/(P+1)$. In addition, due to the existence of a solid wall, the actual fluid domain is shorter than the whole computational domain, therefore the effective wavenumber in the fluid region is higher than $k$. To obtain this effective wavenumber, the wavenumber is re-scaled by the solid ratio $r$ as $ \hat{k} = k / (1-r)$ \citep{kou2021VonNeumann}. We consider a space discretisation with $N = 40$ elements in the computational domain ($h=0.05$). Based on this mesh, we set $\delta = h$ with $r = 1/40$ and $P = 3$ as a representative order for high-order methods. The initial condition with wavenumber $\hat{k} h/(P+1) = 0.3223$ is considered, which lies in the resolved wavenumber region of the scheme. 

\begin{figure*}[htbp]
    \begin{subfigure}{.38\textwidth}
 		\centering
		\includegraphics[width=170pt]{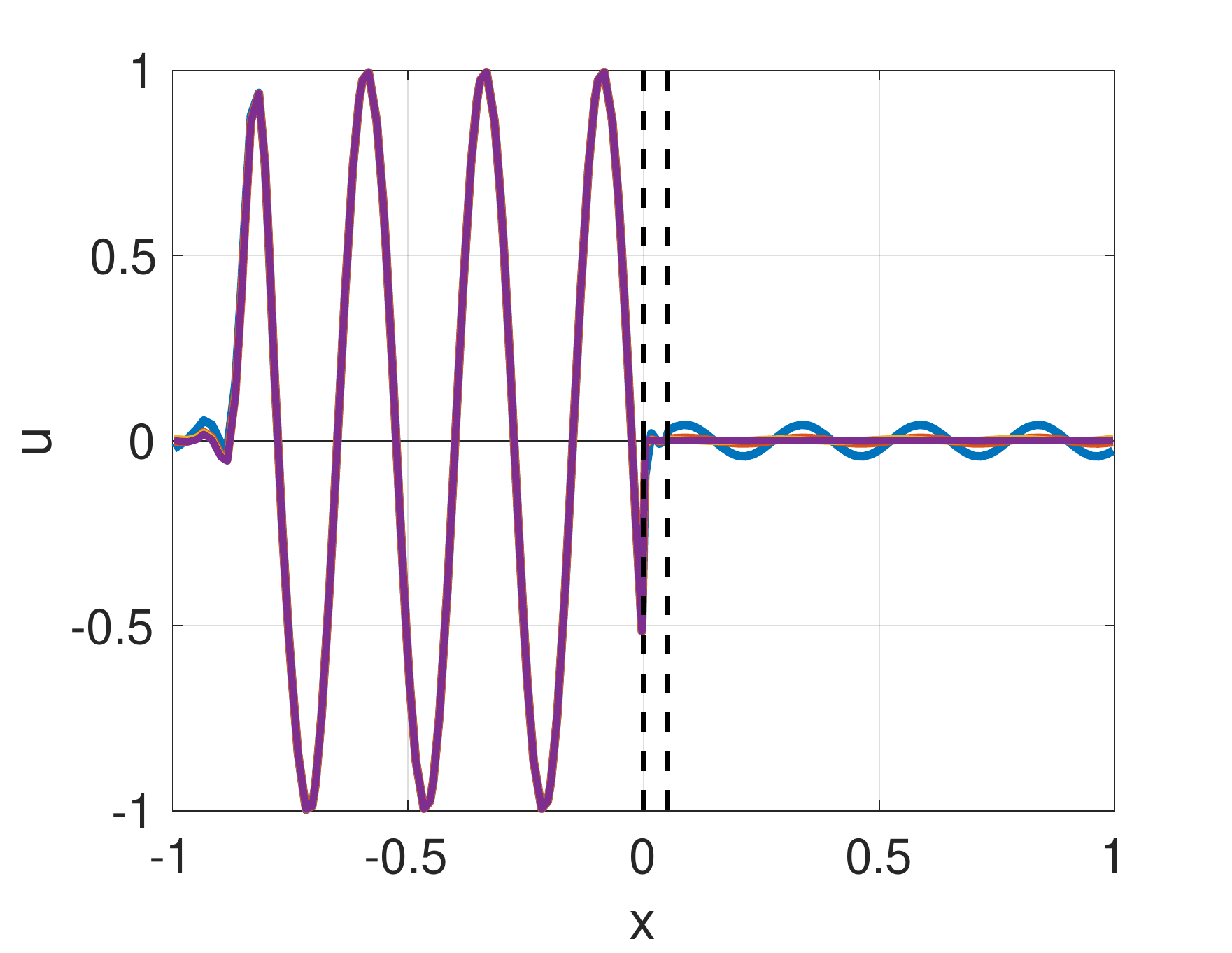}
		\caption{}
	\end{subfigure}
	\begin{subfigure}{.38\textwidth}
 		\centering
		\includegraphics[width=170pt]{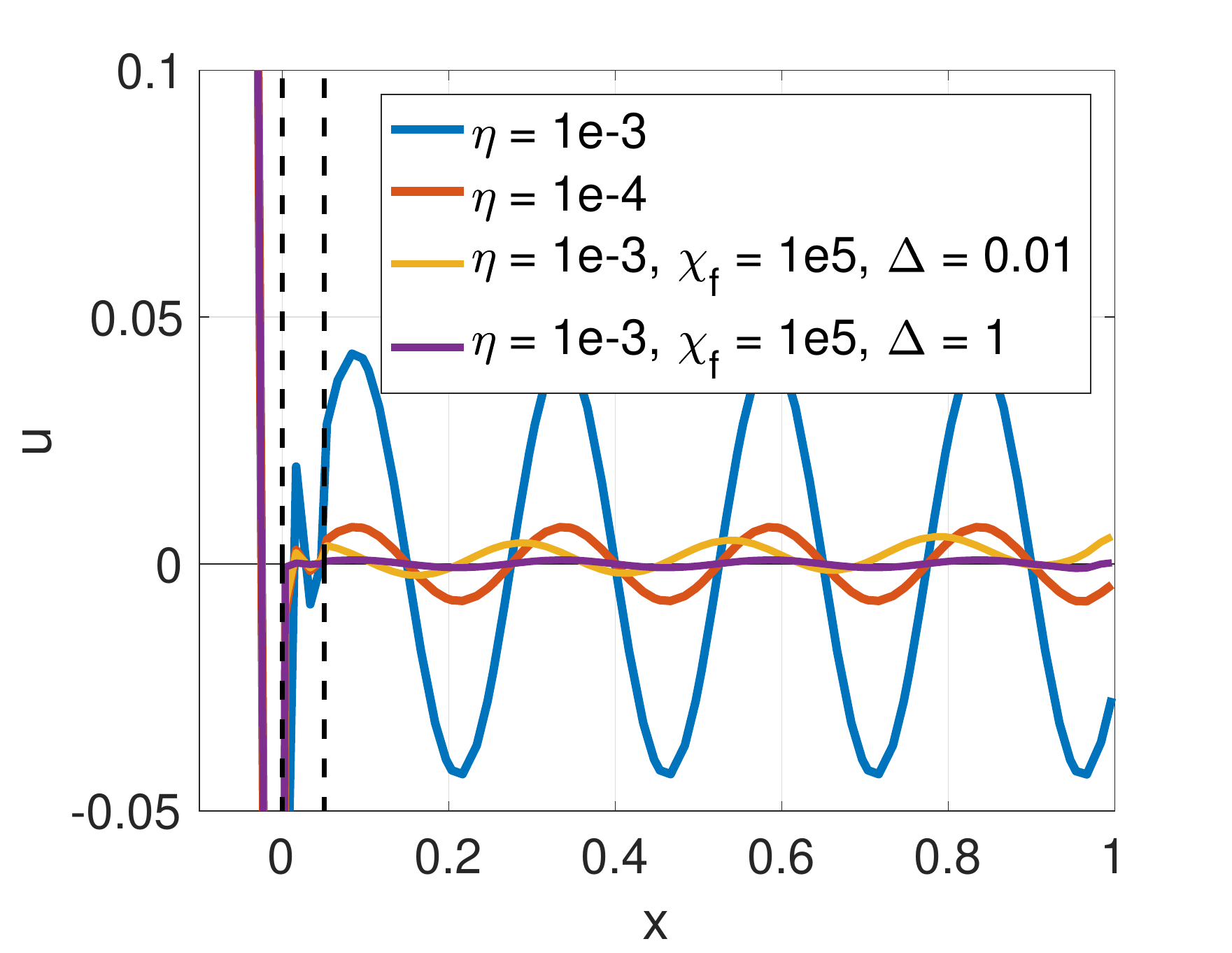}
		\caption{}
	\end{subfigure}
	\centering
	\caption{Simulation under different parameters ($r = 1/40$, $P = 3$, initial wavenumber $\hat{k} h/(P+1) = 0.3223$, $N = 40$). a) Global view. b) Zoom-in view.}
	\label{fig:advection-sim}
\end{figure*}

To compare different combinations of penalization parameter and SFD parameters, we choose four cases: 1) only volume penalization with $\IBMparam = 1 \times 10^{-3}$, 2) only volume penalization with $\IBMparam = 1 \times 10^{-4}$, 3) volume penalization with $\IBMparam = 1 \times 10^{-3}$, and the SFD method with $\chi_{f} = 1 \times 10^{5}$ and $\Delta= 0.01$, 4) volume penalization with $\IBMparam = 1 \times 10^{-3}$, and the SFD method with $\chi_{f} = 1 \times 10^{5}$ and $\Delta= 1$. Figure \ref{fig:advection-sim} shows the transient simulation results of different test cases. The final simulation time is set to $1.1$, in order to guarantee that the solution in the computational domain $x \in [\delta, 1]$ is sufficiently penalized. An explicit third-order Runge-Kutta time-marching scheme is used for time integration. To reduce the temporal error, the time step is set to $\Delta t = 1 \times 10^{-5}$. The SFD treatment is achieved by the encapsulated approach in \citep{jordi2014encapsulated}. From the first two cases, as $\IBMparam$ decreases, a larger penalization source term is imposed, therefore the solution will approach the expected value (zero). Alternatively, we can also keep a relatively larger $\IBMparam = 1 \times 10^{-3}$ but use the SFD method instead, which leads to improved accuracy compared to volume penalization with $\IBMparam = 1 \times 10^{-4}$. This indicates the efficacy of the SFD method to impose the boundary condition. In addition, from the third and the fourth cases, the filter width $\Delta$ will also impact the accuracy of the SFD method. When $\Delta$ is larger, improved accuracy is seen. Therefore, the influence and the selection of these parameters, including $\IBMparam$, $\chi_{f}$ and $\Delta$, should be investigated in detail. 

\begin{figure*}[htbp]
    \begin{subfigure}{.33\textwidth}
		\includegraphics[width=150pt]{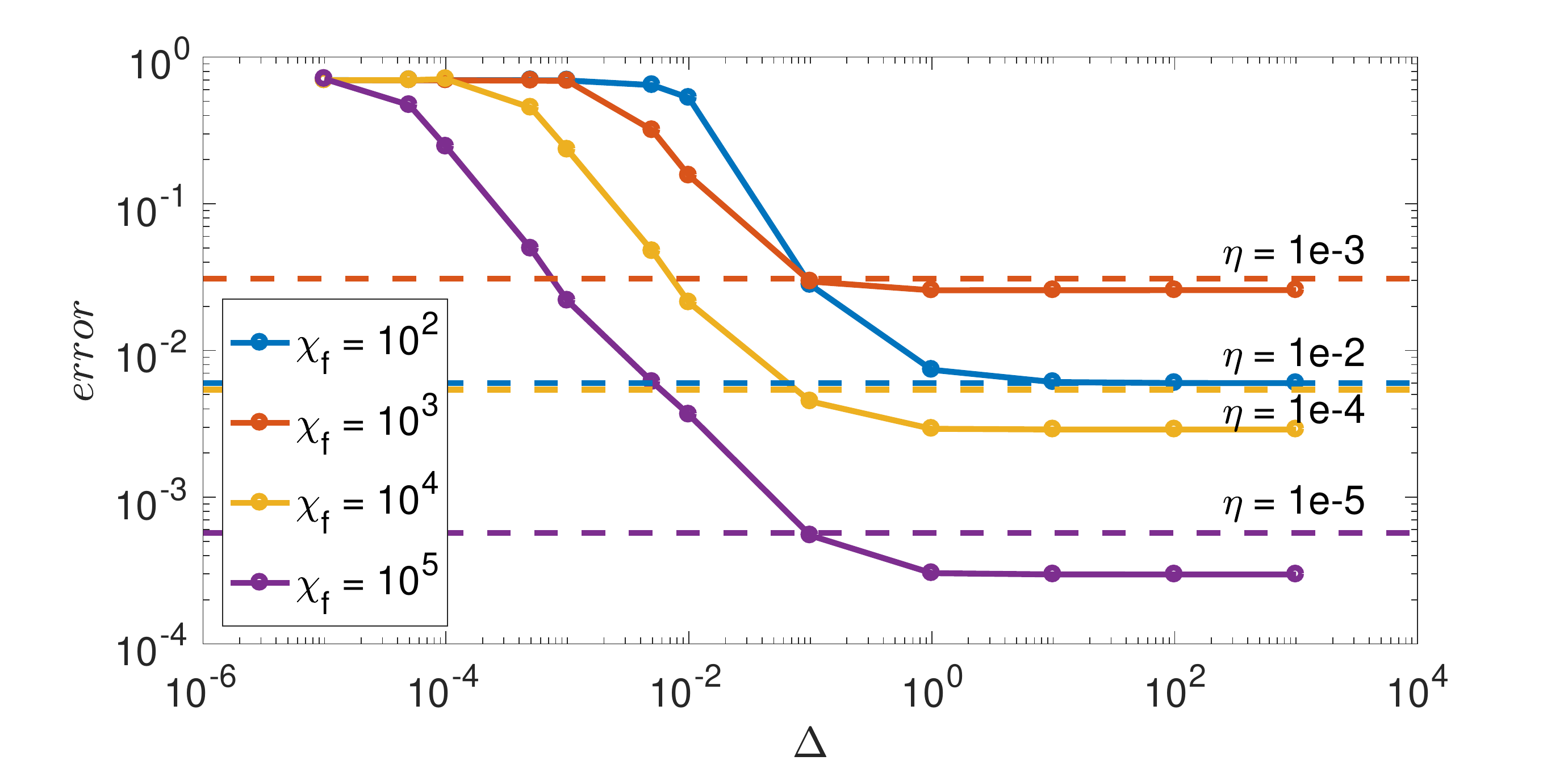}
		\caption{}
	\end{subfigure}
	\begin{subfigure}{.33\textwidth}
		\includegraphics[width=150pt]{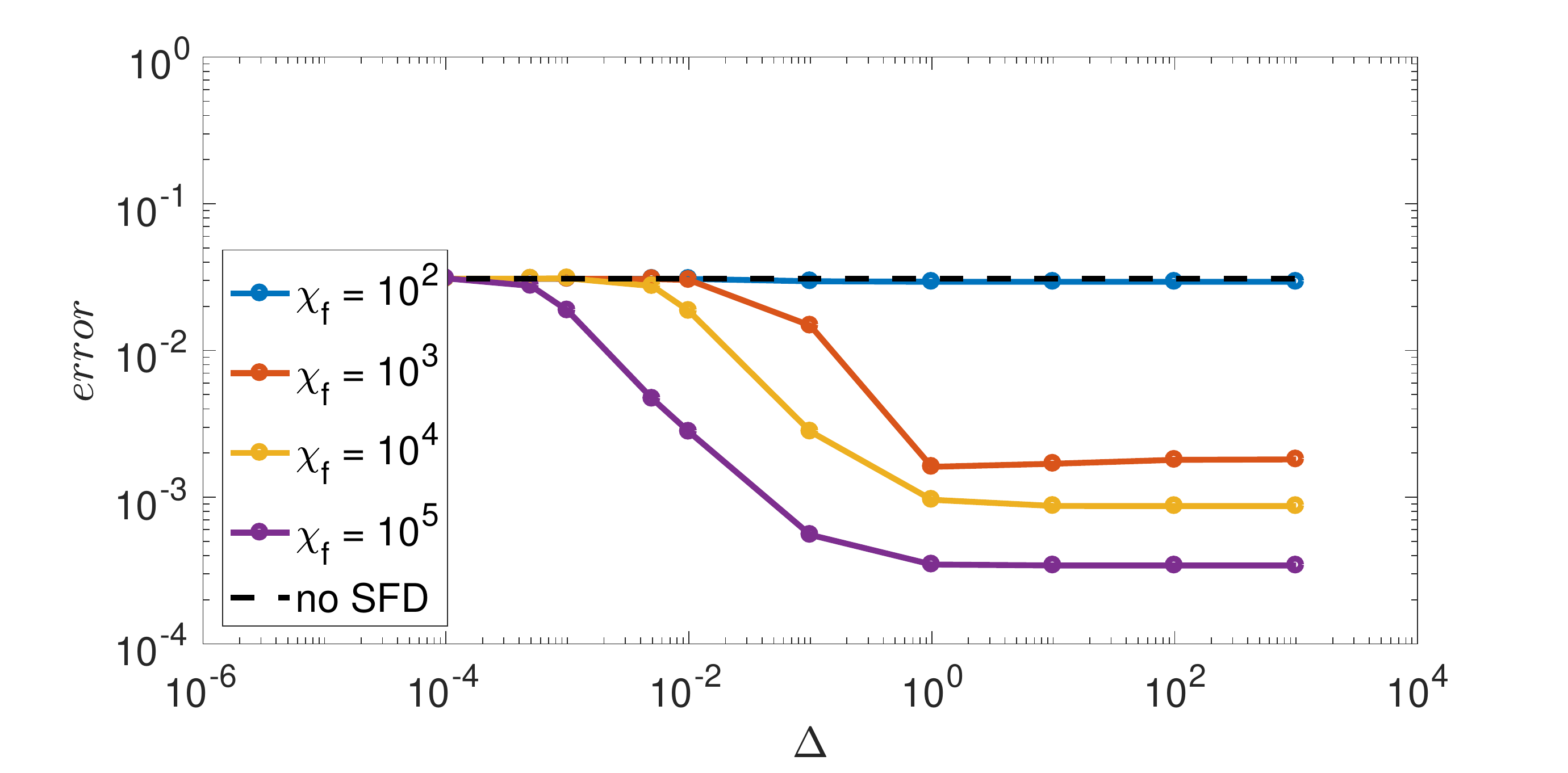}
		\caption{}
	\end{subfigure}
	\begin{subfigure}{.33\textwidth}
		\includegraphics[width=150pt]{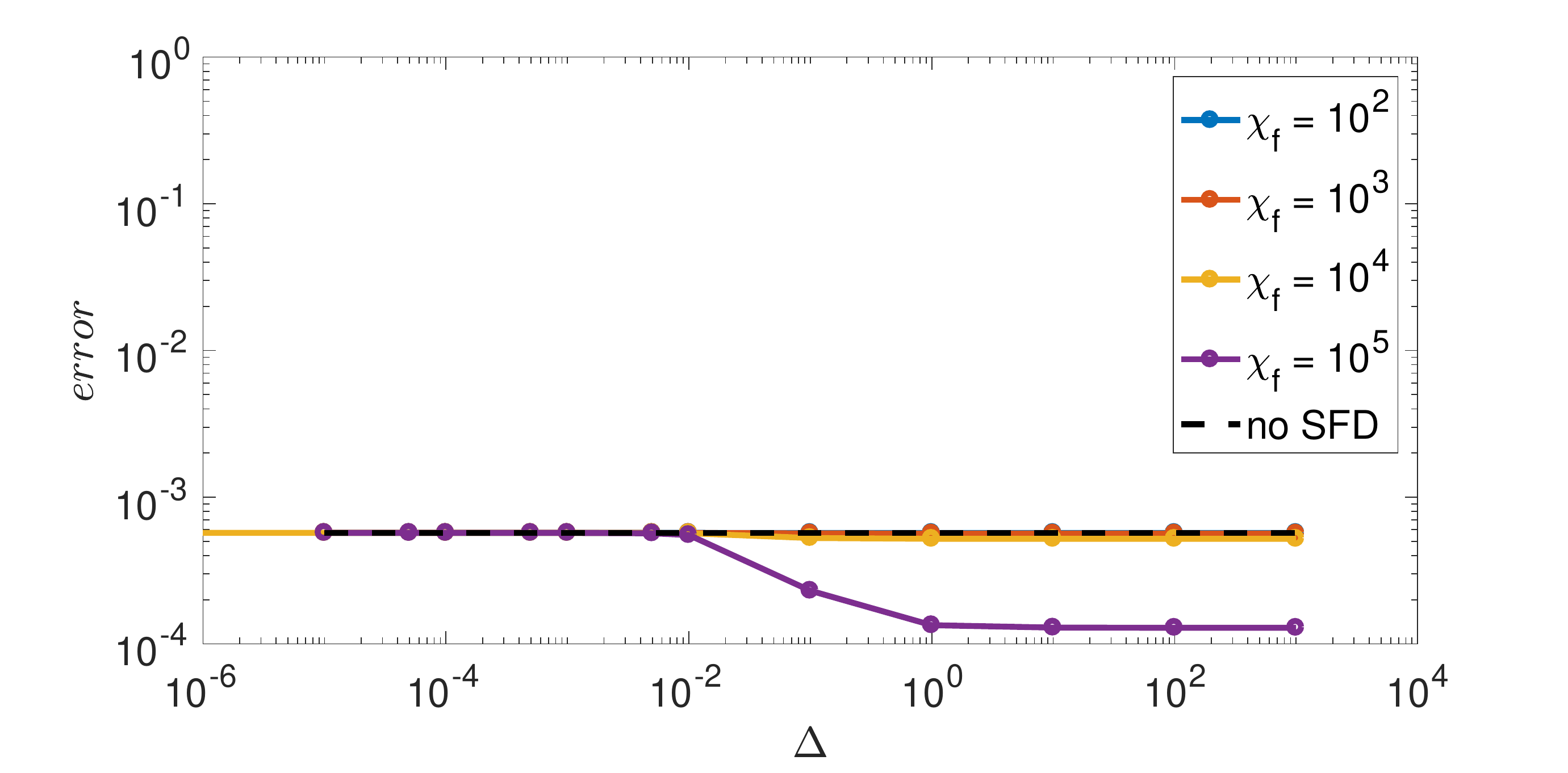}
		\caption{}
	\end{subfigure}
	\begin{subfigure}{.33\textwidth}
		\includegraphics[width=150pt]{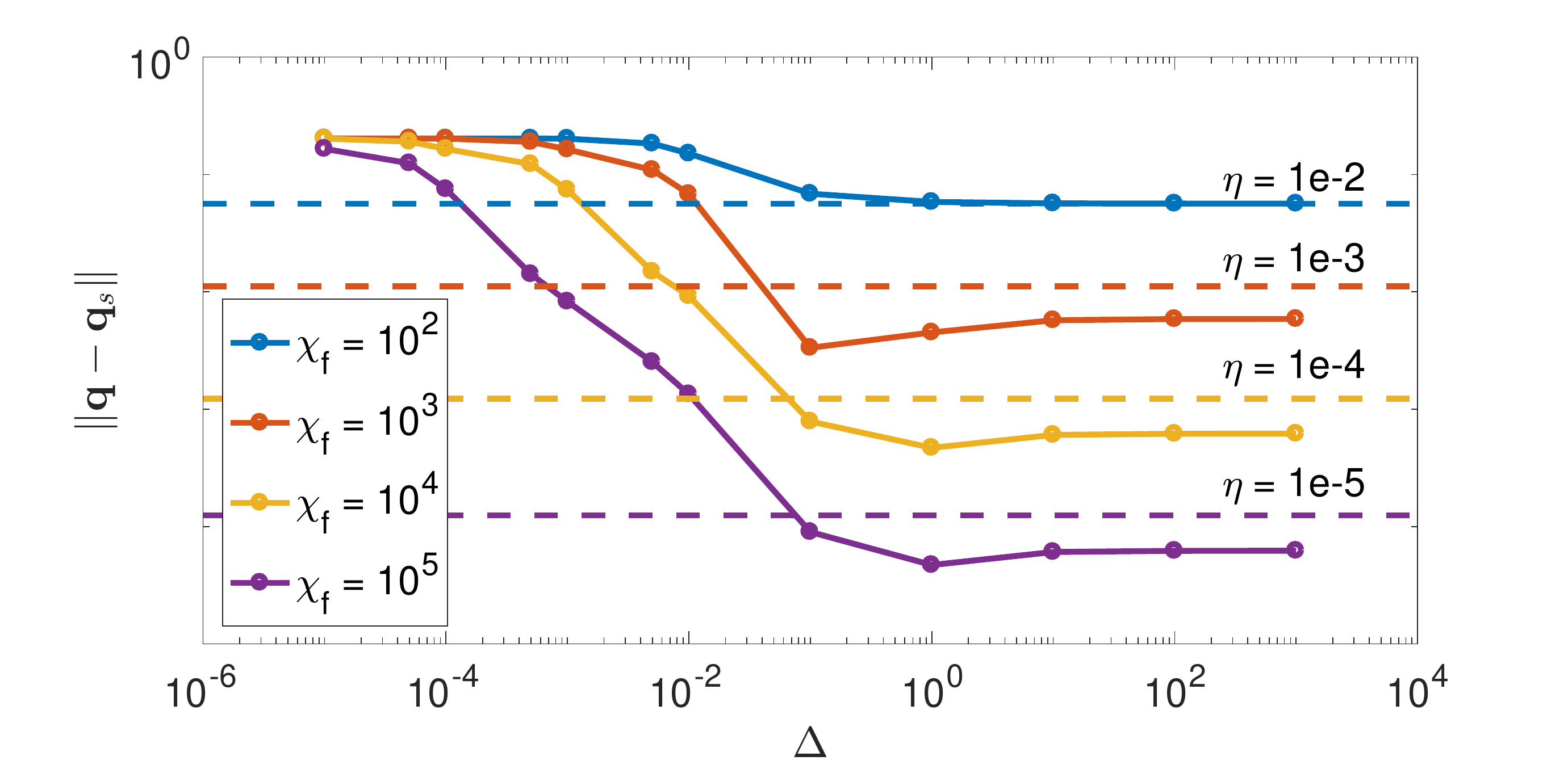}
		\caption{}
	\end{subfigure}
	\begin{subfigure}{.33\textwidth}
		\includegraphics[width=150pt]{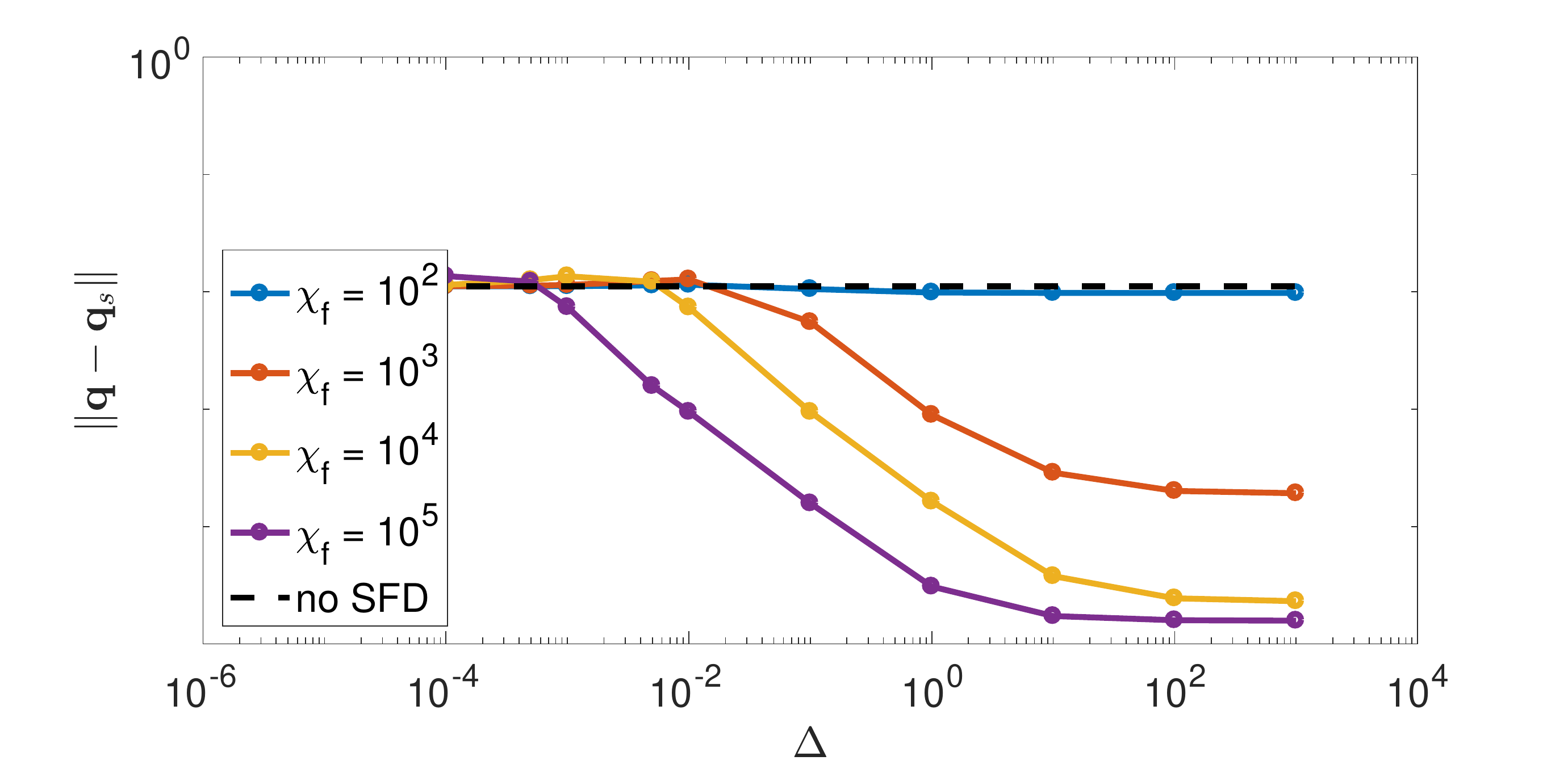}
		\caption{}
	\end{subfigure}
	\begin{subfigure}{.33\textwidth}
		\includegraphics[width=150pt]{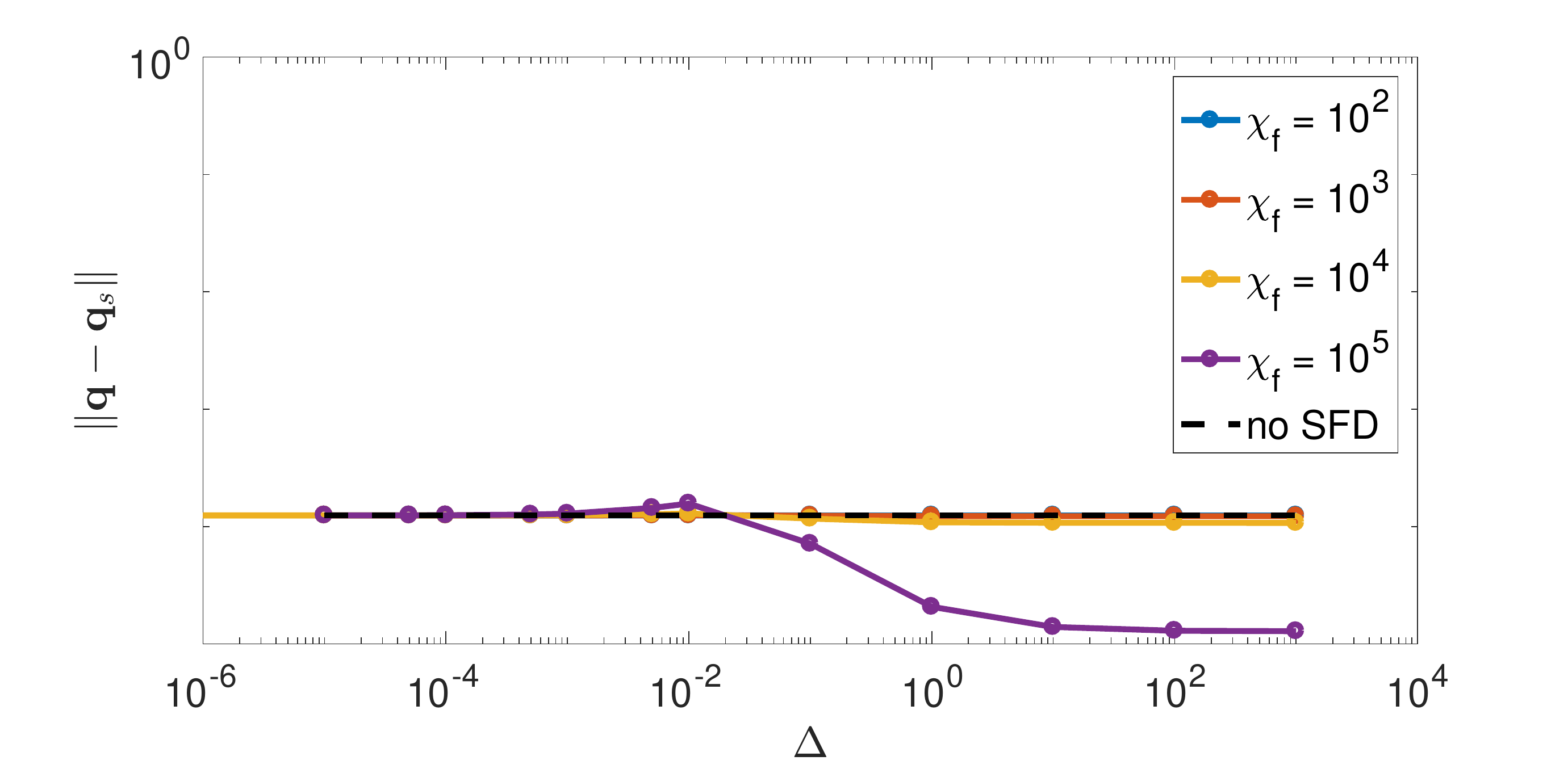}
		\caption{}
	\end{subfigure}
	\centering
	\caption{Error comparison inside flow or solid regions under different parameters. a) Only SFD, error in flow. b) $\IBMparam = 1 \times 10^{-3}$, error in flow. c) $\IBMparam = 1\times 10^{-5}$, error in flow. d) Only SFD, error in solid. e) $\IBMparam = 1 \times 10^{-3}$, error in solid. f) $\IBMparam = 1\times 10^{-5}$, error in solid. Dashed lines refer to the error when only volume penalization is used. }
	\label{fig:advection-err}
\end{figure*}

To further investigate the relationship between accuracy and the related parameters, the simulation error is compared. Both errors in the solid and fluid are compared between the solution and the penalized value $\velx_s = 0$, defined in $x \in [0, \delta]$ and $x \in [\delta, 1]$, respectively. Defining the number of solution points inside the domain of interest as $L$, we have the mean squared error for the flow region
\begin{equation}
    error = \sqrt{\frac{\sum_{i=1}^{L} [\velx(x_i) - \velx^{exact}(x_i)]^2}{L}} \ , \ x_i \in [\delta, 1], \ \velx^{exact} = 0,
\end{equation}
and the mean squared error for the solid region (the vector $\boldsymbol{q}$ in the SFD method)
\begin{equation}
    error_s = ||\boldsymbol{q} - \boldsymbol{q}_s || = \sqrt{\frac{\sum_{i=1}^{L_s} [\velx(x_i) - \velx^{exact}(x_i)]^2}{L}} \ , \ x_i \in [0, \delta], \ \velx^{exact} = 0.
\end{equation}

Again, we choose the wave-like initial condition with wavenumber $\hat{k} h/(P+1) = 0.3223$ and set the final time to $t = 1.1$. A group of simulations is performed, with different combinations of parameters. Generally, three selections of $\IBMparam$ for volume penalization are considered: 1) $\IBMparam \rightarrow \infty$, which means only the SFD method is imposed; 2) $\IBMparam = 1 \times 10^{-3}$, representing a medium level of volume penalization; 3) $\IBMparam = 1 \times 10^{-5}$, representing a high level of volume penalization. In each group of case, variations of $\chi_{f}$ and $\Delta$ are considered. The control coefficient $\chi_f$ varies from $1 \times 10^2$ to $1 \times 10^5$. The filter width $\Delta$ varies from $1 \times 10^{-5}$ to $1 \times 10^3$, indicating the filtered frequencies varying from very large values (only very large frequencies are damped) to small values (almost all frequencies are damped). For each case, to reduce the cost, the time step is calibrated according to the parameters. In addition, we also find that by using the encapsulated scheme  \citep{jordi2014encapsulated} for the SFD method with proper time step, lower errors can be reached. 

All errors are compared in Figure \ref{fig:advection-err}. From this figure, several conclusions can be drawn. It is clearly seen that as the filter width $\Delta$ is increased, smaller errors are obtained and the minimal error is reached when $\Delta > 10$. This suggests that when more frequencies are suppressed, the boundary condition will be imposed more accurately. This is not surprising, since the purpose of using the SFD method here is to damp all the oscillation frequencies inside the solid. This is the main difference between frequency damping for IBM, as opposed to the traditional SFD method (used in instability analyses). In traditional SFD method, only the unstable frequencies need to be damped, therefore $\Delta$ needs to be calibrated carefully \citep{aakervik2006steady,jordi2014encapsulated}. In addition, as shown in Figure \ref{fig:advection-err}a and d, when only the SFD method is used and for a given $\chi_f$, it can reach similar but smaller error compared to that of volume penalization when $\chi_f = 1 / \IBMparam$. This is also reasonable, since the SFD source term works like a penalization source term, it may show similar performance when $\chi_f \approx 1 / \IBMparam$ (but the SFD method still slightly outperforms the volume penalization). Furthermore, when the volume penalization is combined with the SFD method, both errors in the flow and solid can be further reduced as $\IBMparam$ decreases. This is because volume penalization provides additional penalization to impose the boundary conditions. However, it should be noted that when $\chi_f \geq 1 / \IBMparam$, the benefit from the SFD method becomes more evident. This is more pronounced in Figure \ref{fig:advection-err} c and f, where an additional benefit is obtained when $\chi_f = 1 / \IBMparam = 1 \times 10^{5}$. These results indicate that when both volume penalization and the SFD method are applied for IBM treatment, $\chi_f \geq 1 / \IBMparam$ will lead to improved performance compared with only volume penalization. However, it is also observed in these numerical experiments that when $\chi_f = 1 / \IBMparam$, a comparable time step for stability is reached, while when $\chi_f \geq 1 / \IBMparam$, one may need to reduce the time step. Therefore, due to the stability consideration, a guideline of $\chi_f = 1 / \IBMparam$ can be considered. The eigensolution analysis based on this selection will be given in following sections to justify this guideline.

\begin{figure}[htbp]
    \begin{subfigure}{.38\textwidth}
 		\centering
		\includegraphics[width=170pt]{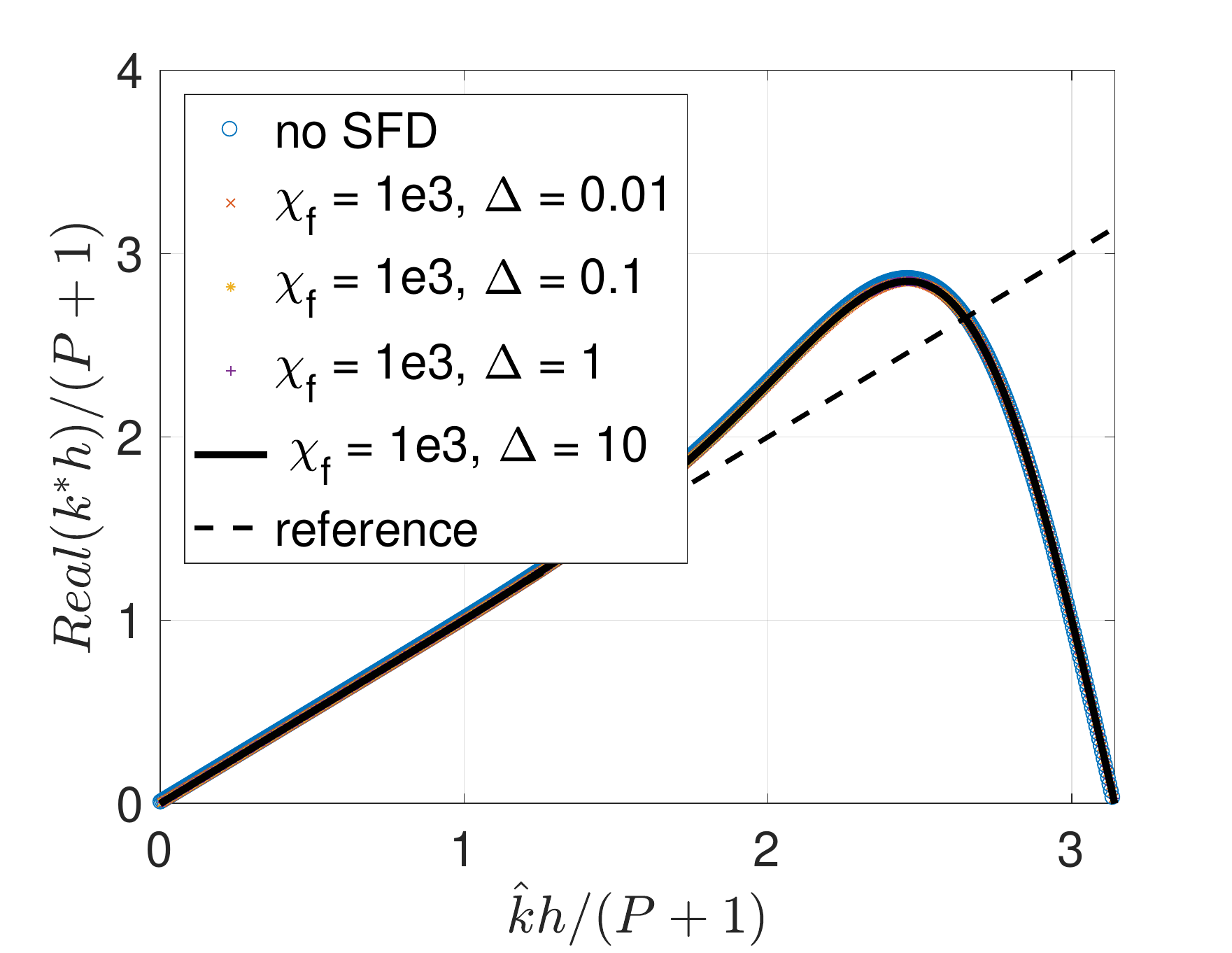}
		\caption{}
	\end{subfigure}
	\begin{subfigure}{.38\textwidth}
 		\centering
		\includegraphics[width=170pt]{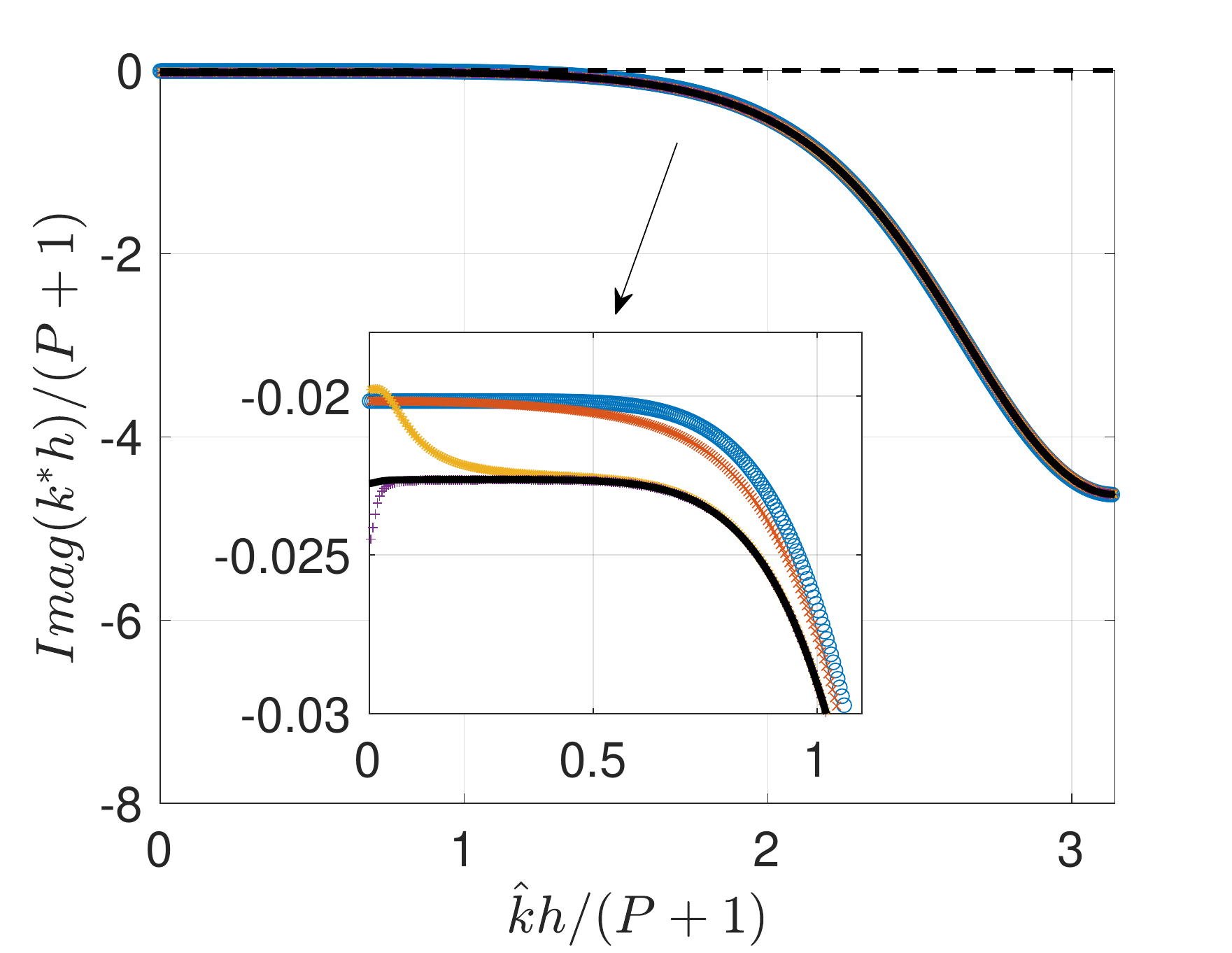}
		\caption{}
	\end{subfigure}
	\centering
	\caption{Semi-discrete dispersion-dissipation curve for volume penalization with and without SFD. $\IBMparam$ for volume penalization is set to $1 \times 10^{-3}$. a) Dispersion. b) Dissipation.}
	\label{fig:advection-semi}
\end{figure}

To better understand the proposed SFD method, an eigensolution analysis \citep{hu1999analysis,moura2015linear,kou2021VonNeumann} is introduced. It applies eigendecomposition to the discretized operator, and extracts the dispersion and dissipation behavior across different wavenumbers. For eigensolution analysis of high-order methods, more details can be found in \citep{moura2015linear,vincent2011insights}. Since IBM source term is only introduced to the solid region, and the global discretisation matrix in Equation \ref{eq:global} is used, the physical mode that recovers $k^* = k$ as $k \rightarrow 0$ needs to be extracted properly \citep{manzanero2018dispersion,kou2021VonNeumann}. We consider the same test condition with $N = 40$, $r = 1/40$ and $P = 3$.

\begin{figure}[htbp]
    \begin{subfigure}{.38\textwidth}
 		\centering
		\includegraphics[width=170pt]{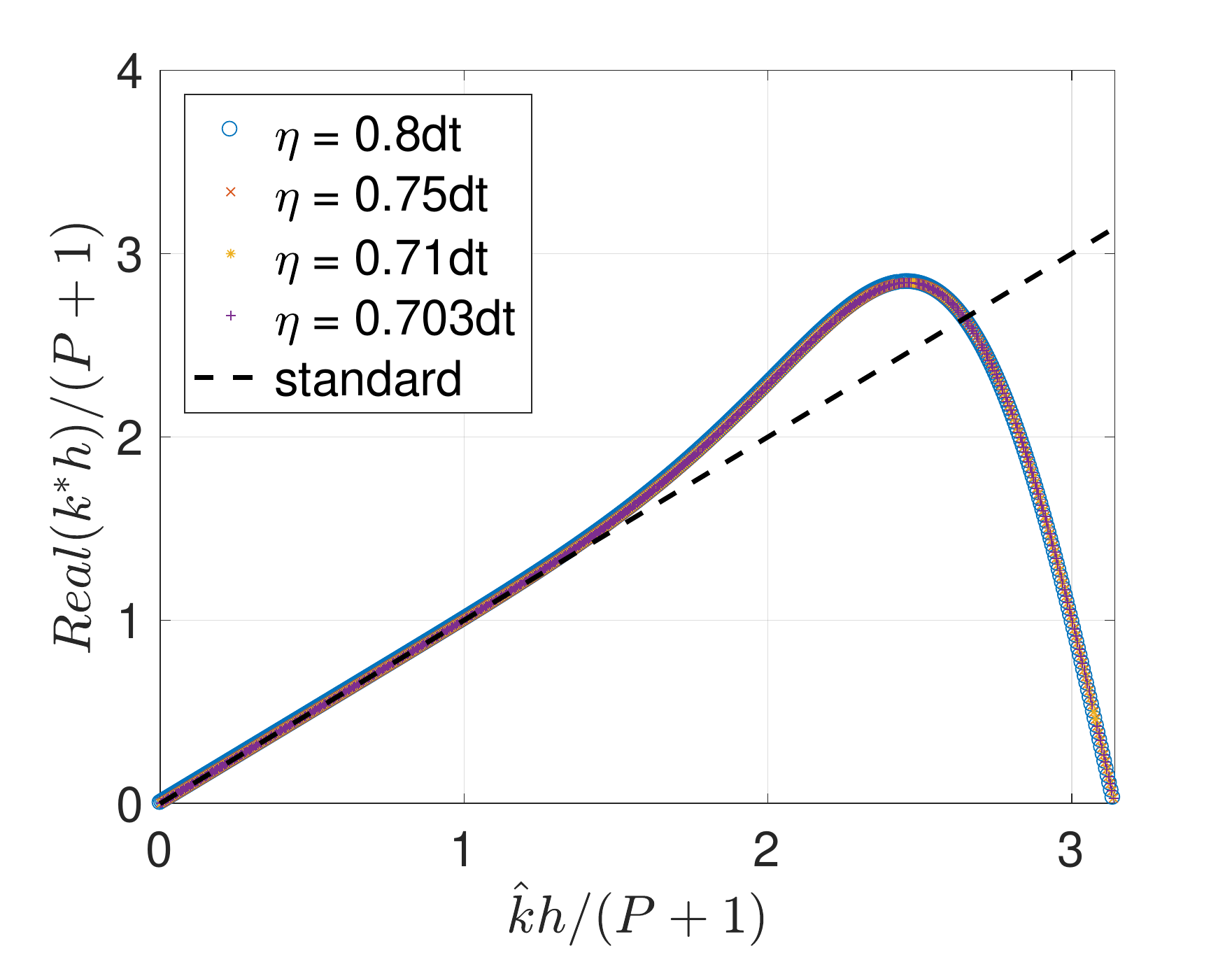}
		\caption{}
	\end{subfigure}
	\begin{subfigure}{.38\textwidth}
 		\centering
		\includegraphics[width=170pt]{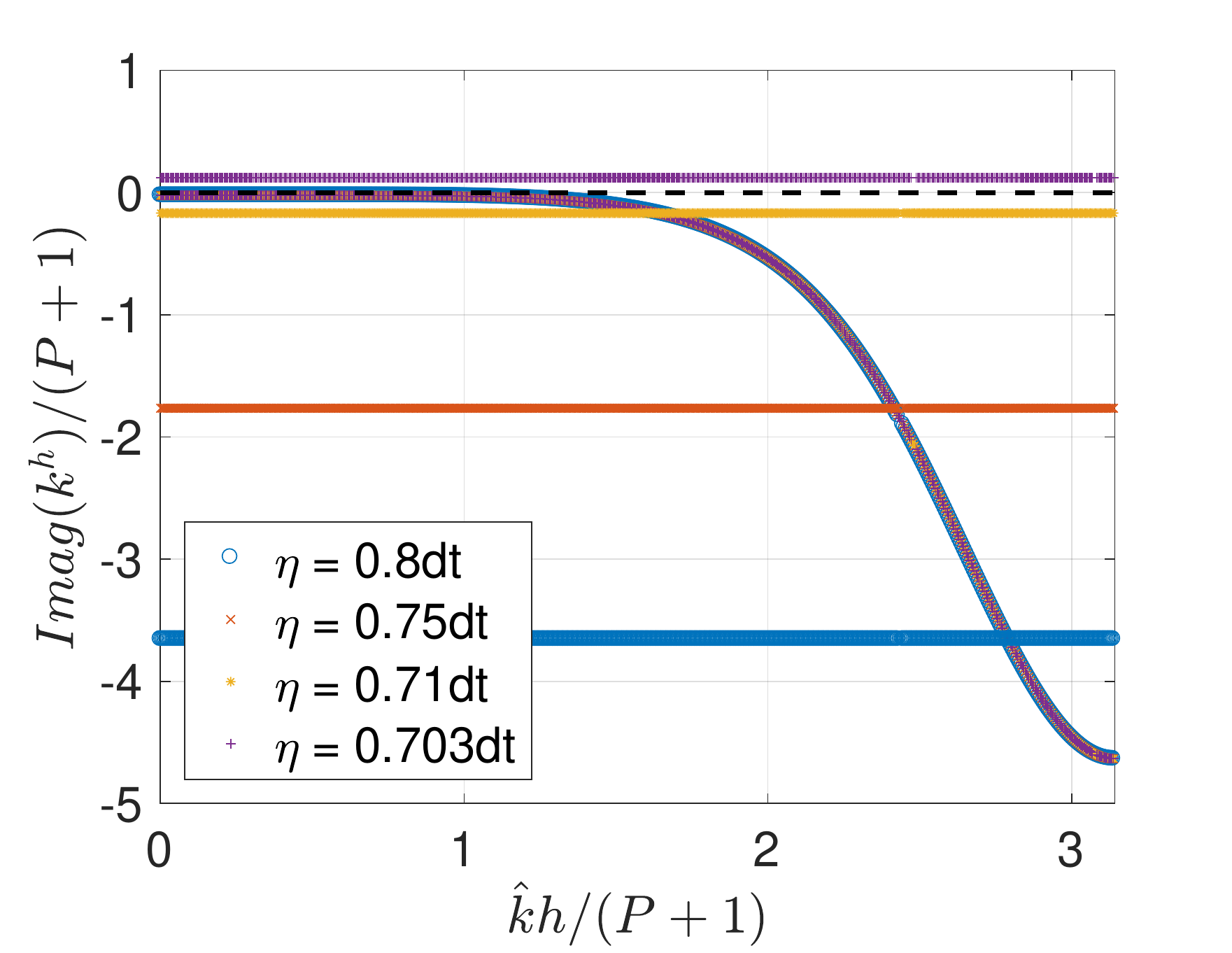}
		\caption{}
	\end{subfigure}
	\centering
	\caption{Fully-discrete dispersion-dissipation curve for the combined approach ($\Delta = 1 \times 10^{-3}$, $\IBMparam = 1 / \chi_f$). a) Dispersion of the physical mode. b) Dissipation of the physical mode and the solid mode.}
	\label{fig:advection-full}
\end{figure}

Firstly, the semi-discrete analysis which only considers the space discretisation is performed and the resulting physical mode is compared in Figure \ref{fig:advection-semi}. A medium penalization parameter $\IBMparam = 1 \times 10^{-3}$ is chosen for volume penalization. In addition, the SFD method with $\chi_f = 1/\IBMparam =  1 \times 10^{3}$ and different filter widths are considered. From Figure \ref{fig:advection-semi}a, it is clear that the additional SFD treatment does not affect the dispersion behavior and maintains good resolution in the resolved wavenumber range. It should be noted that the comparison between the volume penalization and standard advection equation is studied in detail by the authors in \citep{kou2021VonNeumann}. The improved accuracy of combining volume penalization and the SFD method can be observed in the dissipation curve in Figure \ref{fig:advection-semi}b. When the SFD method is added and the filter width is large enough, the physical mode shows additional dissipation effect across the entire wavenumber range. This indicates that the solution will be further damped due to the improved imposition of boundary conditions \citep{kou2021VonNeumann}. In addition, the effect of filter width is evident. For example, when $\Delta = 0.01$, only high frequency component (non-dimensional wavenumber larger than about 0.5) benefits from the SFD method, with very small additional dissipation. When $\Delta = 0.1$, the dissipated wavenumber range is enlarged, where the additional damping is still small for very low wavenumbers. As $\Delta$ keeps increasing, these low wavenumber components can also be damped by the SFD method. This clearly explains why adding the SFD treatment with large filter width will reduce the overall error.

\begin{figure}[htbp]
	\includegraphics[width=200pt]{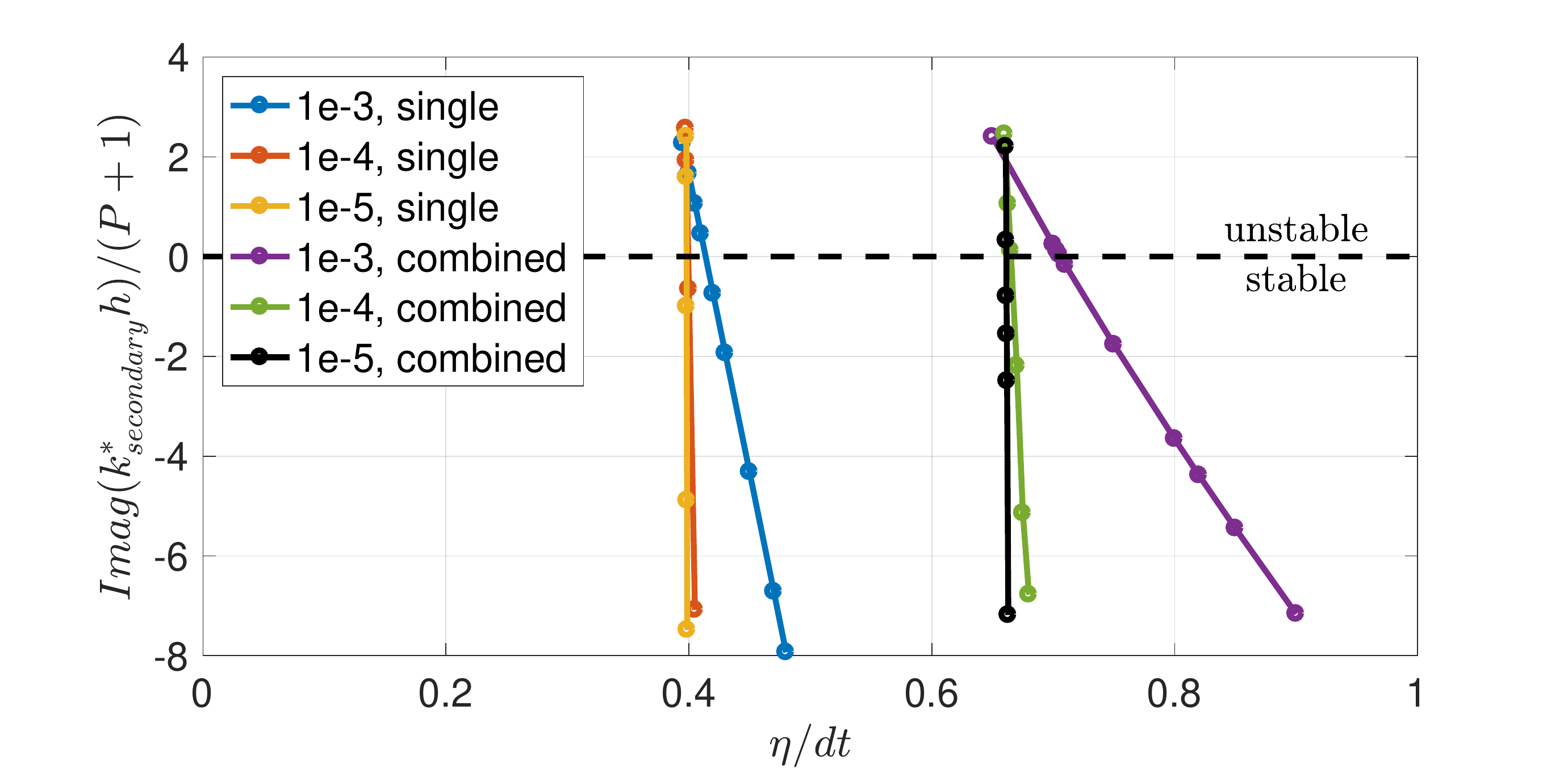}
	\centering
	\caption{Dissipation versus critical time step (stability limit) under different penalization parameters and control coefficients for a given $ \Delta t$. \textit{single} refers to only volume penalization or only the SFD method, while \textit{combined} refers to the combination of volume penalization and the SFD method with $\chi_f = 1/\IBMparam$. The filter width is sufficiently large and is set to $\Delta = 100$. Positive dissipation indicates numerical instabilities.}
	\label{fig:advection-stab}
\end{figure}

Secondly, the fully-discrete analysis is further used to obtain guidelines for the selection of penalization and control parameters $\IBMparam$ and $\chi_f$. Fully-discrete analysis considers both time and space discretisation and is able to evaluate the stability of the space-time scheme. For time discretisation, we choose a standard explicit third-order Runge-Kutta scheme with the following fully-discrete operator matrix \citep{vermeire2017behaviour,he2020dispersion}
\begin{equation}
    \boldsymbol{A} = \boldsymbol{I} + \Delta t \boldsymbol{M} + \frac{1}{2} (\Delta t \boldsymbol{M})^2 + \frac{1}{6} (\Delta t \boldsymbol{M})^3,
\end{equation}
where $M$ is the global matrix in Equation \ref{eq:global} for space discretisation. In this analysis the splitting scheme for both the volume penalization and the SFD is not considered, thus resulting a more strict requirement for stability. We consider the same discretisation used for the semi-discrete analysis, i.e., $N = 40$, $r = 1/40$ and $P = 3$. Through fully-discrete analysis, we will look at how the change in penalization parameter and in SFD parameters affects the stability of the space-time scheme. The stability of the system is reflected by the dissipation curve, where instability is characterised by positive dissipation. Since small $\IBMparam$ introduces stiff source terms that lead to instabilities \citep{schneider2015immersed}, the critical parameters for volume penalization and the SFD method, as well as the combination of the two, need to be carefully studied. As discussed in \citep{kou2021VonNeumann}, when IBM is considered in eigensolution analysis, several solid modes among secondary modes exist, which have constant dispersion and dissipation across all wavenumbers. The dissipation of these solid modes will move towards positive that leads to instabilities. 

We start from the guideline for volume penalization, i.e., $\IBMparam = \Delta t$ \cite{kolomenskiy2009fourier,engels2015numerical}. By fixing the time step $\Delta t$, we gradually change $\IBMparam$ and look at the change of dissipation in the solid modes. In addition, we also consider replacing the volume penalization with the SFD method, and changing the control coefficient $\chi_f$ instead, with a large filter width. The combined scheme of volume penalization and the SFD method with $\chi_f = 1/\IBMparam$ is finally tested to check the critical parameters. A typical result of the combined method is shown in Figure \ref{fig:advection-full}, where dispersion and dissipation behaviors at $\Delta = 1 \times 10^{-3}$ and different $\chi_f = 1/\IBMparam$ are shown. From these curves, it can be observed that changing $\IBMparam$ and $\chi_f$ around the critical values does not lead to big change on the primary physical mode. However, dissipation of the secondary solid mode moves towards the positive plane, where $\IBMparam$ is below $0.703 \Delta t = 0.703 \times 10^{-3}$. Therefore, when combining volume penalization and the SFD method, the empirical guideline $\IBMparam = \Delta t$ still guarantees the stability. 

\begin{figure}[htbp]
	\begin{subfigure}{.45\textwidth}
		\includegraphics[width=200pt]{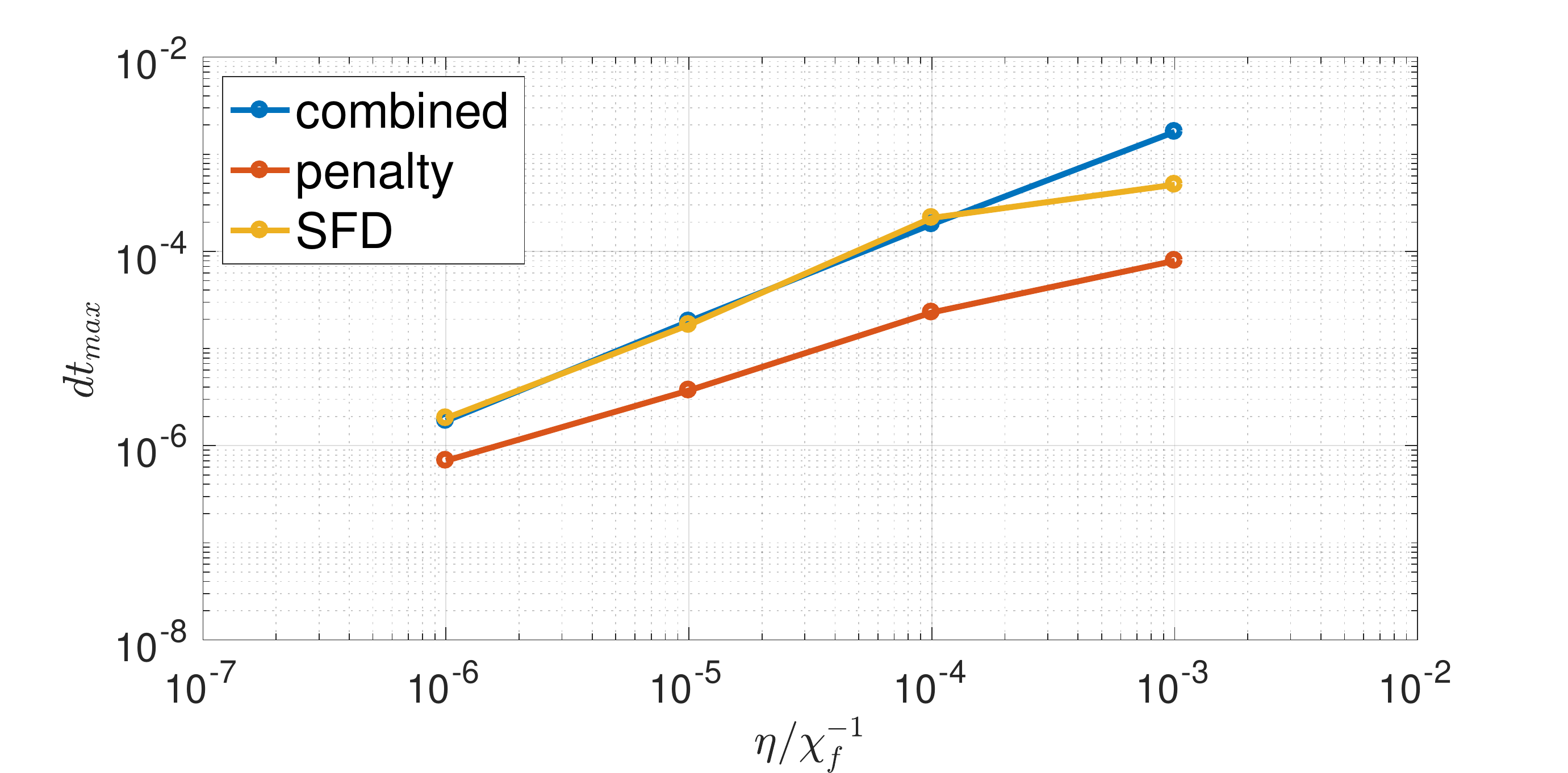}
		\caption{}
	\end{subfigure}
	\begin{subfigure}{.45\textwidth}
		\includegraphics[width=200pt]{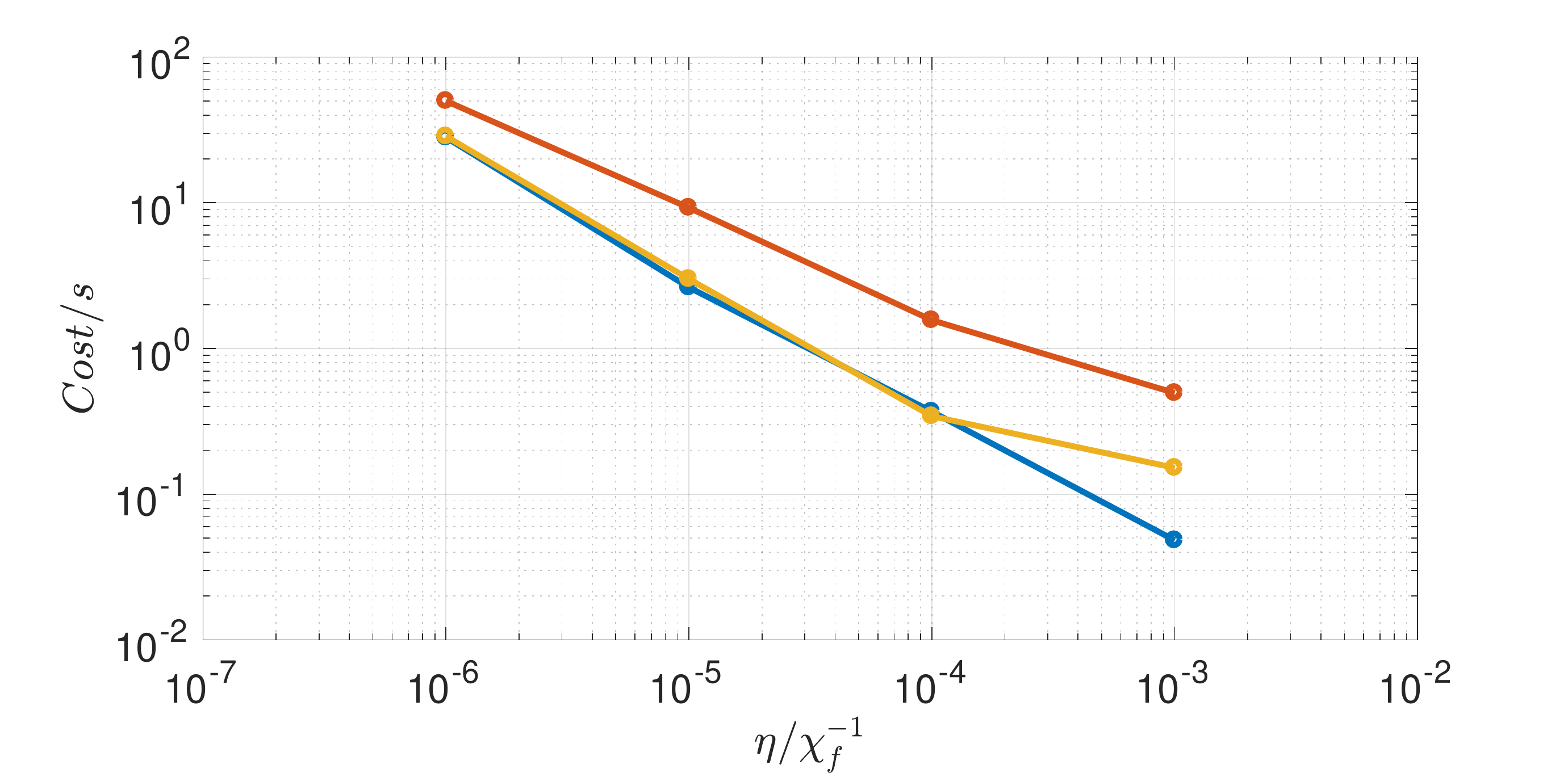}
		\caption{}
	\end{subfigure}
	\centering
	\caption{Comparison of same accuracy with different methods. (a) Maximal time step allowed. (b) Computational cost.}
	\label{fig:advection-crit}
\end{figure}

To further investigate the stability behavior and obtain more general relationship between critical $\IBMparam$ ($\chi_f$) and $\Delta t$, additional fully-discrete analyses are performed. The dissipation of the solid mode close to zero is shown in Figure \ref{fig:advection-stab}. We consider three fixed time steps $1 \times 10^{-3}$, $1 \times 10^{-4}$, and $1 \times 10^{-5}$. For a given time step, we have three different conditions: only volume penalization or only the SFD method, and the combined method with $\chi_f = 1/\IBMparam$. We notice that when only volume penalization or only the SFD method is used, the dissipation is almost the same, therefore only one curve is shown and noted by 'single'. From Figure \ref{fig:advection-stab}, when only one of the methods is considered, the critical parameter $\IBMparam$ or $ 1 / \chi_f$ lies between $0.4 \Delta t$ and $0.42 \Delta t$. If the combined method is used, the critical parameter $\chi_f = 1/\IBMparam$ lies between $0.65 \Delta t$ and $0.71 \Delta t$. When compared to a single method, the stability limit becomes more strict, but improved accuracy can be reached. Furthermore, these analyses still show the benefit of using the combined method, since the general guideline for selecting the penalization parameter (also control parameter of the SFD method since $\chi_f = 1/\IBMparam$) $\IBMparam = \Delta t$ works well while the accuracy is improved. In addition, as the time step decreases, the ratio between critical parameter and time step also decreases, indicating that a smaller $\IBMparam$ or a larger $\chi_f$ can be used. Note that the present fully-discrete analyses consider the direct imposition of the penalization source term, therefore the criterion is more strict that guarantees the generality. If the Strang splitting approach \citep{strang1968construction} is used, more relaxed time step restriction will be obtained.

Finally the computational cost is investigated. For each group of case, we first set the parameters for the combined method with $\chi_f = 1/\IBMparam$ and set $\Delta$ to 100. We simulate the advection equation to final time $1.1$, and save the error in the flow and the maximal allowable time step. Based on this case, we study if  only a single method is used, how much penalty $\IBMparam$ or $\chi_f$, and maximal time step are required to reach the same error. The comparison is shown in Figure \ref{fig:advection-crit}a. It can be concluded that for the same error, the combined method and only SFD method allow similar large time step. However, the maximal time step for the volume penalization method is more limited, indicating increased computational cost, compared in \ref{fig:advection-crit}b. This also shows the advantage of the proposed SFD method (both only SFD and the combined scheme), since in the present simulation, it has a lower error compared with the volume penalization. It should be reminded that the SFD method is only applied to the solution points inside the solid, thus a low computational overhead is included.

\subsection{Flow past a NACA0012 airfoil}
\label{sec:naca}
\begin{figure*}[htbp]
    \begin{subfigure}{.4\textwidth}
		\includegraphics[width=180pt]{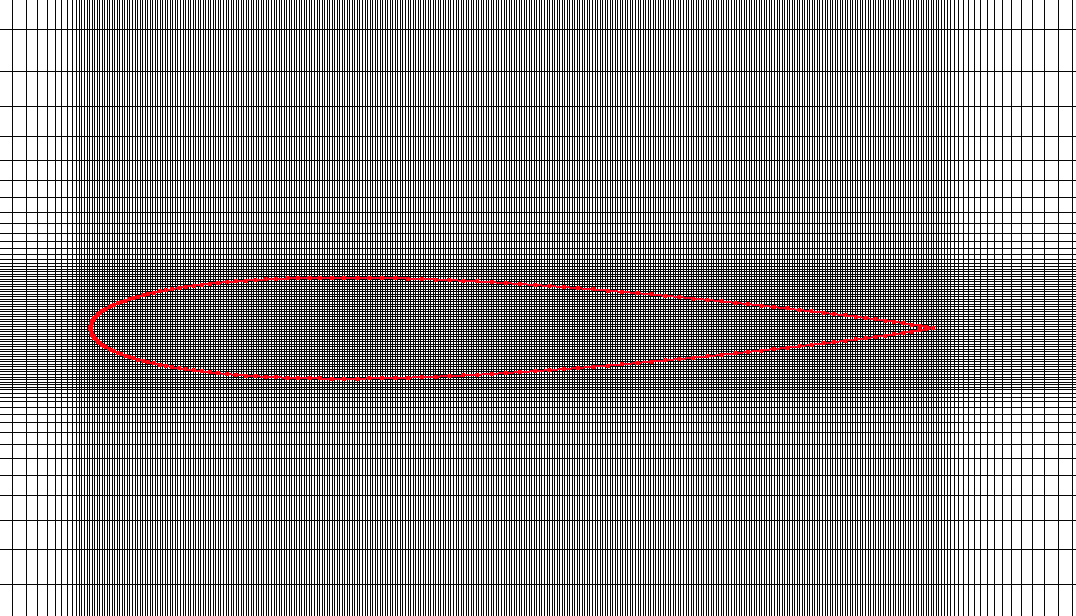}
		\caption{}
	\end{subfigure}
	\begin{subfigure}{.4\textwidth}
		\includegraphics[width=150pt]{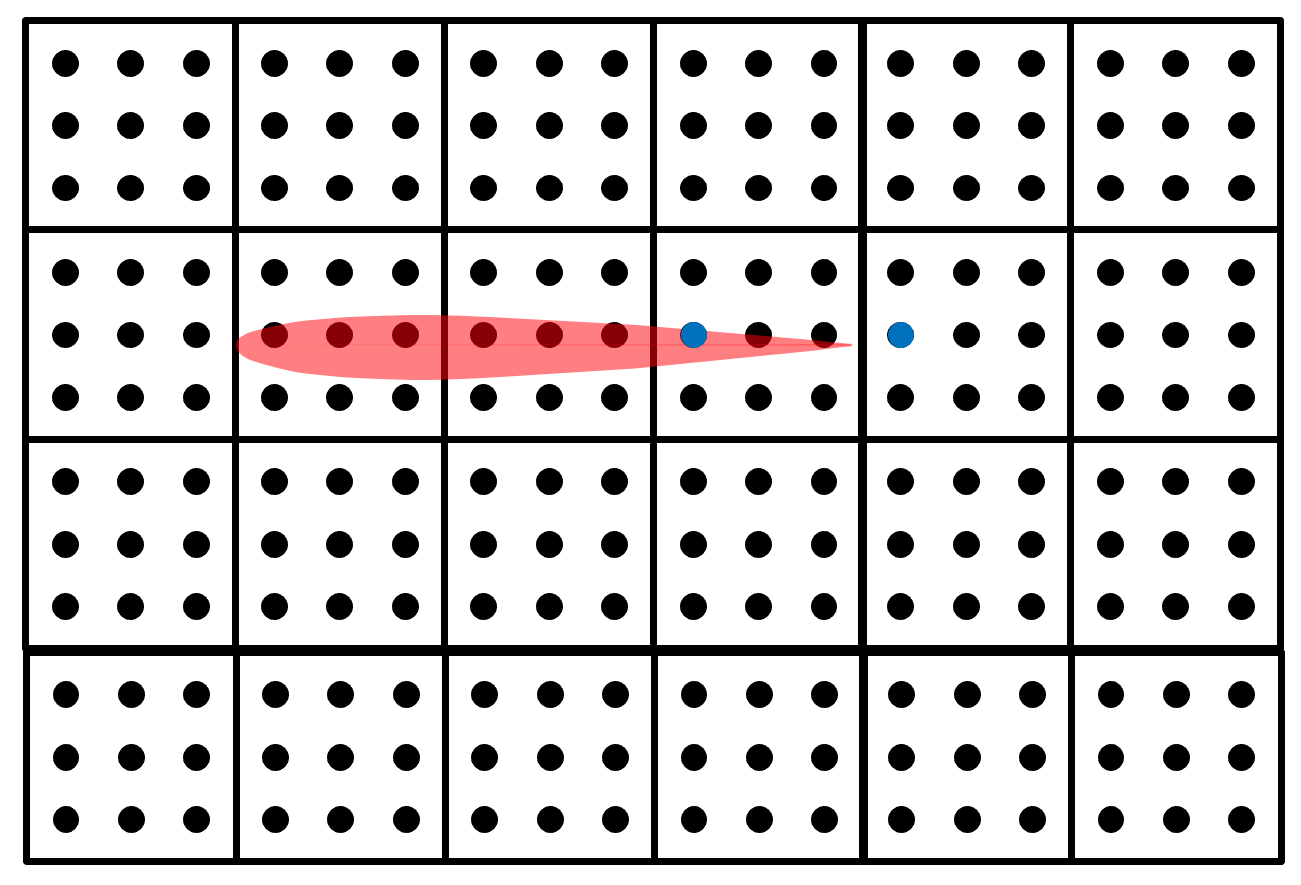}
		\caption{}
	\end{subfigure}
	\centering
	\caption{A typical IBM mesh and sketch of the probe. a) The computational grid with size $0.003$ near the wall. b) Sketch of probe positions (this figure is only an illustration and is not scaled with the real mesh).}
	\label{fig:naca}
\end{figure*}

In this section, we apply the SFD method to impose the no-slip wall boundary condition for the Navier-Stokes equations. In the following simulations, the workflow in Figure \ref{fig:workflow} is used. We consider flow over a NACA0012 airfoil with Mach number 0.5, Reynolds number 5000 (based on the chord length $c$ and free stream velocity) and angle of attack zero. Since we are more interested in the accuracy in the imposition of boundary conditions, we select two probe points in the flow field to observe the temporal evolution of flow velocities. A typical mesh and sketch of the probes are shown in Figure \ref{fig:naca}, where one point is in the solid and the other is in the flow. The airfoil is positioned in $x/c \in [-0.5,0.5]$. The probes are selected as the closest Gaussian solution points to the particular coordinates $(x, y)/c = (0.0225, 0.015)$ for the inside point and $(x, y)/c = (0.555, 0.015)$ for the outside point.

\begin{figure*}[htbp]
    \begin{subfigure}{.33\textwidth}
		\includegraphics[width=150pt]{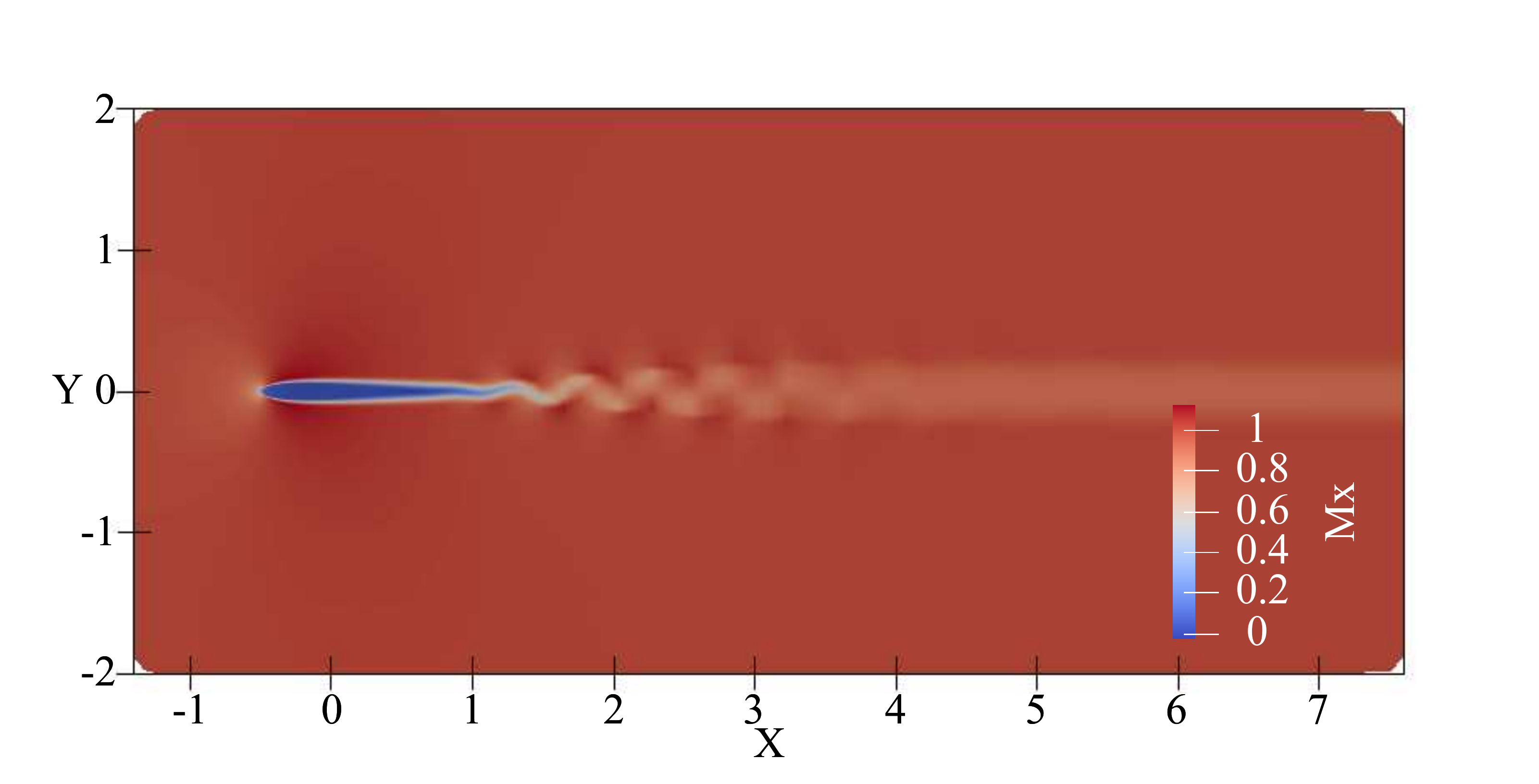}
		\caption{}
	\end{subfigure}
	\begin{subfigure}{.33\textwidth}
		\includegraphics[width=150pt]{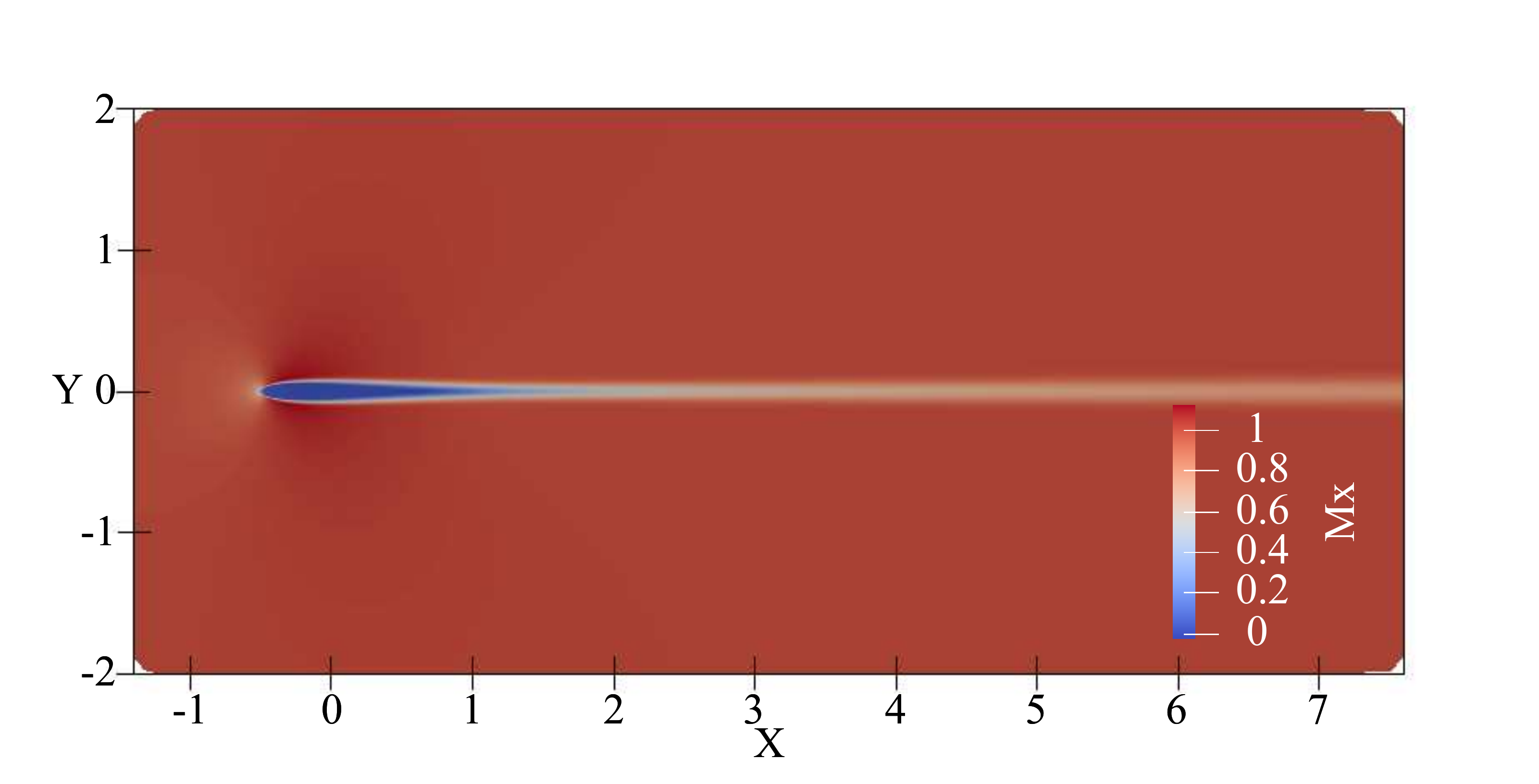}
		\caption{}
	\end{subfigure}
	\begin{subfigure}{.33\textwidth}
		\includegraphics[width=150pt]{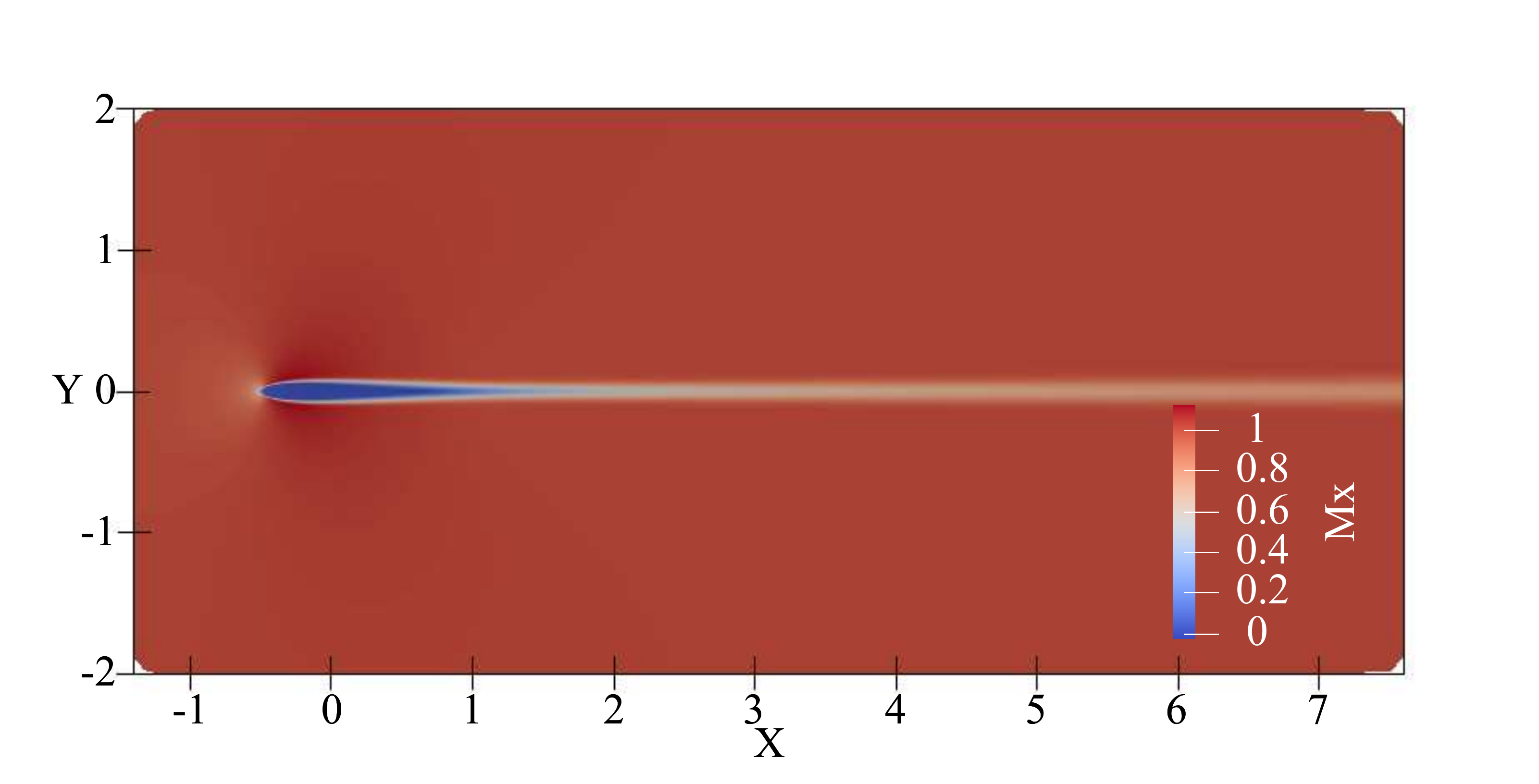}
		\caption{}
	\end{subfigure}
	\centering
	\caption{Flow field comparison: a) Volume penalization with $\IBMparam = 1 \times 10^{-2}$. b) Volume penalization with $\IBMparam = 5 \times 10^{-3}$. c) The combined method with $\IBMparam = 1 \times 10^{-2}, \chi_f = 5 \times 10^{3}, \Delta = 10$.}
	\label{fig:naca-flow}
\end{figure*}

In the first test case, a relatively coarse grid is used to illustrate the effect of volume penalization and the SFD parameters. We generate a uniform grid with size $0.01$ in the domain $x/c \in [-0.51, 0.51]$ and $y/c \in [-0.07, 0.07]$, and gradually increase the mesh size to the far field. The computation domain is defined $x/c \in [-20, 40]$ and $y/c \in [-20, 20]$, which leads to a total number of element $174 \times 68$. We simulate different flows under various combinations of parameters. The order of polynomial is set to $P=2$ and the time step is set to $5 \times 10^{-4}$ given the CFL condition. A comparison of flow fields is shown in Figure \ref{fig:naca-flow}, where two cases from volume penalization and one case from the combined method are shown. For the present flow condition, the flow is steady \cite{swanson2016comparison}. However, when the volume penalization with $\IBMparam = 1 \times 10^{-2}$ is used, the flow becomes unsteady. This indicates that when the penalization term is not sufficiently large (using volume penalization), the flow physics can change and nonphysical oscillations develop. However, meaningful flow fields can be obtained if we keep reducing the penalization parameter, or introduce the SFD method as an additional source term to impose the boundary condition. The effects of these two options are shown in Figure \ref{fig:naca-flow}b and Figure \ref{fig:naca-flow}c, respectively. 

\begin{figure*}[htbp]
    \begin{subfigure}{.48\textwidth}
		\includegraphics[width=200pt]{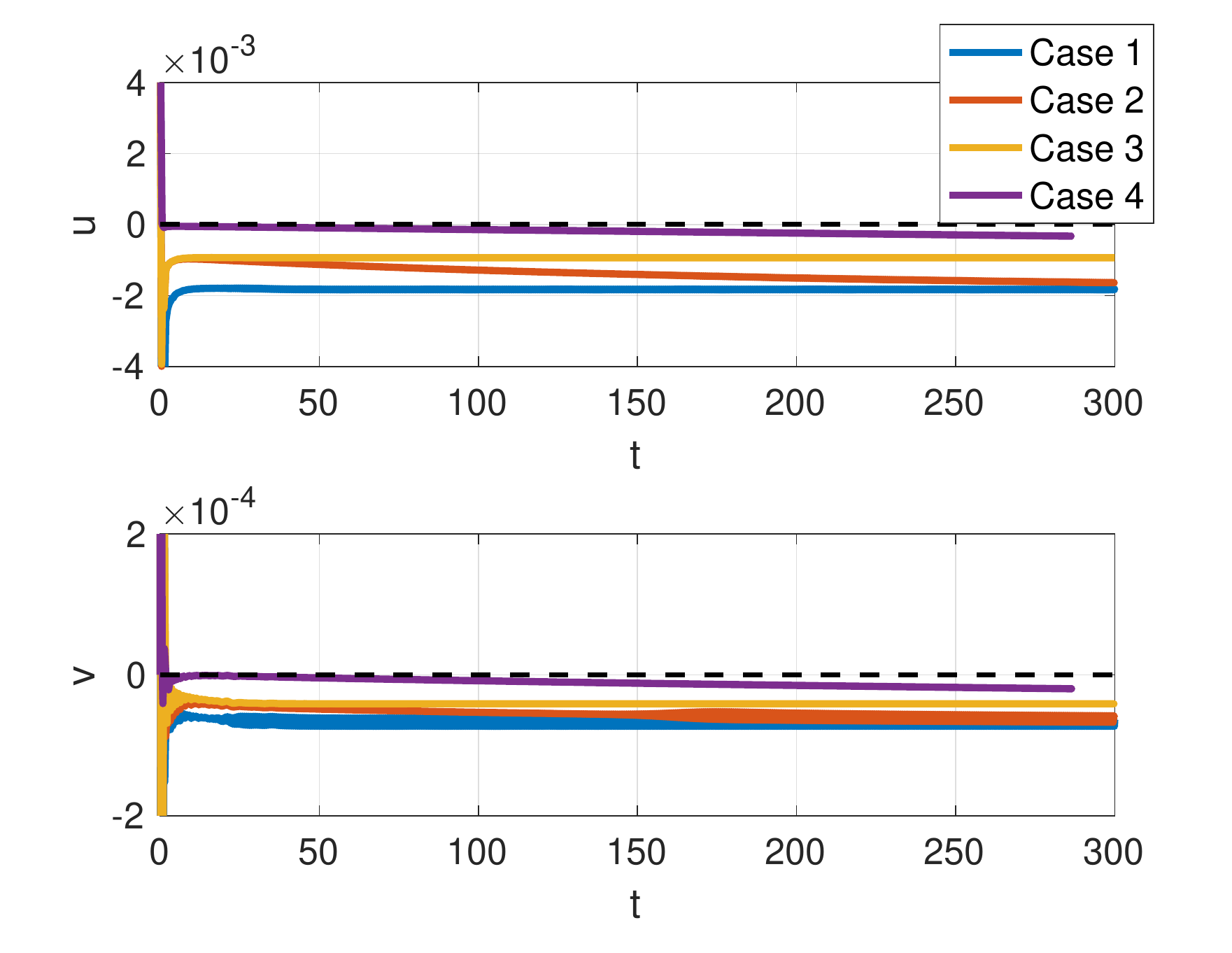}
		\caption{}
	\end{subfigure}
	\begin{subfigure}{.48\textwidth}
		\includegraphics[width=200pt]{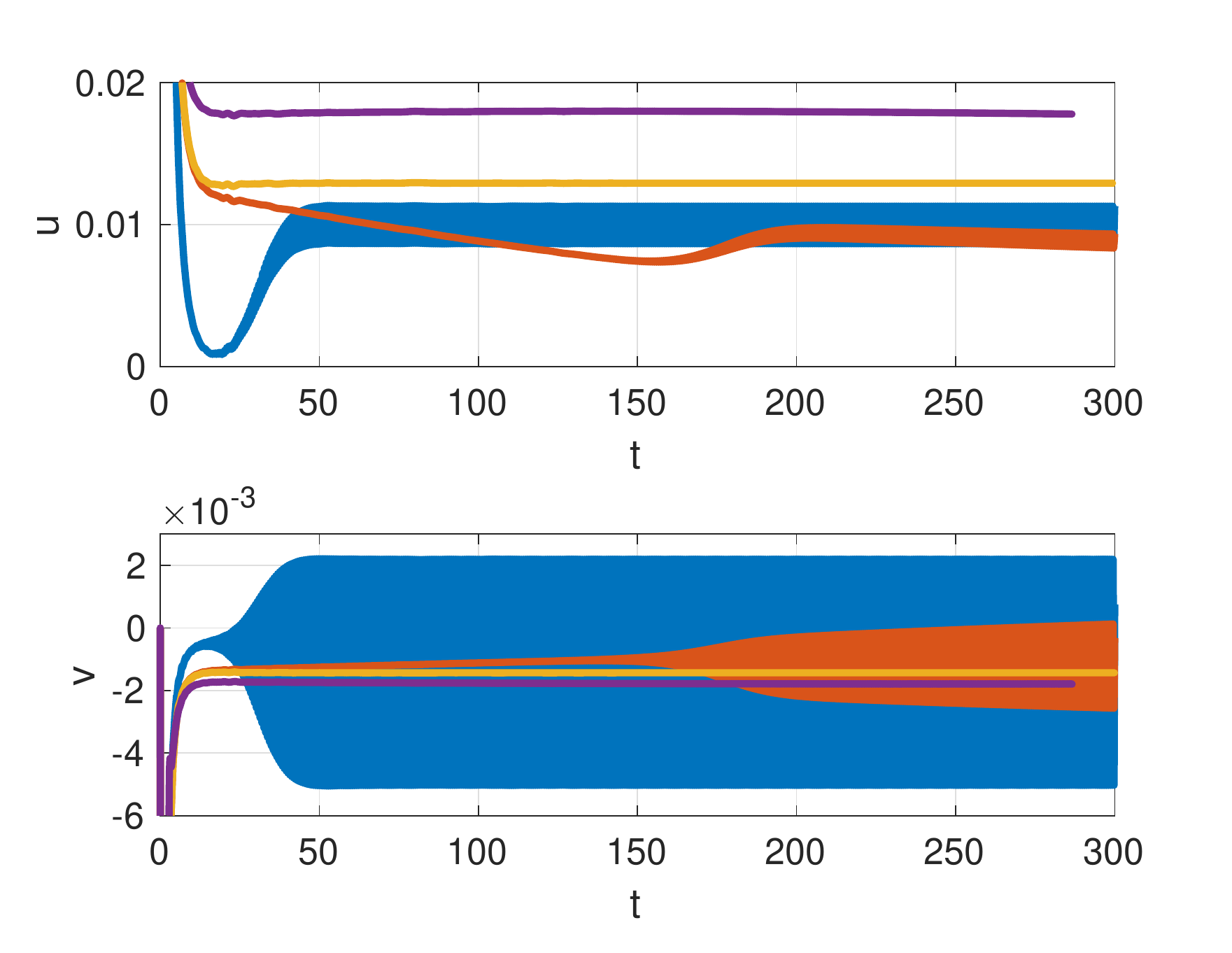}
		\caption{}
	\end{subfigure}
	\centering
	\caption{Evolution of flow velocities at two probe points. a) Inside point. b) Outside point. Case 1: $\IBMparam = 1 \times 10^{-2}$. Case 2: $\IBMparam = 1 \times 10^{-2}, \chi_f = 1 \times 10^{2}, \Delta = 100$. Case 3:$\IBMparam = 5 \times 10^{-3}$. Case 4: $\IBMparam = 5 \times 10^{-3}, \chi_f = 2000, \Delta = 100$.}
	\label{fig:naca-probex}
\end{figure*}

\begin{figure*}[htbp]
    \begin{subfigure}{.48\textwidth}
		\includegraphics[width=200pt]{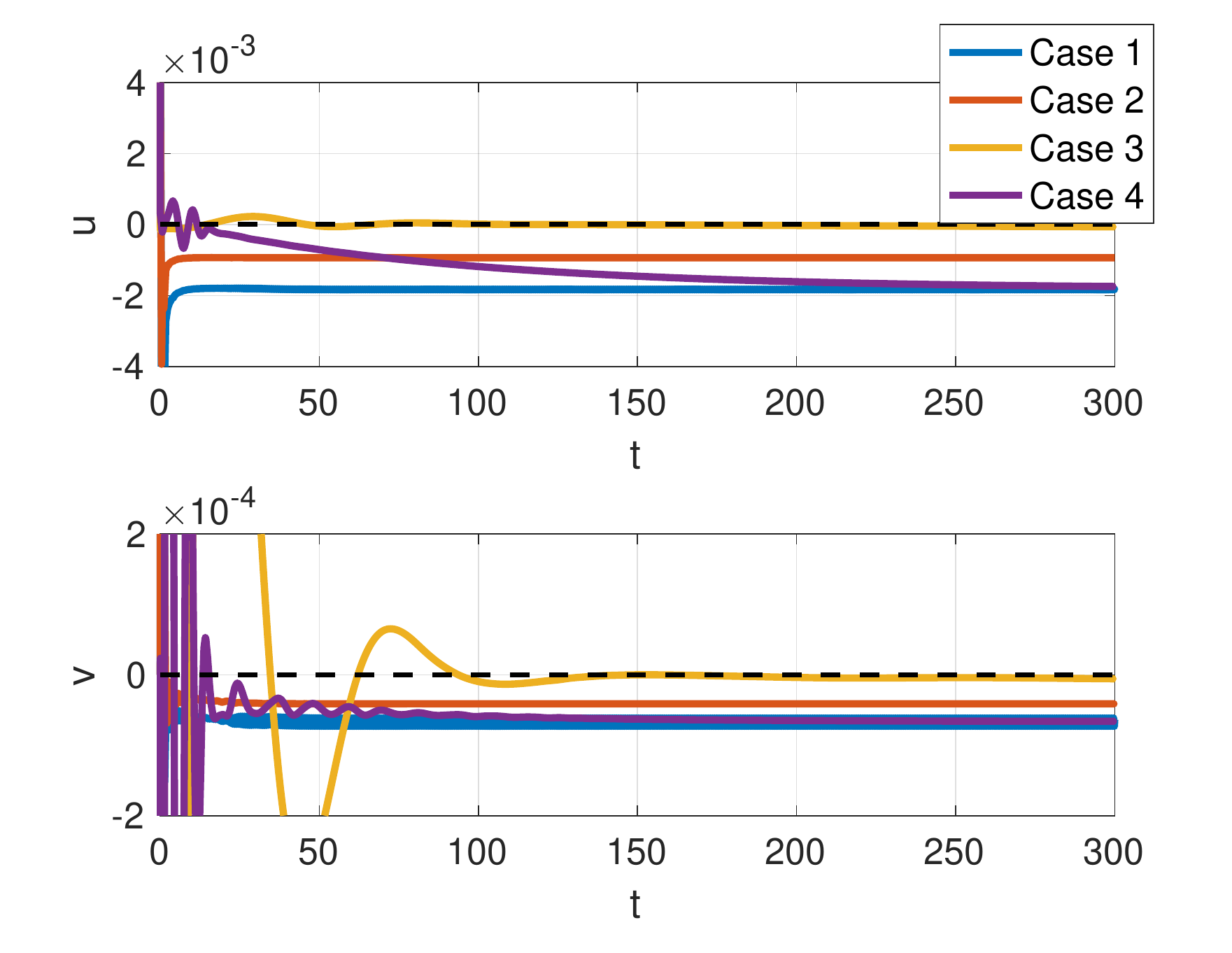}
		\caption{}
	\end{subfigure}
	\begin{subfigure}{.48\textwidth}
		\includegraphics[width=200pt]{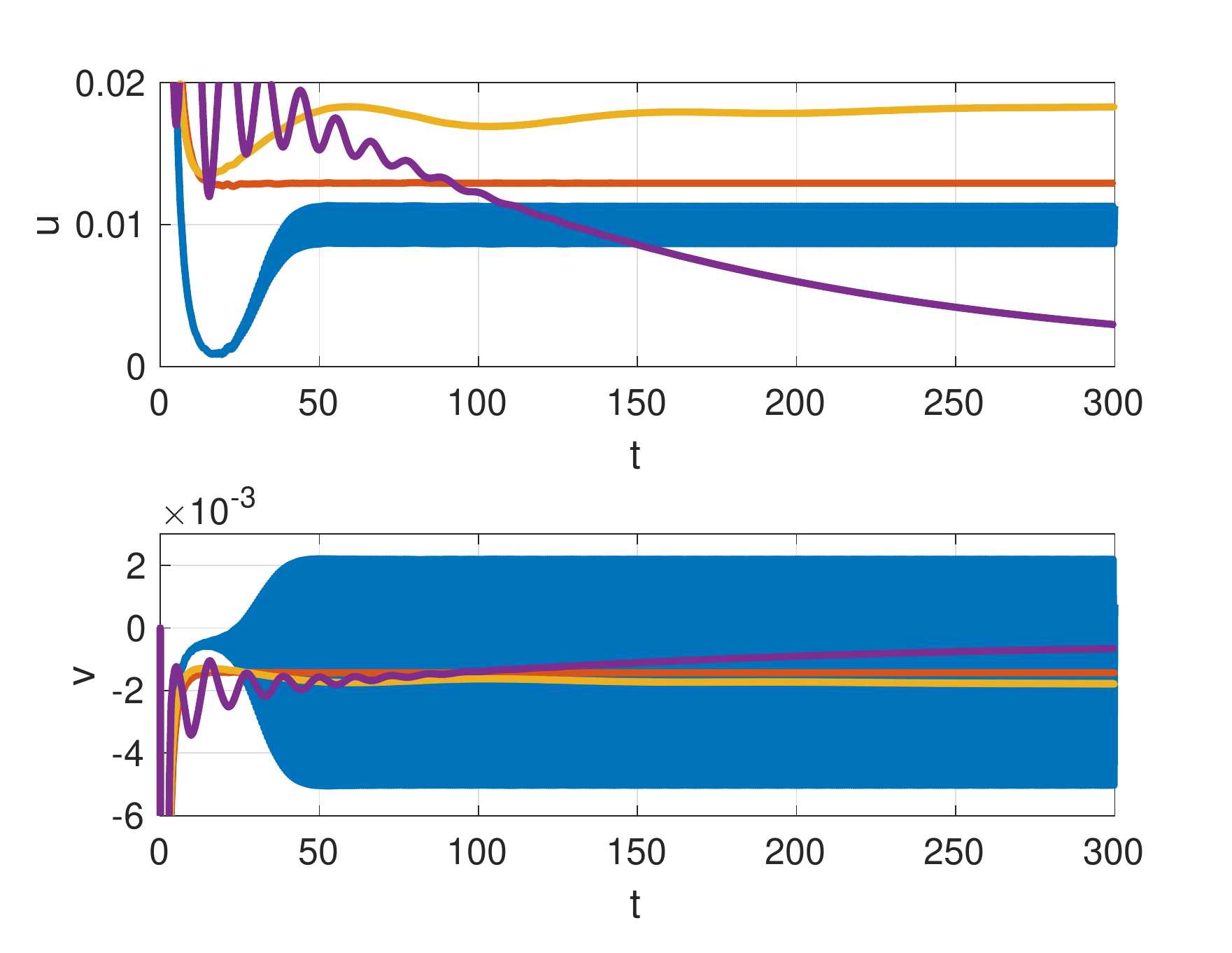}
		\caption{}
	\end{subfigure}
	\centering
	\caption{Evolution of flow velocities at two probe points. a) Inside point. b) Outside point. Case 1: $\IBMparam = 1 \times 10^{-2}$. Case 2: $\IBMparam = 5 \times 10^{-3}$. Case 3:$\IBMparam = 1 \times 10^{-2}, \chi_f = 5\times 10^{3}, \Delta = 100$. Case 4: $\IBMparam = 1\times 10^{-2}, \chi_f = 5\times 10^{3}, \Delta = 2$.}
	\label{fig:naca-probe}
\end{figure*}

Temporal evolutions of flow velocities at two probes inside and outside the body are compared in Figure \ref{fig:naca-probex} and Figure \ref{fig:naca-probe}. Note that inside the solid, the expected horizontal velocity $u$ and vertical velocity $v$ are zero because of the no-slip wall boundary condition. First, comparisons between volume penalization and the combined method are shown in Figure \ref{fig:naca-probex}, where two penalization parameters $1 \times 10^{-2}$ and $5 \times 10^{-3}$ are set, and $\chi_f = 1 / \IBMparam$ is chosen for the combined method. As shown in the figure, when volume penalization with $1 \times 10^{-2}$ is used (Case 1, the flow field is shown in Figure \ref{fig:naca-flow}a), there is a lack of penalization, which leads to oscillations of flow quantities inside the solid. These oscillations are magnified and will result in strong oscillations in the flow field near the trailing edge, leading to an unphysical flow field. When the SFD method is added for the same volume penalization method, the oscillations are reduced. However, since the control parameter in the SFD method is not big enough, weaker oscillations still exist. When the volume penalization method with $5 \times 10^{-3}$ is used, these oscillations are removed, therefore the resulting outer flow becomes steady and physically meaningful. In addition, the improved accuracy of the combined method is also seen in Case 4 of Figure \ref{fig:naca-probex}, where the simulated velocities are the closest to zero. These results show the advantages of adding the SFD method to further reduce the spurious oscillations, and improve the accuracy of the immersed boundaries. 

In addition, another group of results is compared in Figure \ref{fig:naca-probe} to show the capability of the SFD method to improve the accuracy when the penalization source term is small for volume penalization. We add the SFD method with $\chi_f = 5 \times 10^{3}$ based on volume penalization method with $\IBMparam = 1 \times 10^{-2}$ (Case 3 and Case 4). Clearly, when the SFD treatment is included, the oscillations coming from volume penalization are removed and the accuracy is greatly improved. A comparison between Case 3 and Case 4 also shows the effect of the filter width $\Delta$. When $\Delta$ is set to 2, the oscillation from volume penalization with $\IBMparam = 1 \times 10^{-2}$ is removed. However, as shown in Figure \ref{fig:naca-probe}a, constant shift at lower frequencies still exists, making the predicted velocities similar to those computed using volume penalization. This highlights the importance of using a large $\Delta$ to sufficiently filter the oscillations inside the solid. The combined method also leads to improved accuracy compared with volume penalization with a lower penalization parameter (Case 2). 

\begin{figure*}[htbp]
    \begin{subfigure}{.48\textwidth}
		\includegraphics[width=200pt]{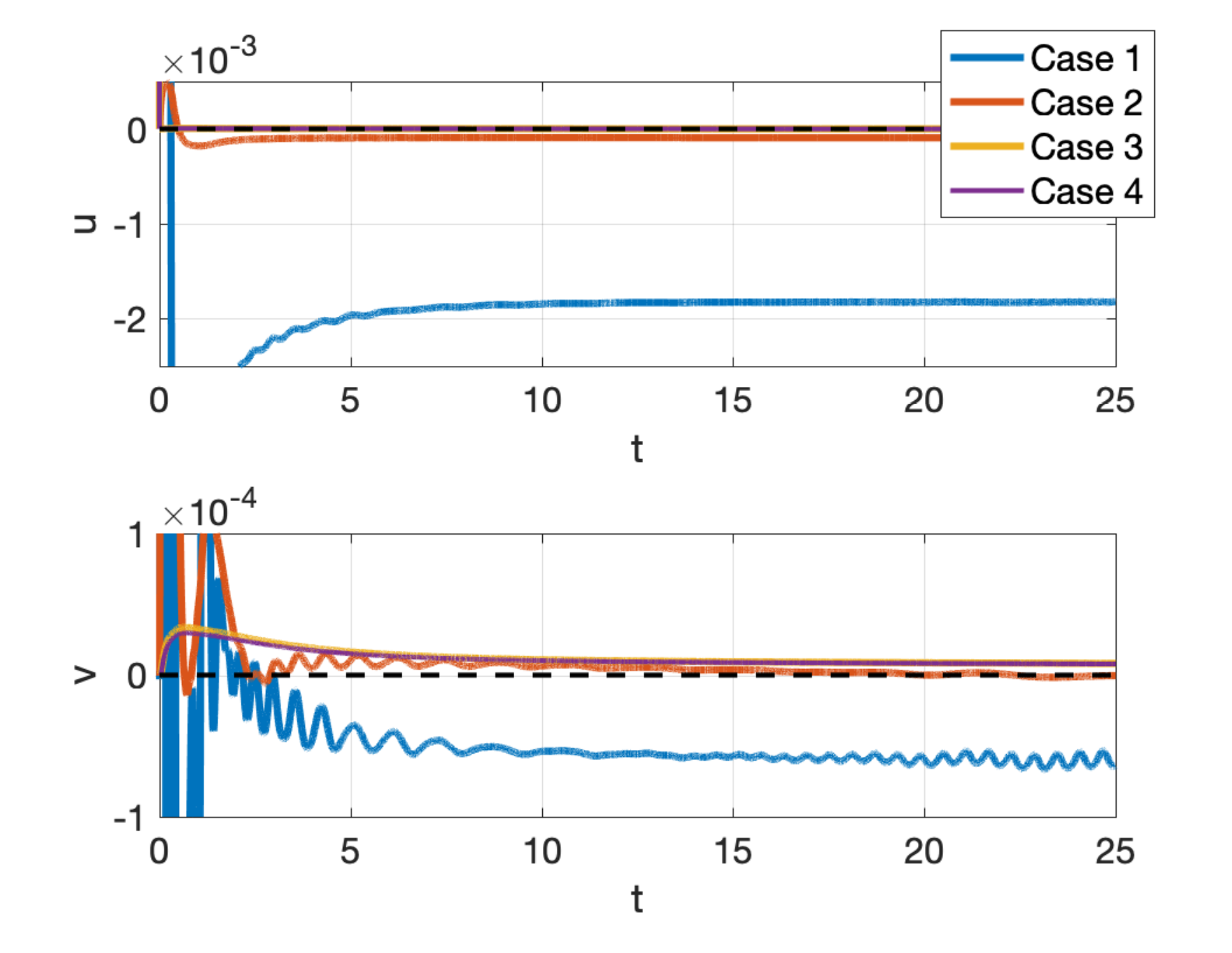}
		\caption{}
	\end{subfigure}
	\begin{subfigure}{.48\textwidth}
		\includegraphics[width=200pt]{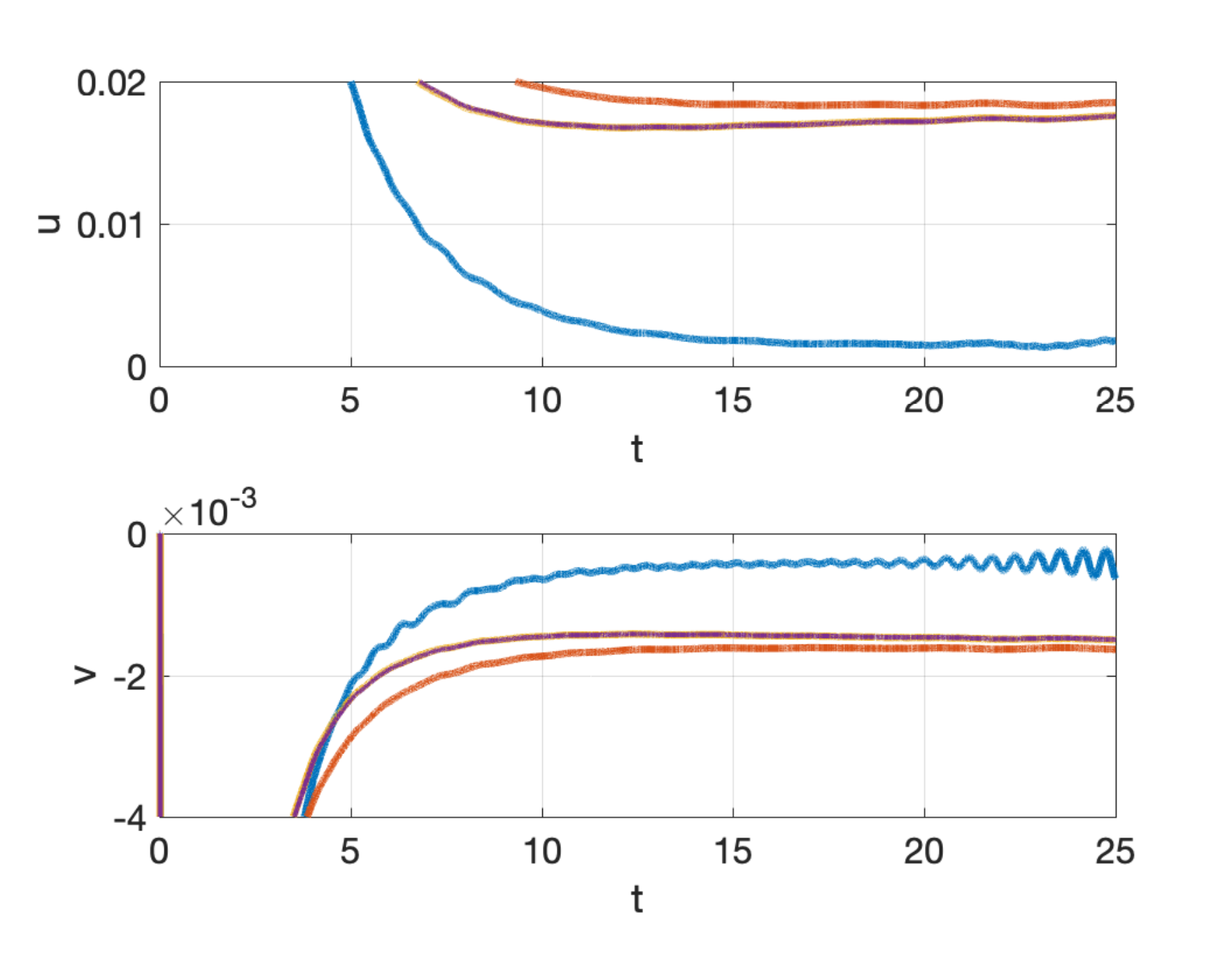}
		\caption{}
	\end{subfigure}
	\centering
	\caption{Evolution of flow velocities at two probe points. a) Inside point. b) Outside point. Cp, h = $3 \times 10^{-3}$. Case 1: $\IBMparam = 1\times 10^{-2}$. Case 2: $\IBMparam = 5\times 10^{-4}$. Case 3: $\IBMparam = \infty, \chi_f = 2\times 10^{4}, \Delta = 100$. Case 4: $\IBMparam = 5\times 10^{-4}, \chi_f = 2 \times 10^{4}, \Delta = 100$.}
	\label{fig:naca-probe2}
\end{figure*}

In the second test case, we analyze the effect of adding the SFD treatment for a refined mesh and improved resolution near the wall. We generate a locally uniform grid with size $0.003$, while the locally refined region and computational domain remain the same, and the total number of element is $453 \times 124$. We simulate different flows under various combinations of parameters. The order of polynomial is set to $P=2$ and the time step is set to $7.5 \times 10^{-5}$ given the CFL condition. Four test cases are shown in Figure \ref{fig:naca-probe2}, where the first two cases only use volume penalization with different $\IBMparam$. The third case only uses the SFD method without volume penalization ($\IBMparam \rightarrow \infty$), while the fourth case is the combined method with $\chi_f = 1 / \IBMparam$ and $\IBMparam = 5\times 10^{-4}$. Similarly to the previous case, when the penalization parameter is large (Case 1), oscillations exist inside the solid, leading to inaccuracy of the boundary condition and the oscillations in the flow region. Volume penalization with $\IBMparam = 5\times 10^{-4}$ (Case 2) helps to remove the oscillation and leads to more accurate boundary condition. When the same parameters for SFD method are considered (Case 3), smoother temporal responses are obtained, and improved accuracy is seen in the horizontal velocity. When the combined method is used (Case 4), small improvement on the accuracy is seen, indicating the solution is sufficiently accurate. Furthermore, for the present test case, the addition of the SFD method also results in faster convergence, as can been seen from the temporal response in Figure \ref{fig:naca-probe2}a.

\begin{figure*}[htbp]
    \begin{subfigure}{.48\textwidth}
		\includegraphics[width=200pt]{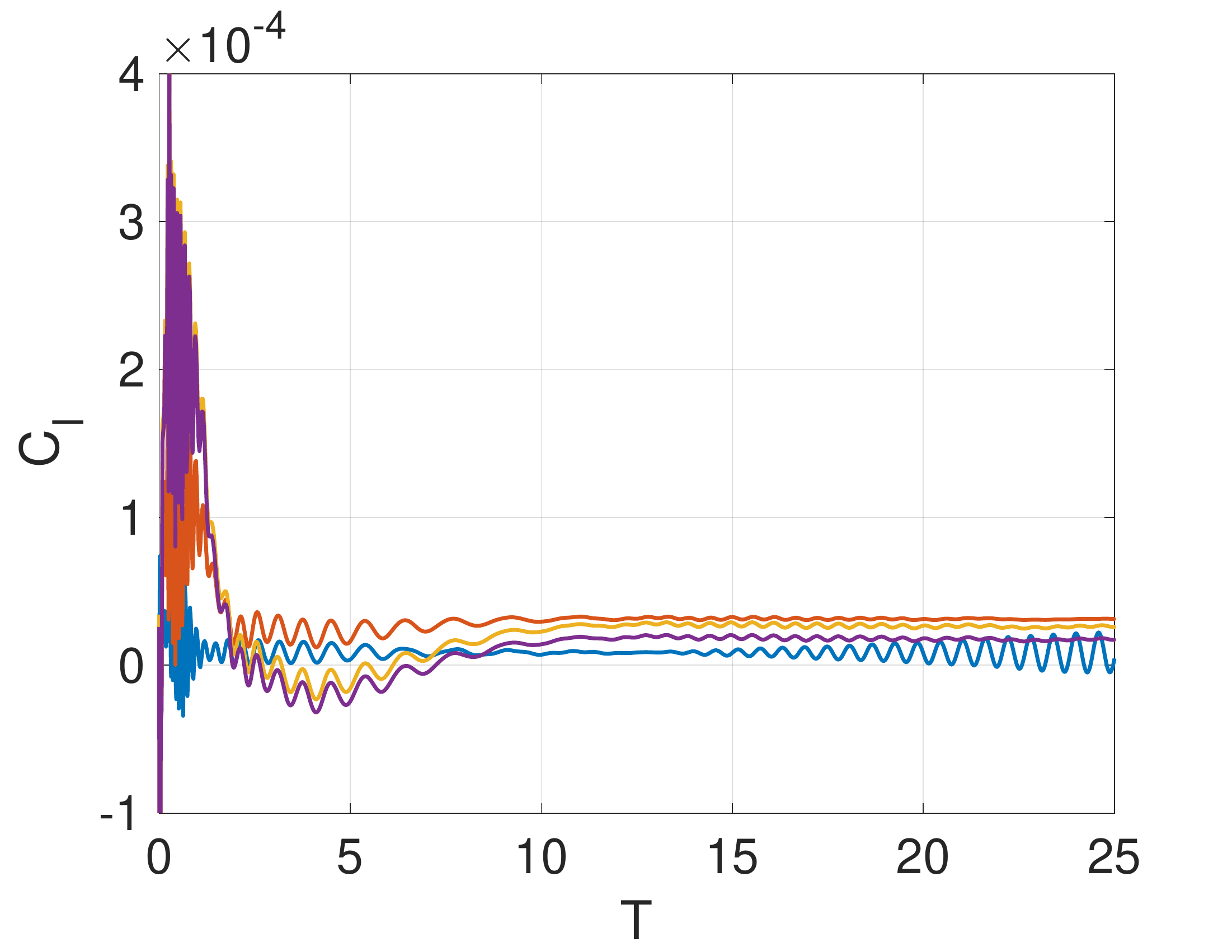}
		\caption{}
	\end{subfigure}
	\begin{subfigure}{.48\textwidth}
		\includegraphics[width=200pt]{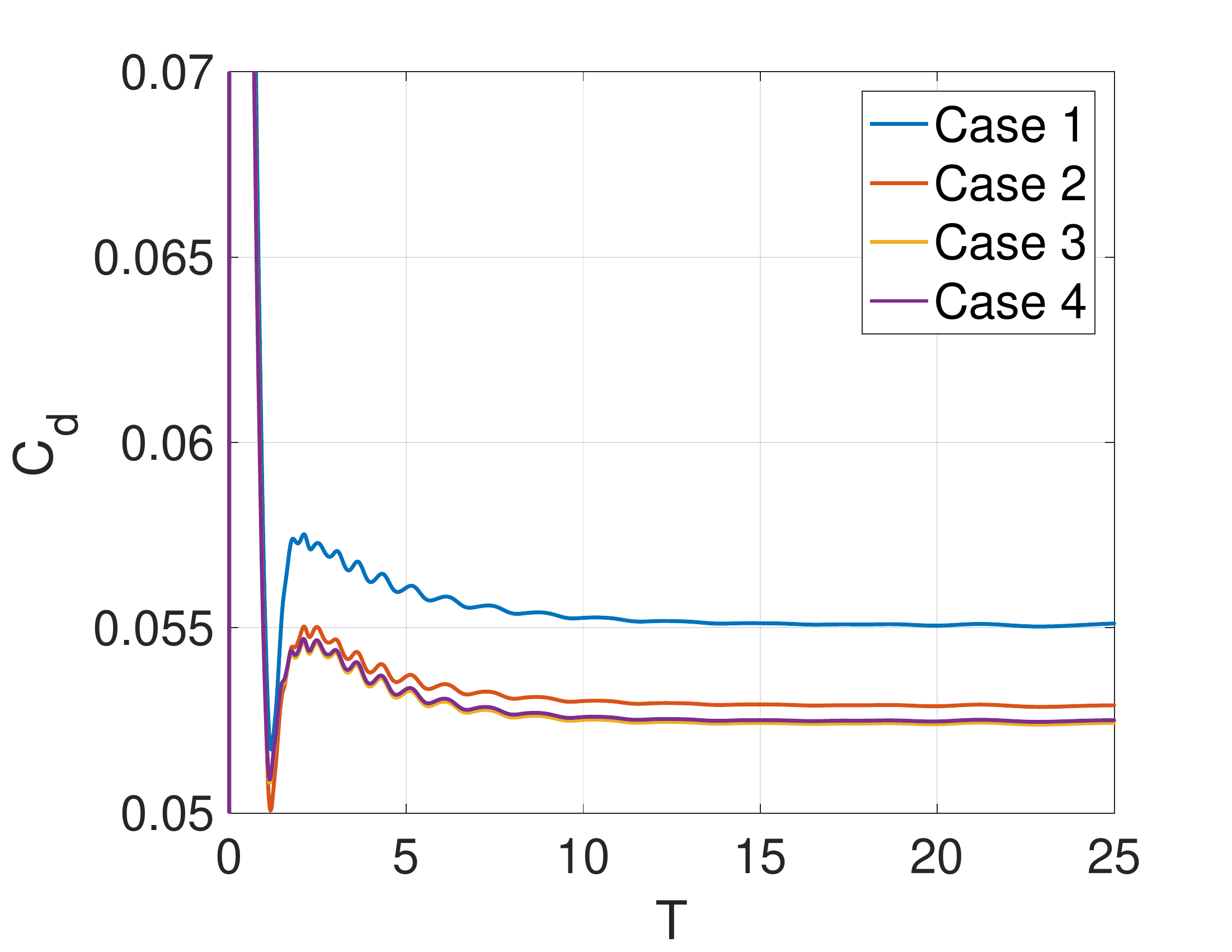}
		\caption{}
	\end{subfigure}
	\centering
	\caption{Evolution of lift and drag coefficients. a) Lift. b) Drag. Cp, h = $3 \times 10^{-3}$. Case 1: $\IBMparam = 1\times 10^{-2}$. Case 2: $\IBMparam = 5\times 10^{-4}$. Case 3: $\IBMparam = \infty, \chi_f = 2\times 10^{4}, \Delta = 100$. Case 4: $\IBMparam = 5\times 10^{-4}, \chi_f = 2\times 10^{4}, \Delta = 100$.}
	\label{fig:force}
\end{figure*}

\begin{figure*}[htbp]
	\begin{subfigure}{.48\textwidth}
		\includegraphics[width=200pt]{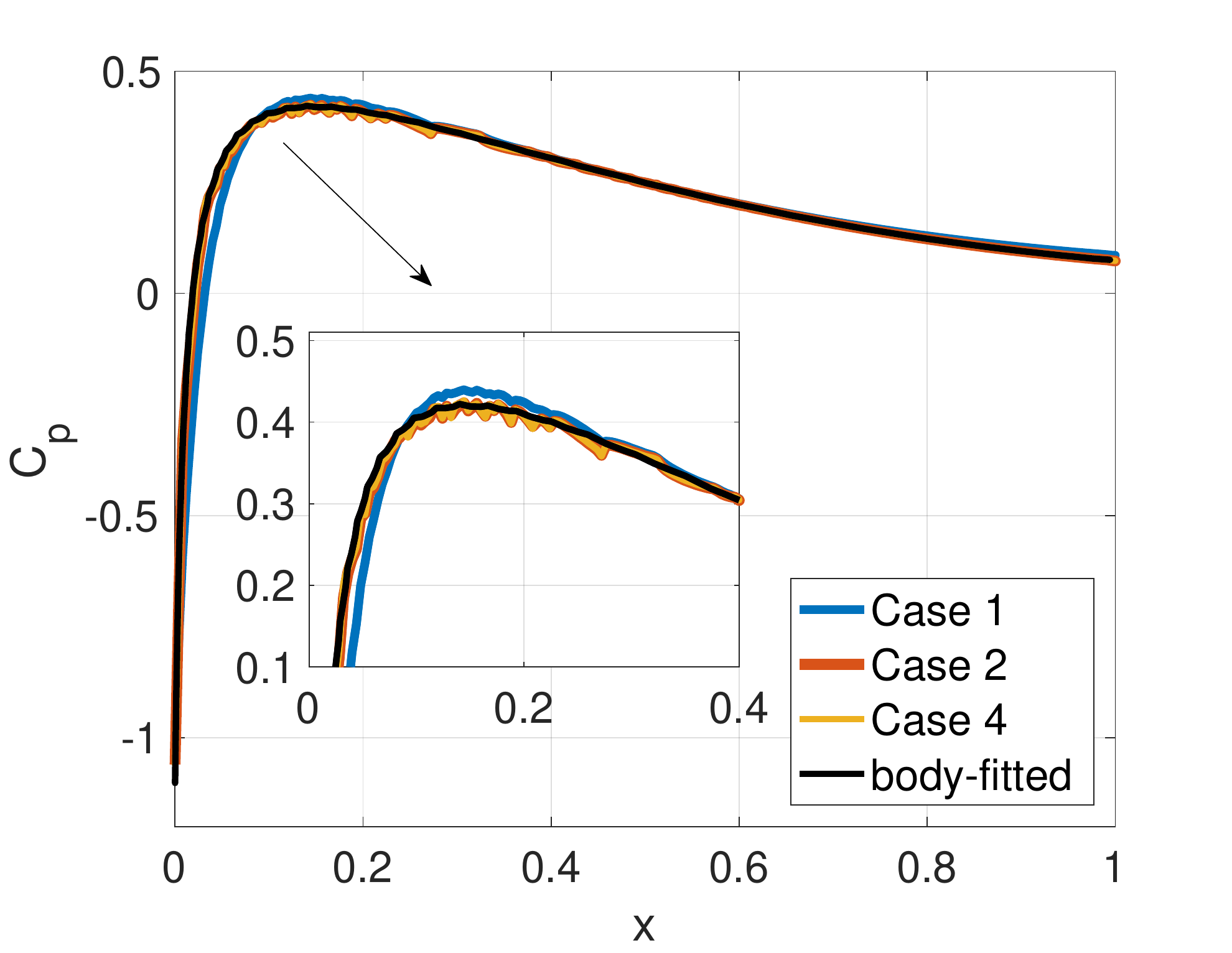}
		\caption{}
	\end{subfigure}
	\centering
	\caption{Cp, h = $3\times 10^{-3}$. Case 1: $\IBMparam = 1\times 10^{-2}$. Case 2: $\IBMparam = 5\times 10^{-4}$. Case 3: $\IBMparam = \infty, \chi_f = 2\times 10^{4}, \Delta = 100$.}
	\label{fig:naca-Cp}
\end{figure*}

The integrated aerodynamic forces are compared in Figure \ref{fig:force} to further justify the use of the SFD method. Figure \ref{fig:force} shows that the oscillations in the flow field, induced by the inaccuracy of the volume penalization with $\IBMparam = 1\times 10^{-2}$, lead to oscillations in the aerodynamic forces. This highlights the importance of the correct imposition of the boundary condition. It should be noted that in Figure \ref{fig:force}a, the expected lift coefficient should be zero because of the symmetric NACA0012 geometry at zero angle of attack. From this perspective, the results are consistent with previous sections, since the SFD method (Case 3) outperforms the volume penalization method (Case 2) when the source terms are at the same level $\chi_f = 1 / \IBMparam$. In addition, the combined method will lead to improved accuracy, showing the advantage of using our new method. The difference between Case 1 and other cases in predicting the drag coefficient also highlights the need of accurate boundary conditions in IBM treatment. Accuracy of the proposed approach is also shown in the comparison of the  surface pressure coefficient in Figure \ref{fig:naca-Cp}. This prediction can be further enhanced if the near-wall resolution is improved by local hp-refinement, as shown in \citep{kou2021IBMFR1}. In summary, the SFD treatment and the combined method are alternatives to the traditional penalty-based IBM approach, and show improved accuracy to compute steady aerodynamic forces. 

\subsection{Flow past a cylinder}
\begin{figure*}[htbp]
	\begin{subfigure}{.4\textwidth}
	\includegraphics[width=180pt]{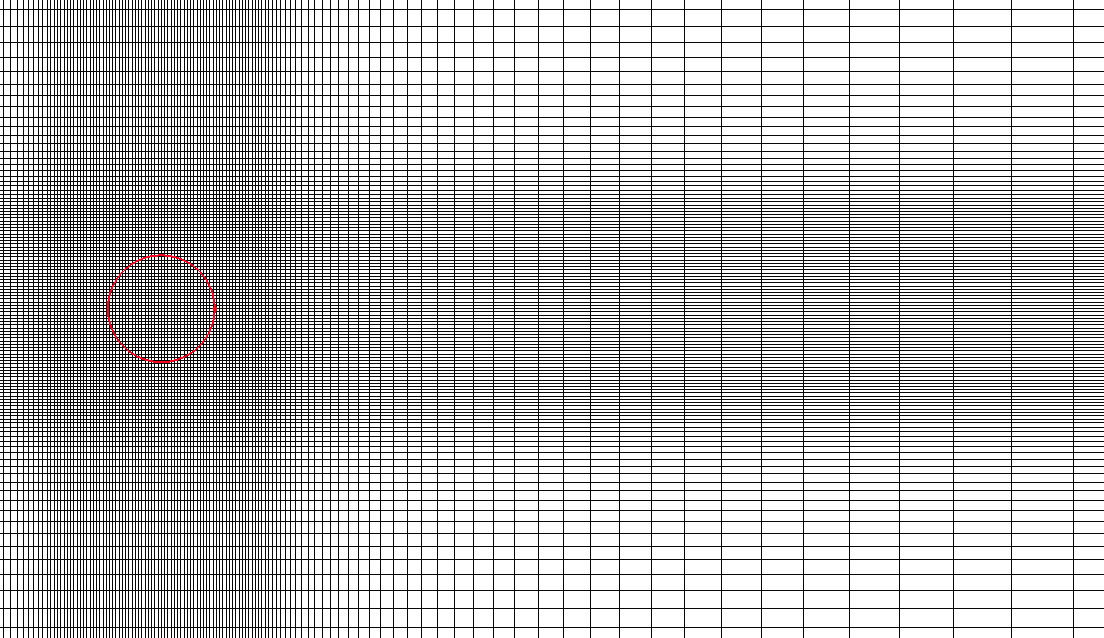}
	\caption{}
	\end{subfigure}
	\begin{subfigure}{.4\textwidth}
		\includegraphics[width=150pt]{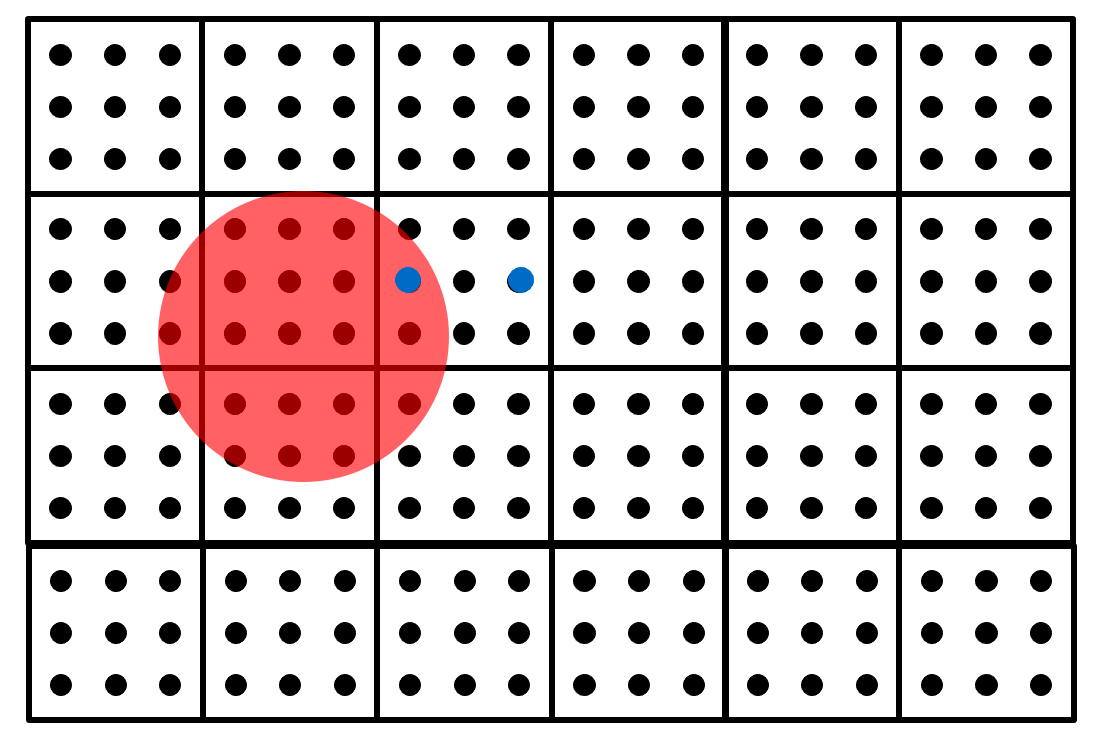}
		\caption{}
	\end{subfigure}
	\centering
	\caption{The mesh for the cylinder and sketch of the probe. a) The locally uniform grid with size $0.03$. b) Sketch of probe positions (this figure is only an illustration and is not scaled with the real mesh).}
	\label{fig:cylinder}
\end{figure*}

This section extends the method to simulate unsteady flow dynamics, exemplified by the flow past a circular cylinder at $Re = 100$ (based on the cylinder diameter $D$ and free stream velocity) and Mach number $0.2$. Again, two probes are selected to monitor the temporal evolution of the velocity. The mesh and probe positions inside and outside are shown in Figure \ref{fig:cylinder}. Again, for the cylinder non-dimensionalized by the diameter $D$, the probes are selected as the closest Gaussian solution points to the particular coordinates $(x,y)/D = (0.36, 0.23)$ for the inside point and $(x,y)/D = (0.75, 0.23)$ for the outside point. 

\begin{figure*}[htbp]
    \begin{subfigure}{.48\textwidth}
		\includegraphics[width=200pt]{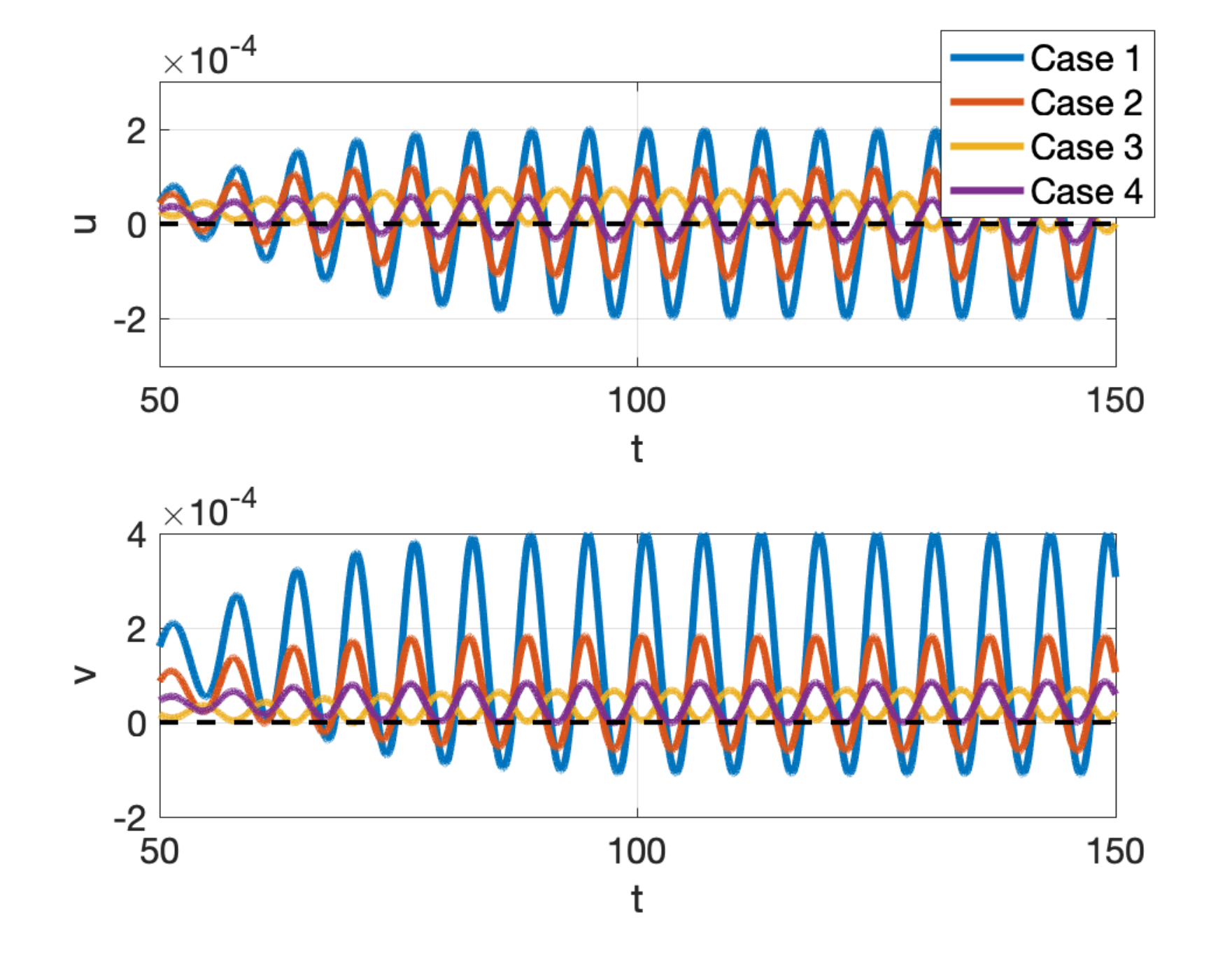}
		\caption{}
	\end{subfigure}
	\begin{subfigure}{.48\textwidth}
		\includegraphics[width=200pt]{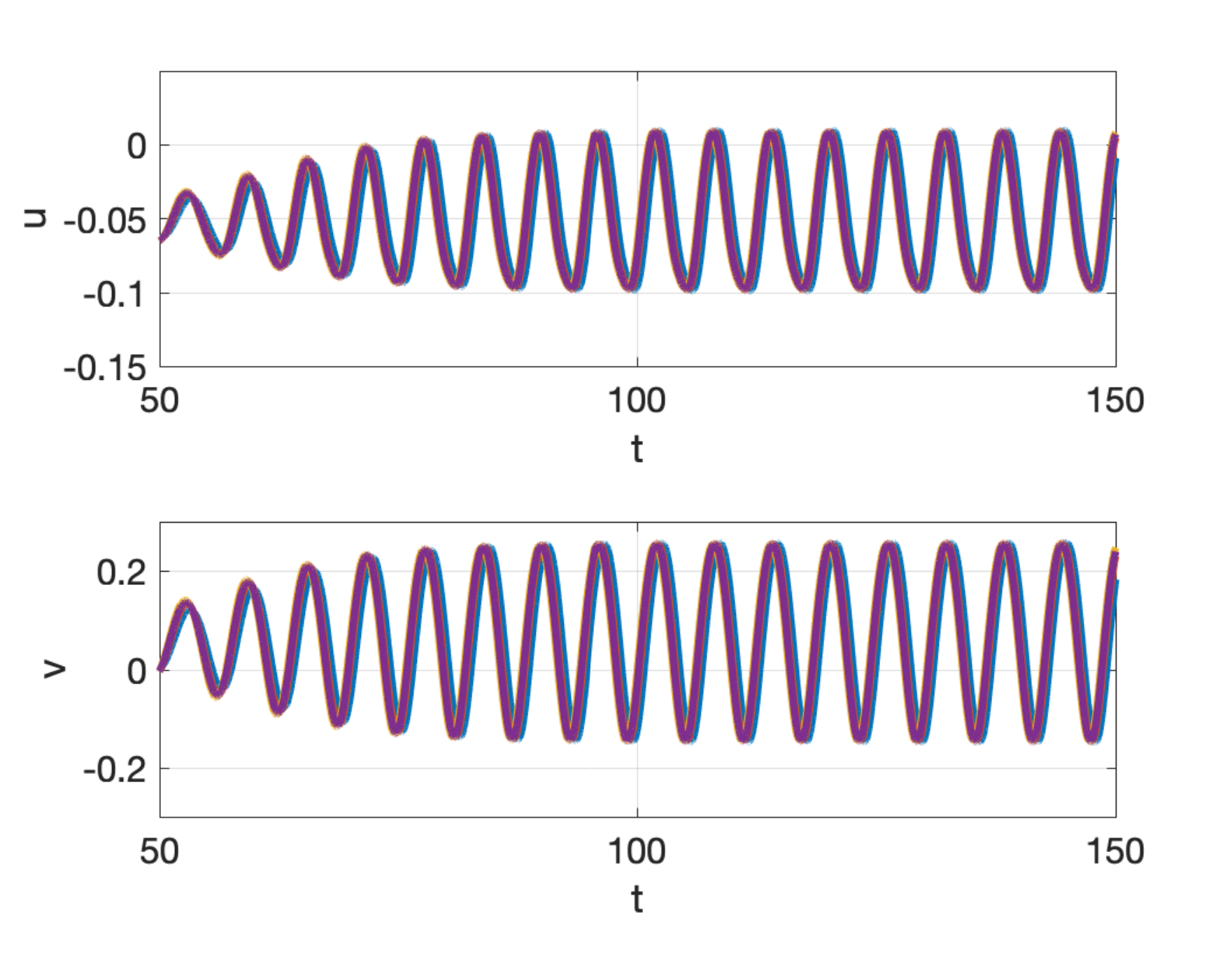}
		\caption{}
	\end{subfigure}
	\centering
	\caption{Evolution of flow velocities at two probe points. a) Inside point. b) Outside point. Case 1: $\IBMparam = 1\times 10^{-3}$. Case 2: $\IBMparam = 5\times 10^{-4}$. Case 3:$\IBMparam = 1\times 10^{-3}, \chi_f = 5\times 10^{3}, \Delta = 100$. Case 4: $\IBMparam = 5\times 10^{-4}, \chi_f = 2\times 10^{3}, \Delta = 100$.}
	\label{fig:cylinder-probe}
\end{figure*}

\begin{figure*}[htbp]
    \begin{subfigure}{.48\textwidth}
		\includegraphics[width=200pt]{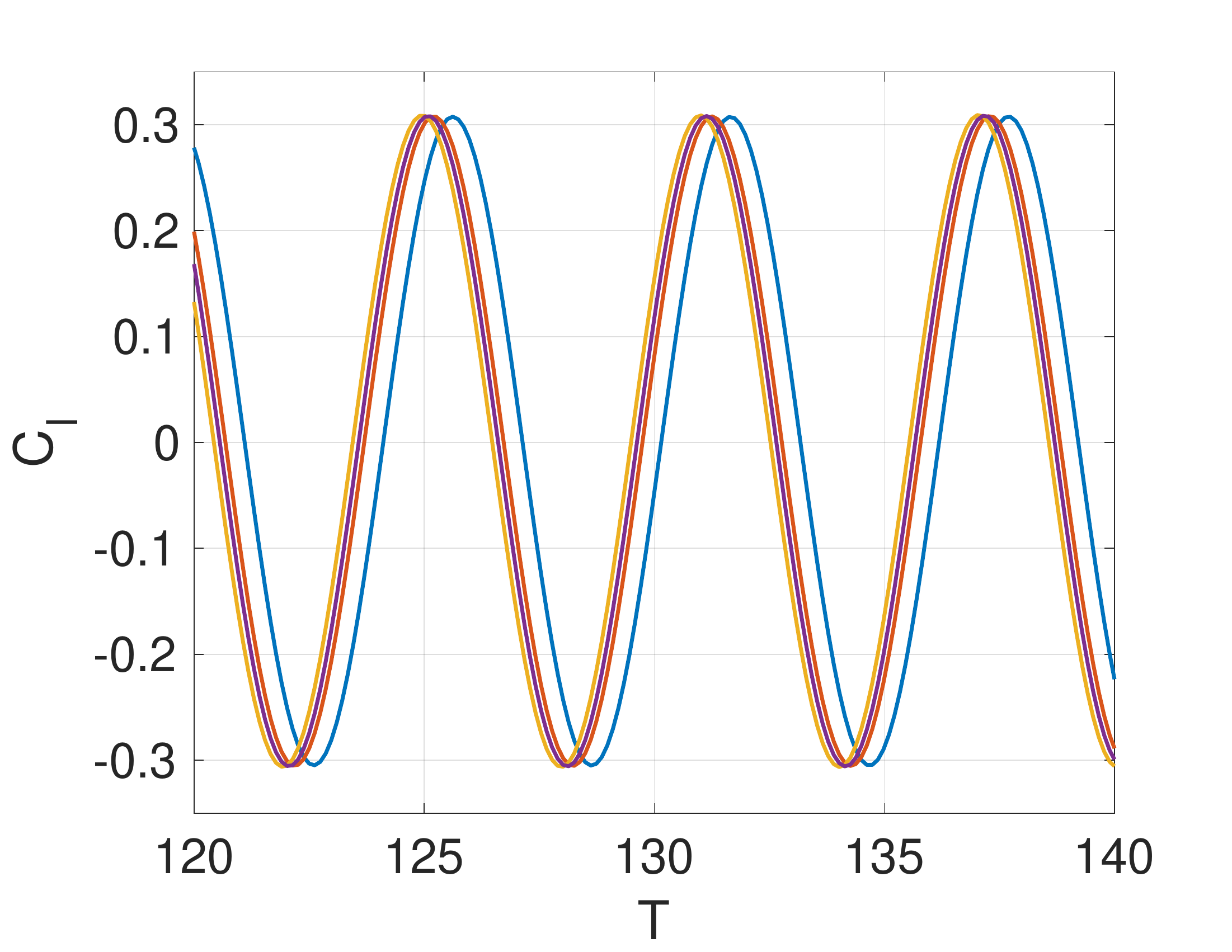}
		\caption{}
	\end{subfigure}
	\begin{subfigure}{.48\textwidth}
		\includegraphics[width=200pt]{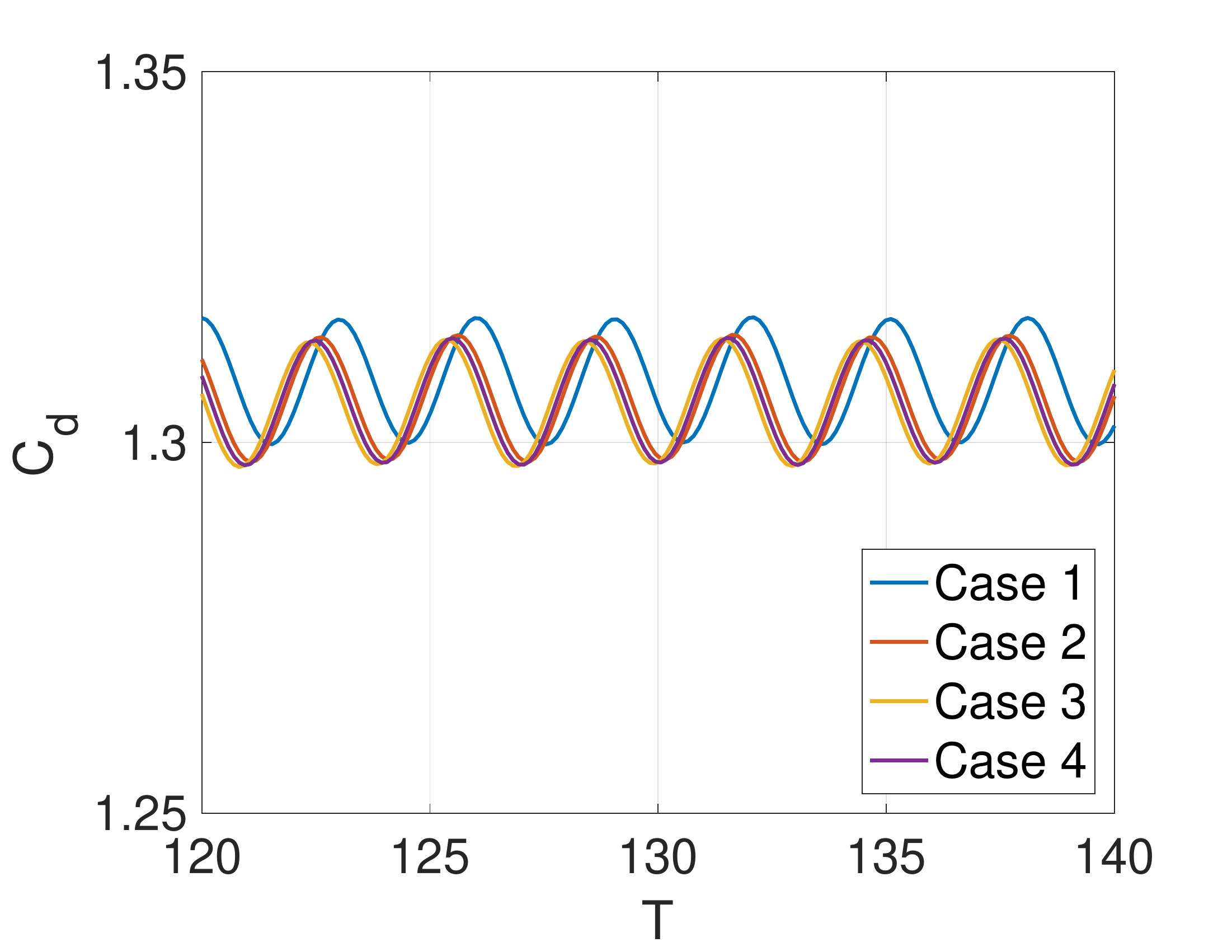}
		\caption{}
	\end{subfigure}
	\centering
	\caption{Evolution of forces. a) Lift coefficient. b) Drag coefficient. Case 1: $\IBMparam = 1\times 10^{-3}$. Case 2: $\IBMparam = 5\times 10^{-4}$. Case 3:$\IBMparam = 1\times 10^{-3}, \chi_f = 5\times 10^{3}, \Delta = 100$. Case 4: $\IBMparam = 5\times 10^{-4}, \chi_f = 2\times 10^{3}, \Delta = 100$.}
	\label{fig:cylinder-force}
\end{figure*}

We generate a uniform grid with size $0.03$ in the domain $x/D \in [-1, 1]$ and $y/D \in [-1, 1]$. The entire computation domain is $x/D \in [-30, 50]$ and $y/D \in [-30, 30]$, which leads to a total number of element $184 \times 178$. The wake region is more refined than in the steady case, since the shedding vortex needs to be captured. Four typical combinations of parameter are tested, and the results are shown in Figure \ref{fig:cylinder-probe} and Figure \ref{fig:cylinder-force}. The order of polynomial is set to $P=2$ and the time step is set to $4 \times 10^{-4}$. The accuracy of imposing the no-slip wall condition is compared in Figure \ref{fig:cylinder-probe}a. The first case is volume penalization with $\IBMparam = 1 \times 10^{-3}$, representing a small penalization source term. Due to the unsteady nature of the problem, a large and periodic temporal oscillation inside the solid is seen (which should still be expected to be zero). If we keep reducing $\IBMparam$, the oscillation inside the body can be damped (Case 2). An alternative is to introduce the SFD method (Case 3) to provide more damping and reduce the oscillations. From the comparison of $u$ in Figure \ref{fig:cylinder-probe}a, the SFD method filters the velocity and reduces the oscillations, but the mean value remains slightly larger when compared to Case 2. This indicates that the volume penalization is also needed to ensure a smaller error in the mean value. Therefore, the optimal performance is given by Case 4, where $\IBMparam = 5 \times 10^{-4}$ and $\chi_f = 1 / \IBMparam$. Then, a small mean value and small oscillations can be achieved, which is more apparent in the evolution of $u$ in Figure \ref{fig:cylinder-probe}a. Looking at Figure \ref{fig:cylinder-probe}b, inaccuracy in imposing the boundary condition can lead to a shift in the phase of the response. In summary, the combined method is advantageous also for unsteady simulations, since the volume penalization helps to correct the mean value and the SFD method helps to reduce the oscillations.

Figure \ref{fig:cylinder-force} compares the temporal evolution of lift and drag coefficients among these four test cases. It can be observed that when the boundary condition is not properly imposed (Case 1), the predicted lift coefficient will have a phase difference, while both the phase and the value differ in the predicted drag coefficient. This indicates that the drag coefficient is more sensitive to boundary conditions. The differences in lift and drag become very small for the other three cases, but the discrepancies in phase and magnitude remain. Therefore, we conclude that the SFD and the combined approach are effective and promising approaches to include immersed boundaries.

\begin{figure*}[htbp]
    \begin{subfigure}{.48\textwidth}
		\includegraphics[width=200pt]{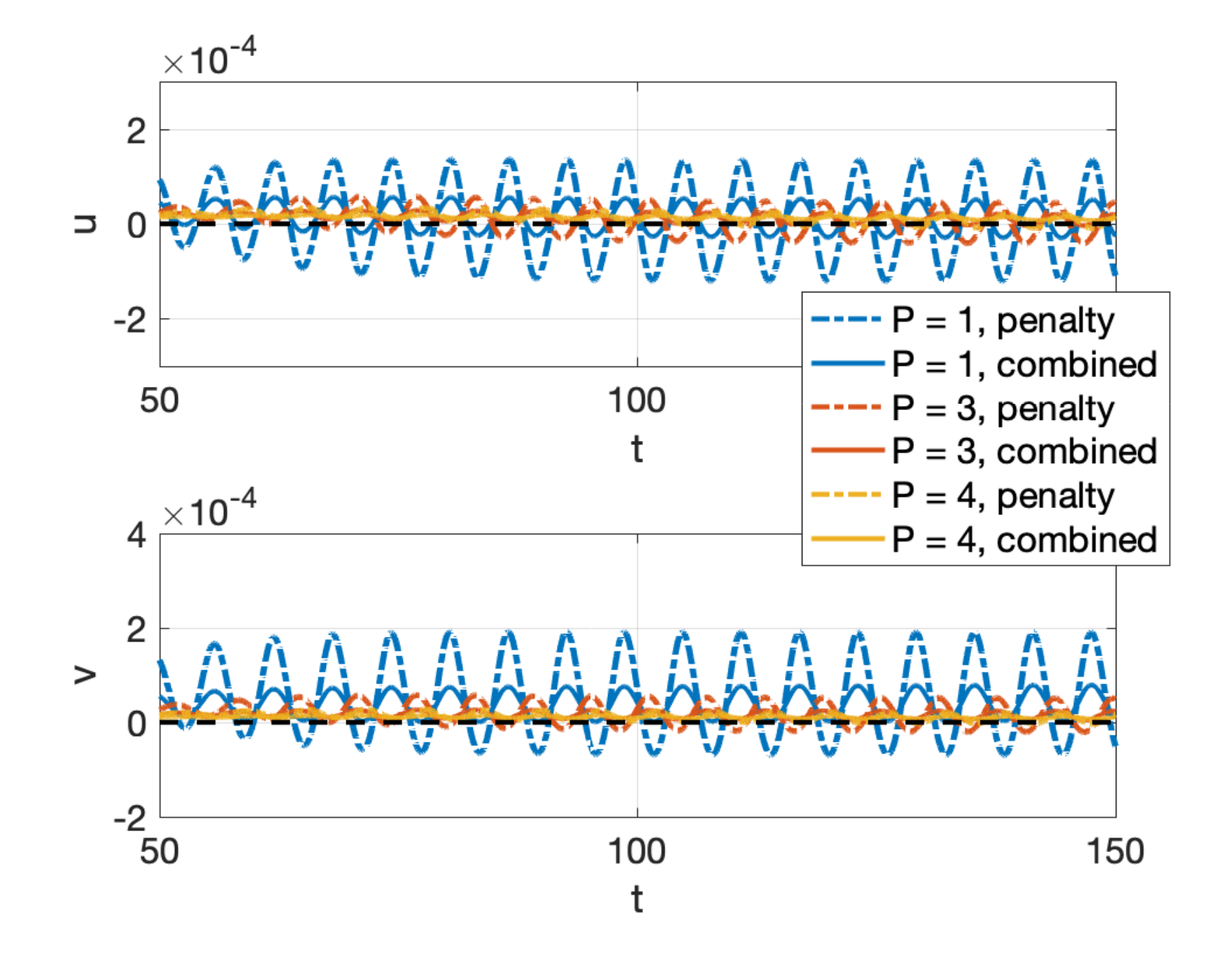}
		\caption{}
	\end{subfigure}
	\begin{subfigure}{.48\textwidth}
		\includegraphics[width=200pt]{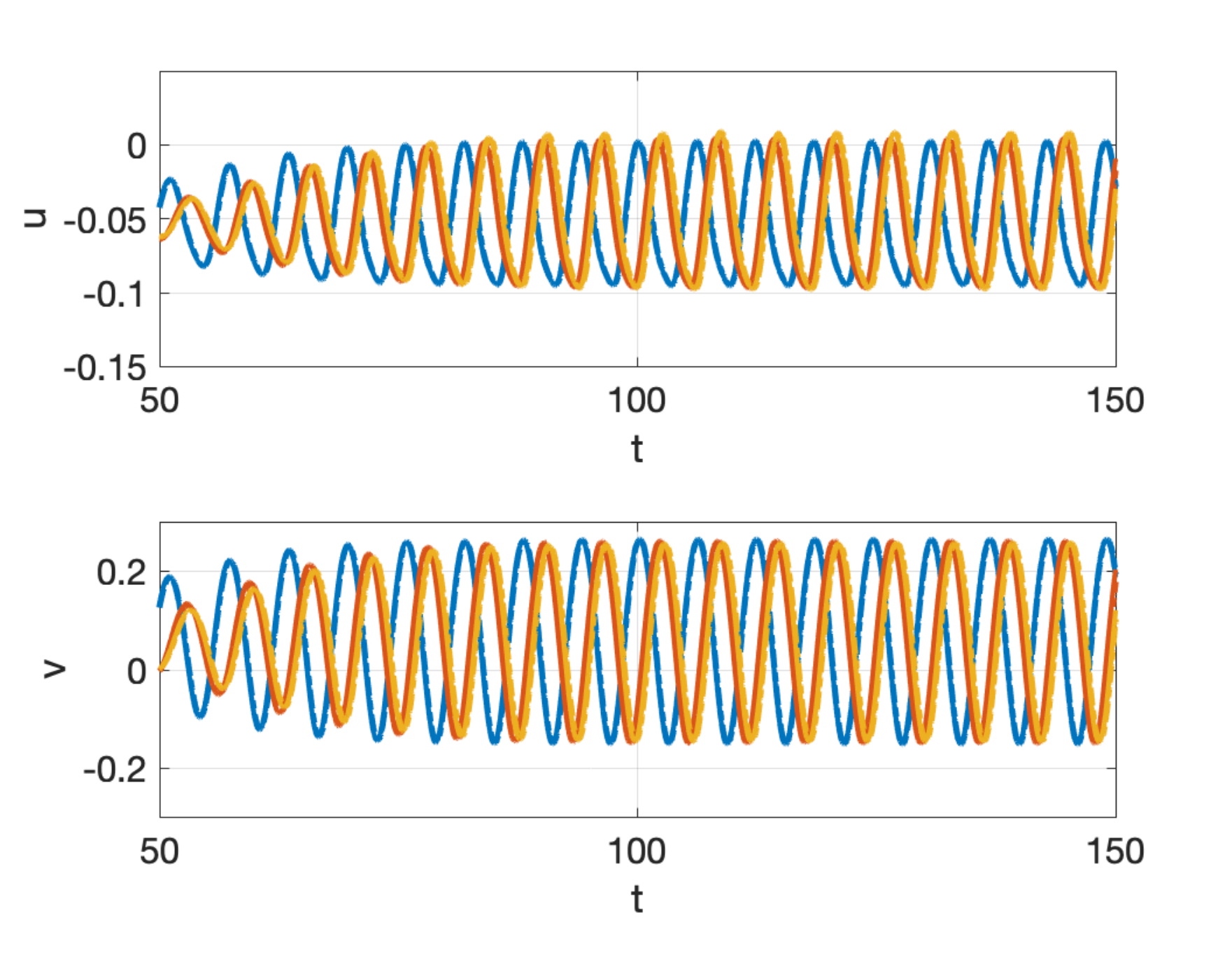}
		\caption{}
	\end{subfigure}
	\centering
	\caption{Evolution of flow velocities at two probe points for different polynomial orders. a) Inside point. b) Outside point. 'penalty' refers to the volume penalization method, while 'combined' refers to the proposed combined volume penalization and the SFD method (with the same penalization parameter). Parameters are set to $\IBMparam = 6 \times 10^{-4}, \chi_f = 1 / \IBMparam, \Delta = 100$ (P=1), $\IBMparam = 1.5 \times 10^{-4}, \chi_f = 1 / \IBMparam, \Delta = 100$ (P=3), $\IBMparam = 5 \times 10^{-5}, \chi_f = 1 / \IBMparam, \Delta = 100$ (P=4).}
	\label{fig:cylinder-probe-highP}
\end{figure*}

\begin{figure*}[htbp]
    \begin{subfigure}{.48\textwidth}
		\includegraphics[width=200pt]{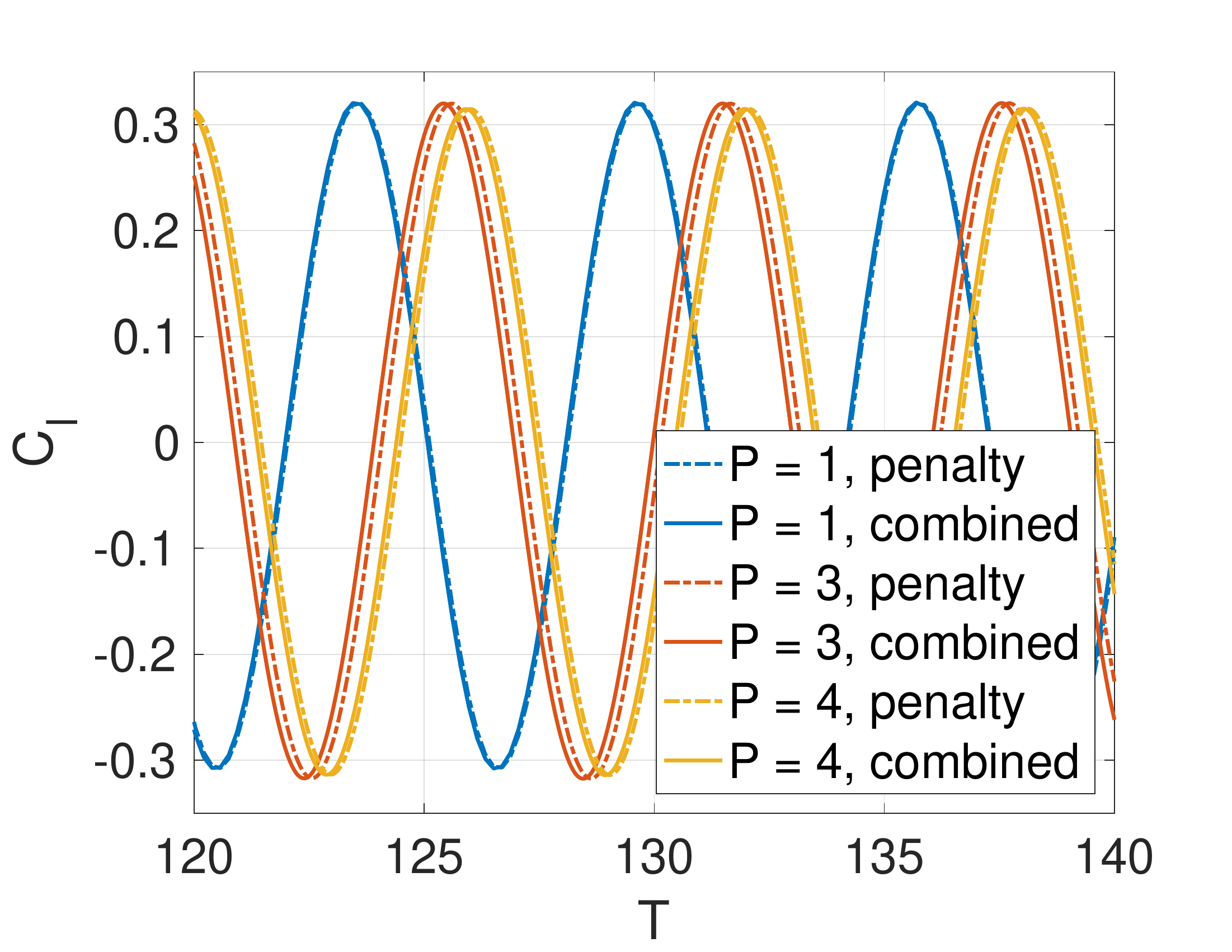}
		\caption{}
	\end{subfigure}
	\begin{subfigure}{.48\textwidth}
		\includegraphics[width=200pt]{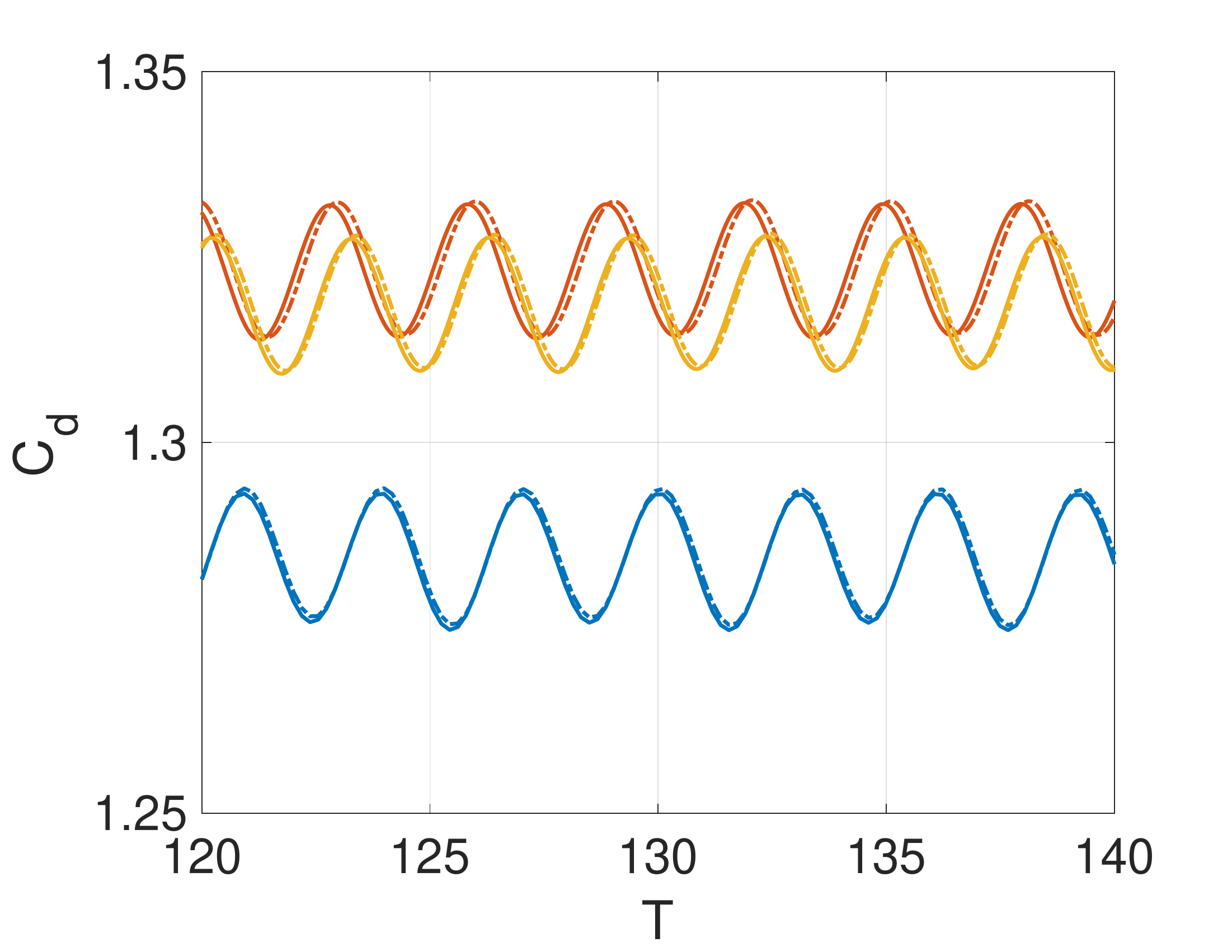}
		\caption{}
	\end{subfigure}
	\centering
	\caption{Evolution of forces for different polynomial orders: a) Lift coefficient. b) Drag coefficient. 
	\textit{penalty} refers to the volume penalization method, while \textit{combined} refers to the proposed combined volume penalization and the SFD method (with the same penalization parameter). Parameters are set to $\IBMparam = 6 \times 10^{-4}, \chi_f = 1 / \IBMparam, \Delta = 100$ (P=1), $\IBMparam = 1.5 \times 10^{-4}, \chi_f = 1 / \IBMparam, \Delta = 100$ (P=3), $\IBMparam = 5 \times 10^{-5}, \chi_f = 1 / \IBMparam, \Delta = 100$ (P=4).}
	\label{fig:cylinder-force-highP}
\end{figure*}

Finally, additional simulations are performed to show the flexibility of the proposed method across different polynomial orders. We use the same grid with polynomial orders $1$, $3$ and $4$. The time step is set to $6 \times 10^{-4}$ ($P = 1$), $1.5 \times 10^{-4}$ ($P = 3$), and $5 \times 10^{-5}$ ($P = 4$). Following the guideline proposed in this paper, the penalization parameter $\IBMparam$ is set to the same value as the time step ($\IBMparam = \Delta t$). In addition, for the combined approach, the control coefficient and filter width are set to $\chi_f = 1 / \IBMparam$ and $\Delta = 100$, respectively. Again, temporal evolution of velocities at two probes and aerodynamic coefficients are compared in Figure \ref{fig:cylinder-probe-highP} and Figure \ref{fig:cylinder-force-highP}, respectively. As shown in Figure \ref{fig:cylinder-probe-highP}a, compared with volume penalization, the proposed combined approach results in reduced oscillation of the velocities inside the solid. When the polynomial order is high (e.g., $P = 3$ or $P = 4$), the combined approach is still able to provide additional damping to better satisfy the no-slip conditions. The improved imposition of boundary condition leads to the phase change of velocities in the flow region, as shown in Figure \ref{fig:cylinder-probe-highP}b. For aerodynamic coefficients, both the differences in the phase and the amplitude are seen, which are more evident in the drag prediction. As shown in Figure \ref{fig:cylinder-force-highP}b, the proposed method can be more efficient at high polynomial orders, since the sensitivity of the drag coefficient to the boundary conditions is more significant. In general, the combination of volume penalization and the SFD approach is very attractive to provide IBM treatment for high-order schemes at different polynomial orders. 

A comparison of aerodynamic coefficients and Strouhal number is shown in Table \ref{table1}. For each polynomial order, the optimal result from the combined approach is shown and is compared against the literature. As the polynomial order increases, we observe converged aerodynamic coefficients and vortex shedding frequencies. The optimal mean drag coefficient is around $1.32$, the optimal fluctuating lift coefficient lies between $0.31$ and $0.32$ while the Strouhal number is about $0.165$. These results are in good agreement with results reported from existing works \citep{braza1986numerical,talley2001experimental,shiels2001flow,gsell2020multigrid}. This indicates that the proposed method can produce accurate prediction for the oscillating aerodynamic loads and frequencies for unsteady flows.

\begin{table}[]
	\caption{Comparison of mean drag coefficient, amplitude of the fluctuating lift coefficient and Strouhal number for flow past a cylinder at Reynolds number $100$. Results of the present study are based on the combined approach.}
	\vspace{20pt}
	\centering
	\begin{tabular}{p{4cm}p{3cm}p{2.5cm}p{2.5cm}p{2.5cm}p{2.5cm}}
		\hline
		Case & $\overline{C_d}$ & $C_{l,max}$ & St \\
		\hline
		Braza et al. \cite{braza1986numerical}  		&  	 1.28 & 0.29 &     0.16          \\
		Talley et al. \cite{talley2001experimental} 		&    1.34 	&   0.33     &    0.16        \\
		Shiels et al. \cite{shiels2001flow} 		&    1.33  	&     0.30    &      0.167	        \\
		Gsell et al. \citep{gsell2020multigrid} & 1.41 &     0.34        &       0.17        \\
		P = 1 		&   1.29 		&     0.31       & 0.164               \\
		P = 2  	&   1.31 		&     0.31        & 0.165             \\
		P = 3  	&   1.32 		&     0.32         & 0.165               \\  
		P = 4  	&   1.32 		&     0.31         & 0.165               \\
		\hline       
	\end{tabular}
	\label{table1}
\end{table}

\section{Conclusions}
This work introduces the SFD method as an alternative and a complementary approach to include immersed boundaries for in high-order solvers. The basic idea of the SFD method is to add a temporal filter inside the solid, to penalize the solution and impose the desired boundary conditions, while filtering out the temporal oscillations of the velocity. Several main conclusions are drawn:

1. A detailed implementation of volume penalization and the SFD method is introduced. The SFD treatment is only applied to the momentum equations of the solid body to impose no-slip boundary conditions for the velocity, and the encapsulated formulation \citep{jordi2014encapsulated} enables very little modification of the underlying time-stepping code.

2. Eigensolution analyses are performed in the proposed scheme, where semi-discrete analysis helps to understand this scheme and provides guidance on the selection of filter width of the SFD method. The fully-discrete analysis indicates that the proposed combined approach ensures good stability, where the classic rule $\IBMparam = \Delta t$ applies. These findings have been analyzed by simulating the one-dimensional advection equation.

3. Steady and unsteady flow simulations based on the Navier-Stokes equations have been performed with different combinations of parameters. The combined approach shows advantages in providing more accurate immersed boundaries, resulting in minimal unphysical oscillations and improved accuracy in the aerodynamic forces.

4. The combined approach (volume penalization and the SFD method) provides more accurate results for unsteady cases and for all polynomial orders tested ($1 \leq P \leq 4$).

\section*{Acknowledgments}
This project has received funding from the European Union’s Horizon 2020 research and innovation programme under the Marie Skłodowska-Curie grant agreement (MSCA ITN-EID-GA ASIMIA No 813605).

\appendix
\section{discretisation and matrix formulation of one-dimensional advection equation}
\label{append2}
When considering high-order discretisations, we have the following matrix form of each element for standard advection equation \citep{moura2015linear,kou2021VonNeumann}:

\begin{equation}
    \frac{d\boldsymbol{u}_n}{dt} = \boldsymbol{L} \boldsymbol{u}_{n-1} +  \boldsymbol{C} \boldsymbol{u}_{n} +  \boldsymbol{R} \boldsymbol{u}_{n+1},
\end{equation}
where the left, middle and right operators for the current element $\boldsymbol{u}_n$ are defined as $\boldsymbol{L}$, $\boldsymbol{C}$, and $\boldsymbol{R}$, respectively. The solution is represented by the polynomial representation in terms of a nodal basis ($\Tilde{u}(r_{i}), i=0,1,...,P$) with order $P$ defined at $P+1$ solution points $\Tilde{u}_n(r) = \sum_{i=0}^{P}l_{i}(r)\Tilde{u}_n(r_{i}) = \boldsymbol{l}^T \Tilde{\boldsymbol{u}}_n$, where $\boldsymbol{l}$ is a Lagrange interpolation operator and $l_{i}(r)$ denotes the Lagrange polynomial defined at a solution point $r_{i}$. These matrix can be derived based on flux reconstruction: 

\begin{equation}
    \boldsymbol{L} = - \frac{1}{h} \boldsymbol{g}^L_r \left [ (c+\lambda \left | c \right |)\boldsymbol{l}^T(1) \right ]
\end{equation}

\begin{equation}
   \boldsymbol{C} = - \frac{2}{h} \left [ c \boldsymbol{D} - \frac{1}{2} \boldsymbol{g}^L_r  (c+\lambda \left | c \right | ) \boldsymbol{l}^T(-1) + \frac{1}{2} \boldsymbol{g}^R_r ( \lambda \left | c \right | - c ) \boldsymbol{l}^T(1)  \right]
\end{equation}

\begin{equation}
    \boldsymbol{R} = - \frac{1}{h} \boldsymbol{g}^R_r \left [ (c-\lambda \left | c \right | )\boldsymbol{l}^T(-1) \right ]
\end{equation}
where the gradient vector of correction function at all solution points is defined as 
\begin{equation}
    \boldsymbol{g}^{dir}_r = \left [ \frac{\partial g^{dir}}{\partial r}(r_{0}), \frac{\partial g^{dir}}{\partial r}(r_{1}),..., \frac{\partial g^{dir}}{\partial r}(r_{P})  \right]^T \ , \ dir = L (left) \ or \ R (right).
\end{equation}
where $\lambda$ is the upwind parameter of the Lax–Friedrichs type numerical flux. The local gradient operator is defined as $\boldsymbol{D}$ with $\boldsymbol{D}[i,j]=\frac{ \partial l_{j}(r_i) }{\partial r}$. For standard linear advection equation with periodic boundary conditions, we have the following matrix form
\begin{equation}
    \frac{d}{dt} \begin{pmatrix} \boldsymbol{\velx}_{1}\\  \boldsymbol{\velx}_{2}\\ \boldsymbol{\velx}_{3}\\ \vdots\\ \boldsymbol{\velx}_{N} \end{pmatrix} = \begin{pmatrix} \boldsymbol{C}-\frac{\chi}{\eta}\boldsymbol{I}
    & \boldsymbol{R} & \boldsymbol{0} & \cdots & \text{exp}[-2ikT] \boldsymbol{L}\\ \boldsymbol{L} & \boldsymbol {C}-\frac{\chi}{\eta} \boldsymbol{I} & \boldsymbol{R} & \cdots & \boldsymbol{0} &\\ \boldsymbol{0} & \boldsymbol{L} &\boldsymbol{C}-\frac{\chi}{\eta}\boldsymbol{I} & \cdots & \boldsymbol{0}\\ \vdots & \vdots & \vdots & \ddots &\vdots\\ \text{exp}[2ikT] \boldsymbol{R} & \boldsymbol{0} & \boldsymbol{0} & \cdots &\boldsymbol{C}-\frac{\chi}{\eta}\boldsymbol{I} \end{pmatrix}
    \begin{pmatrix} \boldsymbol{\velx}_{1}\\  \boldsymbol{\velx}_{2}\\ \boldsymbol{\velx}_{3}\\ \vdots\\ \boldsymbol{\velx}_{N} \end{pmatrix}.
\end{equation}

In addition, considering a no-slip wall in the middle imposed by SFD and volume penalization, we have the matrix form for the advection equation

\begin{equation}
\label{eq:global}
    \frac{d}{dt} \begin{pmatrix} \boldsymbol{\velx}_{1}\\  \boldsymbol{\velx}_{2} \\ \vdots\\ \boldsymbol{\velx}_{IBM} \\ \vdots \\ \boldsymbol{\velx}_{N} \\ \bar{\boldsymbol{q}}\end{pmatrix} = \begin{pmatrix} \boldsymbol{C} & \boldsymbol{R} & \boldsymbol{0} & \cdots & \cdots & \cdots & \text{exp}[-2ikT] \boldsymbol{L} & \boldsymbol{0}\\ \boldsymbol{L} & \boldsymbol{C} &\boldsymbol{R} & \cdots & \cdots & \cdots & \boldsymbol{0} & \boldsymbol{0} \\ \vdots & \vdots & \vdots & \ddots & \ddots & \ddots & \vdots & \vdots\\ \vdots & \vdots & \vdots & \boldsymbol{L} & \boldsymbol{C}-\frac{1}{\eta}\boldsymbol{I} - \chi_f \boldsymbol{I} & \boldsymbol{R}& \vdots & \chi_f \boldsymbol{I}\\ \vdots & \vdots & \vdots & \ddots & \ddots & \ddots & \vdots & \vdots\\
    \text{exp}[2ikT] \boldsymbol{R} & \boldsymbol{0} & \boldsymbol{0} & \cdots & \cdots & \cdots &\boldsymbol{C} & \boldsymbol{0}\\ \boldsymbol{0} & \boldsymbol{0} & \boldsymbol{0} & \cdots & \boldsymbol{I} / \Delta & \cdots & \boldsymbol{0} & -\boldsymbol{I} / \Delta \end{pmatrix}
    \begin{pmatrix} \boldsymbol{\velx}_{1}\\  \boldsymbol{\velx}_{2} \\ \vdots\\ \boldsymbol{\velx}_{IBM} \\ \vdots \\ \boldsymbol{\velx}_{N} \\ \bar{\boldsymbol{q}} \end{pmatrix}.
\end{equation}

Note that the extension to the advection equation with only volume penalization or only SFD is straightforward through removing the corresponding terms in the above equation.

\section{Navier-Stokes Equations}
\label{append1}
The governing equations for a compressible viscous fluid are written as
\begin{equation}
    \frac{\partial \bm{U}}{\partial t} + \mathbf{\nabla} \cdot \bm{F} =  \frac{\partial \bm{U}}{\partial t} + \frac{\partial \bm{F}_x}{\partial x} + \frac{\partial \bm{F}_y}{\partial y} + \frac{\partial \bm{F}_z}{\partial z} = 0\,,
\end{equation}
where $\bm{U}$ denotes the vector of conserved variables $\bm{U} = (\rho , \rho u, \rho v, \rho w, E)^T$. $\rho$ is the density, $u\,, v\, \text{and}\, w$ are the velocity components and $E$ is the total energy. The equations are closed by the ideal gas equation-of-state:
\begin{equation}
    E = \frac{P}{\gamma - 1}+\frac{1}{2} \rho (u^2+v^2+w^2)\,,
\end{equation}
where $P$ is the static pressure and $\gamma$ is the ratio of specific heats. The flux vectors $\bm{F}_{x}$, $\bm{F}_{y}$, $\bm{F}_{z}$ contain the inviscid and viscous fluxes and are written as
\begin{equation}
\label{eq:flux-definition}
    \bm{F}_{x} = \begin{pmatrix} \rho u\\  \rho u^2+P\\ \rho u v\\ \rho u w\\ u (E+P) \end{pmatrix}
    -
    \begin{pmatrix} 0\\ \tau_{{x} {x}} \\ \tau_{{x} {y}}\\ \tau_{{x} {z}}\\ u\tau_{{x} {x}} + v \tau_{{x}{y}} + w \tau_{{x}{z}} + q_{x} \end{pmatrix} = \bm{F}_{x,inv} + \bm{F}_{x,vsc} 
\end{equation}

\begin{equation}
    \bm{F}_{y} = \begin{pmatrix} \rho v\\  \rho u v\\ \rho v^2+P\\ \rho v w\\ v (E+P) \end{pmatrix}
    -
    \begin{pmatrix} 0\\ \tau_{{y}{x}} \\ \tau_{{y}{y}}\\ \tau_{{y}{z}}\\ u \tau_{{y}{x}} + v \tau_{{y}{y}} + w \tau_{{y}{z}} + q_{y} \end{pmatrix} =  \bm{F}_{y,inv} + \bm{F}_{y,vsc} 
\end{equation}

\begin{equation}
   \bm{F}_{z} = \begin{pmatrix} \rho w\\  \rho u w\\ \rho v w\\ \rho w ^ 2 + P \\ w (E+P) \end{pmatrix}
    -
    \begin{pmatrix} 0\\ \tau_{{z}{x}} \\ \tau_{{z}{y}}\\ \tau_{{z}{z}}\\ u \tau_{{z}{x}} + v \tau_{{z}{y}} + w \tau_{{z}{z}} + q_{z} \end{pmatrix} = \bm{F}_{z,inv} + \bm{F}_{z,vsc} .
\end{equation}
In these equations, $\tau_{{i}{j}} = \mu (\frac{\partial v_{i}}{\partial x_{j}} + \frac{\partial v_{j}}{\partial x_{i}} - \frac{2}{3} \delta _{ij}\frac{\partial v_{k}}{\partial x_{k}})$ is the viscous stress tensor with $\mu$ denoting the dynamic viscosity. The heat flux vector $ \nabla \bm{q}_h$ is given by
\begin{equation}
    \frac{\partial \bm{q}_h}{\partial x_{i}} = \lambda \frac{\partial T}{\partial x_{i}},
\end{equation}
where $\lambda$ is the thermal conductivity and $T$ is the static temperature. The equations are solved in non-dimensional form with the introduction of the Prandtl number $\Pr = \mu \frac{C_{p}}{\lambda}$, the Reynolds number $ Re = \rho_{\text{ref}} V_{ref}L_{ref}/ \mu_{ref}$ and the Mach number $ M = V_{ref}/\sqrt{\gamma R_{gas}T_{ref}}$, with $C_{p}$ being the specific heat capacity at constant pressure and $R_{gas}$ being the gas constant. Finally, $V_{ref}, L_{ref}, T_{ref}$ are reference velocity, length and temperature, respectively. The discretisation of these equations with the flux reconstruction method is described next.

\section{Flux reconstruction for general conservation laws}
\label{append3}
Without loss of generality, we consider the following second-order hyperbolic conservation law:
\begin{equation}
\begin{aligned}
	\frac{\partial \bm{U}}{\partial t} + \mathbf{\nabla} \cdot \bm{F} (\bm{Q},\bm{U}) &= \bm{S}, \\
	\bm{Q} - \mathbf{\nabla}\bm{U} &= 0,
\end{aligned}
\end{equation}
where $\bm{Q}$, $\bm{F}$ and $\bm{S}$ refer to the gradient, the flux and the source term respectively, and $\bm{F}$ includes both inviscid (advection) and viscous (diffusion) terms in all directions. FR is a popular high-order scheme to discretize the conservation law in space, which can recover existing high-order schemes including nodal DG and SD, as shown in \citep{huynh2007FR,vincent2011ESFR}.

Like other high-order methods, the preliminaries of FR includes the discretisation of the computational domain, the definition of solution and flux points, as well as the space transformation. The domain is firstly discretized by non overlapping elements. For each element, $Np$ internal degrees of freedom and $Nf$ points on the element interface are defined. The solution $\bm{U}$ and flux $\bm{F}$ are represented by a polynomial of degree $P$ and $P+1$. More details are given in \citep{huynh2007FR,vincent2011ESFR,williams2013energy,castonguay2013energy}. For the general conservation law, the FR procedure generally includes seven stages, where more details can be found in \citep{williams2013energy,castonguay2013energy}. In the first three stages, the gradient $\boldsymbol{Q}$ is obtained in order to form the viscous flux. The polynomial interpolation of solution $\boldsymbol{U}$ is firstly obtained
\begin{equation}
\bm{U}^{\delta\textrm{D}} (\bm{\xi}) = \sum_{i=1}^{N_p} \bm{U}^{\delta\textrm{D}}_{i} {\mathcal{I}^{P}_{i}}(\bm{\xi}),
\end{equation}
where $\mathcal{I}^{P}_{i}$ refers to the nodal basis function defined at each solution point with polynomials of degree $P$, and $\bm{U}_{i}^{\delta\textrm{D}} $ is the solution at the $i$th solution point. ${\delta\textrm{D}}$ indicates that the solution is discrete and discontinuous within the element. This formulation allows interpolating the solution to the flux points at the interface $\bm{L}_{f}\bm{U}^{\delta\textrm{D}}$, where $\bm{L}_{f}$ is the interpolation operator to obtain values at the flux point. In the second stage, the common solution values at the interfaces and flux points $\bm{U}^{Comm}(\bm{\xi}_{flux})$ are obtained from the approximate discontinuous solution on the left and right of the interface. These common solution values can be solved by different approaches such as the Central Flux \citep{hesthaven2007nodal}, Local Discontinuous Galerkin (LDG) \citep{cockburn1998LDG}, Internal Penalty \citep{arnold1982interior}, etc. In the third stage, the correction function, which is the key component in FR, is introduce to leverage the discontinuous values to be continuous values across the element. This results in the following continuous solution fluxes:
\begin{equation}
\bm{U}^{\delta\textrm{C}}(\bm{\xi}) = \bm{U}^{\delta\textrm{D}}(\bm{\xi}) + (\bm{U}^{Comm}(\bm{\xi}_{flux}) - \bm{L}_{f}\bm{U}^{\delta\textrm{D}}) {\bm{h}}(\bm{\xi}),
\end{equation}
where $\bm{h}$ serves as a “lifting” operator \cite{williams2013energy} that will transfer the correction values to the quadrature points within the element. Consequently, the second term in these equations are defined as the correction term which is the main ingredient of FR (denoted with superscript Corr, $\bm{U}^{\textrm{Corr}} = \bm{U}^{Comm}(\bm{\xi}_{flux}) - \bm{L}_{f}\bm{U}^{\delta\textrm{D}}$). From this expression, the gradient in the reference space can be obtained as follow:
\begin{equation}
\boldsymbol{Q}_x = \Tilde{\mathbf{\nabla}}_x \bm{U} = \sum_{i=1}^{N_p} \bm{U}_{i}^{\delta\textrm{D}} \Tilde{\mathbf{\nabla}}_x {\mathcal{I}^{P}_{i}}(\bm{\xi}) + \sum_{f=1}^{N_{face}} \sum_{j=1}^{N_f} \Tilde{\mathbf{\nabla}}_x \mathcal{C}^{P+1}_{f,j}(\bm{\xi}) \cdot \bm{U}_{f,j}^{\textrm{Corr}},
\end{equation}
where $\Tilde{\mathbf{\nabla}}_x$ is the discrete gradient in the reference space in $x$ direction. The function $\mathcal{C}^{P+1}_{f,j}$ is the component of correction function $\boldsymbol{h}$, which is of polynomial order $ P + 1$. Note that the resulting gradient and flux should be transformed to the physical space.

The divergence of flux is solved in the remaining four stages, which follows similar procedures as computing the gradient. After substituting the solution $\boldsymbol{U}$ and gradient $\boldsymbol{Q}$ into the flux $\boldsymbol{F}$ to obtain the flux at each solution point, the flux can also be represented as the same polynomial of degree $P$ defined in the nodal basis:
\begin{equation}
\bm{F}^{\delta\textrm{D}} (\bm{\xi}) = \sum_{i=1}^{N_p} \bm{F}^{\delta\textrm{D}}_{i} {\mathcal{I}^{P}_{i}}(\bm{\xi}),
\end{equation}

Therefore, the interpolated discontinuous flux at the flux points are obtained from the this formulation in the fourth stage, through the interpolation matrix $\boldsymbol{L}_f$. In the fifth stage, the numerical interaction flux $\boldsymbol{F}^{Int}$ at flux points along the edges of the element is obtained. At each flux point, this continuous interface flux is constructed from the multiply defined discontinuous fluxes from both sides. This flux, including both inviscid (advection) and viscous (diffusion) parts, is approximated by the Riemann solvers. For example, a Lax–Friedrichs flux is used for the advection term of the advection-diffusion equation, while a Roe \citep{roe1981approximate} or a Rusanov \citep{rusanov1961calculation} type approximate Riemann solver can be considered for the inviscid flux of the Navier-Stokes equations. The viscous flux is then obtained using one of the aforementioned approach for the common flux. The sixth stage involves forming the transformed flux correction based on the correction function $\boldsymbol{h}$ which extends the values defined at the interface flux points to the quadrature points within the element and corrects the values at the flux points to be the interaction flux. This leads to the following formulation for the continuous flux
\begin{equation}
\bm{F}^{\delta\textrm{C}}(\bm{\xi}) = \bm{F}^{\delta\textrm{D}}(\bm{\xi}) + (\bm{F}^{Int}(\bm{\xi}_{flux}) - \bm{L}_{f}\bm{F}^{\delta\textrm{D}}) \cdot \vec{n} {\bm{h}}(\bm{\xi}),
\end{equation}
where the second term represents the correction flux along each direction, and the correction function $\boldsymbol{h}$ is approximated by a polynomial of degree $P+1$. In the final stage, the divergence of this continuous flux is obtained as

\begin{equation}
\Tilde{\mathbf{\nabla}} \cdot \bm{F} =  \sum_{k=1}^{\mathcal{D}} \sum_{i=1}^{N_p} \Tilde{\mathbf{\nabla}}_{k} {\mathcal{I}^{P}_{i}}(\bm{\xi}) \cdot \bm{F} _{i,k} ^ {\delta\textrm{D}} +  \sum_{k=1}^{\mathcal{D}} \sum_{f=1}^{N_{face}} \sum_{j=1}^{N_f} \Tilde{\mathbf{\nabla}}_{k} {\mathcal{C}^{P+1}_{f,j}}(\bm{\xi}) \cdot \bm{F}_{f,j,k}^{\textrm{Corr}},
\end{equation}
where $\mathcal{D}$ refers to the space dimension, $N_{face}$ refers to number of faces for each element, and $h$ is the uniform mesh size. $\Tilde{\mathbf{\nabla}}$ is the discrete gradient in the reference space. The correction flux is $\bm{F}^{\textrm{Corr}}$ where $\bm{F}^{\textrm{Corr}} = (\bm{F}^{Int}(\bm{\xi}_{flux}) - \bm{L}_{f}\bm{F}^{\delta\textrm{D}}) \cdot \vec{n}$. The function $\mathcal{C}^{P+1}_{f,j}$ is the component of correction function $\boldsymbol{h}$, which is of polynomial order $P + 1$ and can be different from that used for computing the gradient. Finally, the resulting gradient and flux should be transformed to the physical space. In the complete FR workflow, six factors \citep{williams2013energy} will influence the performance of FR, including the location of solution and flux points, the Riemann solvers used for computing the common solution values and the interaction fluxes, and the form of the correction functions for solution and flux values. In the present study, we choose the Gaussian quadrature points as the solution and the flux points. Riemann solvers based on the Rusanov and the LDG approaches are selected for the inviscid and viscous fluxes, respectively. The boundary conditions, in particular the characteristic boundary condition for the far field, are imposed in the computation of common and interaction fluxes \citep{mengaldo2014guide}. For time integration, explicit time integration scheme using the third-order TVD Runge-Kutta method \citep{gottlieb1998total} or the low-storage five-stage fourth-order explicit Runge-Kutta method (LSERK) \citep{carpenter1994fourth,hesthaven2007nodal} method is considered. Note that both methods produce identical predictions but the latter allows a relatively larger time step.

\bibliographystyle{model1-num-names}
\bibliography{refs}

\end{document}